\documentclass{amsart}
\usepackage{amssymb, latexsym}

\def\bi{\begin{itemize}}
\def\bs{\begin{split}}
\def\es{\end{split}}
\def\ba{\begin{align}}
\def\bas{\begin{align*}}
\def\ea{\end{align}}
\def\eas{\end{align*}}
\newcommand{\util}{\tilde{u}}
\def\Im{{\hbox{Im}}}

\def\Re{{\hbox{Re}}}

\def\C{{\mathbb C}} 
\def\R{{{\mathbb R}}}
\def\Z{{{\mathbb Z}}}
\def\N{{{\mathbb N}}}

\def\O{{{\mathcal O}}}

\def\emph#1{{\it #1}}
\def\textbf#1{{\bf #1}}

\newcommand{\lqlr}{{L^q_t L^r_x}}
\newcommand{\ir}{{I \times \R^n}}
\newcommand{\rr}{{\R \times \R^n}}
\newcommand{\lli}{{\lqlr (\ir)}}
\newcommand{\llr}{{\lqlr (\rr)}}
\newcommand{\sz}{{\dot S^0}}
\newcommand{\szi}{{\sz (\ir)}}

\newcommand{\so}{{\dot S^1}}
\newcommand{\soi}{{\so (\ir)}}
\newcommand{\sor}{{\so (\rr)}}

\newcommand{\ho}{{\dot H^1}}

\newcommand{\is}{{I_* \times \R^n}}

\newcommand{\sois}{{\so (\is)}}
\newcommand{\szis}{{\sz (\is)}}

\newcommand{\eps}{{\varepsilon}}
\newcommand{\uhi}{{u_{hi}}}
\newcommand{\ulo}{{u_{lo}}}
\newcommand{\uhip}{{u_{hi'}}}
\newcommand{\ecrit}{{E_{crit}}}
\newcommand{\pl}{{P_{lo}}}
\newcommand{\ph}{{P_{hi}}}
\newcommand{\pmed}{{P_{med}}}

\theoremstyle{plain}
\newtheorem{theorem}{Theorem}
\newtheorem{definition}[theorem]{Definition}
\newtheorem{remark}[theorem]{Remark}
\newtheorem{proposition}[theorem]{Proposition}
\newtheorem{lemma}[theorem]{Lemma}
\newtheorem{corollary}[theorem]{Corollary}

\numberwithin{equation}{section} \numberwithin{theorem}{section}

\begin{document}

\title[Global well-posedness for the energy-critical NLS]%
{The defocusing energy-critical nonlinear Schr\"odinger equation in higher dimensions}
\author{Monica Visan}
\address{University of California, Los Angeles}
\email{mvisan@math.ucla.edu}

\vspace{-0.3in}
\begin{abstract}
We obtain global well-posedness, scattering, and global $L^{\frac{2(n+2)}{n-2}}_{t,x}$ spacetime bounds for
energy-space solutions to the energy-critical nonlinear Schr\"odinger equation in $\R_t\times \R^n_x$, $n\geq 5$.
\end{abstract}

\maketitle

\section{Introduction}

We study the initial value problem for the defocusing energy-critical nonlinear Schr\"odinger equation in
$\R\times\R^n$, $n\geq 5$,
\begin{equation}\label{schrodinger equation}
\begin{cases}
i u_t +\Delta u = |u|^{\frac{4}{n-2}} u\\
u(0,x) = u_0(x)
\end{cases}
\end{equation}
where $u(t,x)$ is a complex-valued function in spacetime $\R_t\times \R^n_x$.

This equation has Hamiltonian
\begin{equation}\label{energy}
E(u(t))=\int_{\R^n} \Bigl(\tfrac{1}{2}|\nabla u (t,x)|^2+\tfrac{n-2}{2n}|u(t,x)|^{\frac{2n}{n-2}}\Bigr)dx.
\end{equation}
Since \eqref{energy} is preserved by the flow corresponding to \eqref{schrodinger equation}, we shall refer to it as the
\emph{energy} and often write $E(u)$ for $E(u(t))$.

We are interested in this particular nonlinearity because it is critical with respect to the energy norm.  That is,
the scaling $u \mapsto u^\lambda$ where
\begin{equation}\label{scaling}
u^\lambda (t,x) := \lambda^{-\frac{n-2}{2}} u\bigl(\lambda^{-2}t,\lambda^{-1}x\bigr)
\end{equation}
maps a solution of \eqref{schrodinger equation} to another solution of \eqref{schrodinger equation}, and $u$ and
$u^\lambda$ have the same energy.

A second conserved quantity we will occasionally rely on is the mass $\|u(t)\|_{L^2_x(\R^n)}^2$.  However, since the
equation is $L_x^2$-supercritical with respect to the scaling \eqref{scaling}, we do not have bounds on the mass that are
uniform across frequencies. Indeed, the low frequencies may simultaneously have small energy and large mass.

It is known that if the initial data $u_0$ has finite energy, then \eqref{schrodinger equation} is locally well-posed
(see, for instance \cite{cwI}, \cite{TV}). That is, there exists a unique local-in-time solution that lies in
$C_t^0 \ho_x \cap L^{\frac{2(n+2)}{n-2}}_{t,x}$ and the map from the initial data to the solution is uniformly
continuous in these norms. If in addition the energy is small, it is known that the solution exists globally in
time\footnote{One should compare this result to the case of the focusing energy-critical NLS, where an
argument of Glassey \cite{glassey} shows that certain Schwartz data will blow up in finite time; for instance, this will
occur whenever the potential energy exceeds the kinetic energy.} and scattering occurs (see \cite{TV}); that is,
there exist solutions
$u_\pm$ of the free Schr\"odinger equation $(i \partial_t +\Delta)u_\pm = 0$ such that
$\|u(t) - u_\pm (t)\|_{\ho_x} \rightarrow 0$ as $t \rightarrow \pm \infty$. However, for initial data with large energy,
the local well-posedness arguments do not extend to give global well-posedness.

Global well-posedness in $\dot{H}^1_x(\R^3)$ for the energy-critical NLS in the case of
large (but finite) energy, radially-symmetric initial data was first obtained by J. Bourgain (\cite{borg:scatter},
\cite{borg:book}) and subsequently by M. Grillakis, \cite{grillakis:scatter}. T. Tao, \cite{tao: gwp radial},
settled the problem for arbitrary dimensions (with an improvement in the final bound due to a simplification of the
argument), but again only for radially symmetric data. A major breakthrough in the field was made by J.~Colliander,
M.~Keel, G.~Staffilani, H.~Takaoka, and T.~Tao in \cite{ckstt:gwp} where they obtained global well-posedness
and scattering for the
energy-critical NLS in dimension $n=3$ with arbitrary initial data. Recently, E.~Ryckman and M.~Visan, \cite{RV},
obtained global well-posedness and scattering for the energy-critical NLS in dimension $n=4$;
the argument follows closely the one in \cite{ckstt:gwp}, but the derivation of the frequency-localized interaction
Morawetz inequality is significantly simpler and yields an improvement in the final bound.

The main goal of this paper is to extend to dimensions $n\geq 5$ the global well-posedness results
for \eqref{schrodinger equation} in the energy space:

\begin{theorem}\label{THE THEOREM}
For any $u_0$ with finite energy (i.e., $E(u_0) < \infty$) there exists a unique global solution
$u \in C_t^0 \ho_x \cap L^{\frac{2(n+2)}{n-2}}_{t,x}$ to \eqref{schrodinger equation}
such that
\begin{equation}\label{l6 bounds}
\int_{-\infty}^{\infty} \int_{\R^n} |u(t,x)|^{\frac{2(n+2)}{n-2}} dx dt \leq C(E(u_0))
\end{equation}
for some constant $C(E(u_0))$ depending only on the energy.
\end{theorem}

As is well known (see, for instance, \cite{TV}), the $L_{t,x}^{\frac{2(n+2)}{n-2}}$ bound
also gives scattering and asymptotic completeness:

\begin{corollary}
Let $u_0$ have finite energy and let $u$ be the unique global solution in
$C_t^0 \ho_x \cap L^{\frac{2(n+2)}{n-2}}_{t,x}$ to \eqref{schrodinger equation}. Then there exist finite energy solutions
$u_{\pm}$ to the free Schr\"odinger equation $(i\partial_t+\Delta)u_{\pm}=0$ such that
$$
\|u_{\pm}(t)-u(t)\|_{\dot{H}^1_x}\rightarrow 0 \ \ \text{as} \ \ t\rightarrow \pm \infty.
$$
Furthermore, the maps $u_0 \mapsto u_{\pm}(0)$ are homeomorphisms from $\dot{H}^1(\R^n)$ to $\dot{H}^1(\R^n)$.
\end{corollary}

As we will see, treating dimensions $n\geq 5$ introduces several new difficulties relative to \cite{ckstt:gwp} and
\cite{RV}. For the most part, these stem from the small power of the nonlinearity in \eqref{schrodinger equation}.
For example, $|u|^{\frac{4}{n-2}}u$ is not a smooth function of $u, \bar{u}$ for any $n\geq 5$, which immediately
implies the failure of persistence of regularity.
Moreover, as the power of the nonlinearity is no longer an even integer in dimensions $n\geq 5$, the difference
of two nonlinearities cannot be written as a polynomial. Instead, we will make use of the following inequalities:
Let $F:\C\to\C$ be given by $F(z)=|z|^{\frac{4}{n-2}}z$. Then,
\begin{align}\label{derivate}
F_z(z)=\tfrac{n}{n-2} |z|^{\frac{4}{n-2}}\quad \text{and} \quad F_{\bar z}(z)=\tfrac{2}{n-2} |z|^{\frac{4}{n-2}}\tfrac{z^2}{|z|^2},
\end{align}
where $F_z$, $F_{\bar{z}}$ are the usual complex derivatives
$$
F_z := \frac{1}{2}\Bigl(\frac{\partial F}{\partial x} - i\frac{\partial F}{\partial y}\Bigr),
\quad F_{\bar{z}} := \frac{1}{2}\Bigl(\frac{\partial F}{\partial x} + i\frac{\partial F}{\partial y}\Bigr).
$$
Note that in dimensions $n>6$, the functions $z\mapsto F_z(z)$ and $z\mapsto F_{\bar z}(z)$ are no longer Lipschitz
continuous; however, they are H\"older continuous of order $\tfrac{4}{n-2}$. Thus, writing
\begin{equation}\label{diff integral form}
F(u+v)-F(v)
= \int_0^1 \bigl[F_{z}\bigl(v+\theta u\bigr)u+F_{\bar{z}}\bigl(v+\theta u\bigr)\bar{u}\bigr]d\theta,
\end{equation}
we estimate
\begin{align}\label{diff1}
\bigl|F(u+v)-F(v)\bigr|\lesssim |u||v|^{\frac{4}{n-2}} + |u|^{\frac{n+2}{n-2}}
\end{align}
and
\begin{align}\label{diff2}
\bigl|F(u+v)-F(u)-F(v)\bigr|\lesssim
\begin{cases}
|u||v|^{\frac{4}{n-2}}, & |u|\leq |v|\\
|v||u |^{\frac{4}{n-2}}, & |v|< |u|.
\end{cases}
\end{align}
Moreover, by the chain rule and the Lipschitz/H\"older continuity of the derivatives $F_z$ and $F_{\bar z}$, we get
\begin{align}\label{diff3}
\bigl| \nabla \bigl(F(u+v)-F(u) - F(v)\bigr)\bigr| \lesssim
\begin{cases}
\ \bigl(|\nabla u||v|+ |\nabla v||u|\bigr) \bigl(|u|^{\frac{1}{3}}+|v|^{\frac{1}{3}}\bigr), & n=5\\
\ |\nabla u||v|^{\frac{4}{n-2}} + |\nabla v| |u|^{\frac{4}{n-2}}, &
n\geq 6.
\end{cases}
\end{align}

\subsection{Outline of the proof of Theorem \ref{THE THEOREM}}

Our argument follows the scheme of \cite{ckstt:gwp} and we summarize it below.

For an energy $E \geq 0$ we define the quantity $M(E)$ by
$$
M(E):=\sup\ \|u\|_{L^{\frac{2(n+2)}{n-2}}_{t,x}(I\times\R^n)},
$$
where $I \subset \R$ ranges over all compact time intervals and $u$ ranges over all
$\dot{S}^1$ solutions\footnote{See Sections 1.2 and 2.1 for the notation and definitions appearing in the outline of the proof.}
to \eqref{schrodinger equation} on $I\times \R^n$ with $E(u)\leq E$.  For $E<0$ we define $M(E)=0$ since,
of course, there are no negative energy solutions.

{}From the local well-posedness theory (see Lemma \ref{lemma long time}), we know that \eqref{schrodinger equation}
is locally wellposed in $\dot{S}^1$. Moreover, from the global well-posedness theory for small initial data,
we see that $M(E)$ is finite for small energies $E$. Our task is to show that $M(E) < \infty$ for all $E>0$
as Theorem \ref{THE THEOREM} follows from this claim by a standard argument. More precisely, given initial data
$u_0$ with energy $E$, we let
$$
\Omega_1=\{T: \exists u \in \dot{S}^1([0,T]\times\R^n)\ \text{solving \eqref{schrodinger equation} with}\
   \|u\|_{\dot{S}^1([0,T]\times\R^n)}\leq C_1(E)\}
$$
and
$$
\Omega_2=\{T: \exists u \in \dot{S}^1([0,T]\times\R^n)\ \text{solving \eqref{schrodinger equation} with}\
   \|u\|_{\dot{S}^1([0,T]\times\R^n)}<\infty\}.
$$
Here $C_1(E)=C(E, M(E))E$ and $C(E, M(E))$ is the constant from Lemma~\ref{persistence of regularity}.

Note that by definition and Fatou's lemma, $\Omega_1$ is a closed set. By the local well-posedness theory
(see Lemma~\ref{blow}), if $T\in \Omega_1$ then there exists $\eps$ sufficiently small such that
$[T,T+\eps]\subset \Omega_2$. In particular, as $0$ lies in $\Omega_1$, we get that a small neighbourhood of $0$,
say $[0,\eps]$, lies in $\Omega_2$. Hence, to obtain a global solution to \eqref{schrodinger equation} it suffices to
see that $\Omega_2\subset\Omega_1$. By the definition of $M(E)$, given $T\in \Omega_2$ we immediately get
$$
\|u\|_{L_{t,x}^{\frac{2(n+2)}{n-2}}([0,T]\times\R^n)}\leq M(E).
$$
Combining this estimate with Lemma \ref{persistence of regularity} we obtain $T\in \Omega_1$.

We will prove that $M(E)<\infty$ by contradiction. Assume $M(E)$ is not always finite.  From perturbation theory
(see Lemma \ref{lemma long time}), we see that the set $\{E : M(E)<\infty\}$ is open. Since it is also connected and
contains zero, there exists a \emph{critical energy}, $0 < E_{crit} < \infty$, such that $M(E_{crit}) = \infty$ but
$M(E) < \infty$ for all $E < E_{crit}$. From the definition of $E_{crit}$ and Lemma~\ref{persistence of regularity}, we get

\begin{lemma}[Induction on energy hypothesis]\label{lemma induct on energy}
Let $t_0 \in \R$ and let $v(t_0)$ be an $\dot{H}^1_x$ function with $E(v(t_0)) \leq E_{crit}-\eta$ for some $\eta > 0$.
Then there exists a global $\dot{S}^1$ solution $v$ to \eqref{schrodinger equation} on $\rr$ with initial data $v(t_0)$
at time $t_0$, such that
$$
\|v\|_{L^{\frac{2(n+2)}{n-2}}_{t,x}(\R\times\R^n)} \leq M(E_{crit} - \eta) .
$$
Moreover, we have
$$
\|v\|_{\sor} \leq C( E_{crit}-\eta, M(E_{crit} - \eta)).
$$
\end{lemma}

We will need a few small parameters for the contradiction argument. Specifically, we will need
$$
1 \gg \eta_0 \gg \eta_1 \gg \eta_2 \gg \eta_3 \gg \eta_4 \gg \eta_5> 0
$$
where each $\eta_j$ is allowed to depend on the critical energy and any of the larger $\eta$'s. We will choose
$\eta_j$ small enough such that, in particular, it will be smaller than any constant depending on the
previous $\eta$'s used in the argument.

As $M(E_{crit})$ is infinite, given any $\eta_5> 0$ there exist a compact interval $I_* \subset \R$ and an $\dot{S}^1$
solution $u$ to \eqref{schrodinger equation} on $\is$ with
$ E(u) \leq E_{crit}$ but
\begin{equation}\label{HUGE}
\|u\|_{L_{t,x}^{\frac{2(n+2)}{n-2}} (\is)} > 1/\eta_5.
\end{equation}
Note that we may assume $E(u)\geq \frac{1}{2}E_{crit}$, since otherwise we would get
$$
\|u\|_{L_{t,x}^{\frac{2(n+2)}{n-2}} (\is)} \leq M(\tfrac{1}{2}E_{crit})<\infty
$$
and we would be done.

This suggests we make the following definition:

\begin{definition}
A \emph{minimal energy blowup solution} to \eqref{schrodinger equation} is an $\dot{S}^1$ solution $u$ on a time interval
$I_*\subset \R$ with energy
\begin{equation}\label{minimal energy blowup solution}
\tfrac{1}{2}E_{crit} \leq E(u(t)) \leq E_{crit}
\end{equation}
and huge $L_{t,x}^{\frac{2(n+2)}{n-2}}$-norm in the sense of \eqref{HUGE}.
\end{definition}


Note that conservation of energy together with \eqref{minimal energy blowup solution} and Sobolev embedding imply
\begin{equation}\label{potential energy bound}
\|u\|_{L^\infty_t L^{\frac{2n}{n-2}}_x (\is)} \lesssim 1
\end{equation}
and also
\begin{equation}\label{kinetic energy bound}
\|u\|_{L^\infty_t \ho_x (\is)} \sim 1,
\end{equation}
where, following our standard convention, the constants are allowed to depend on $E_{crit}$.

In Section 2 we recall the standard linear Strichartz estimates that we will use throughout the proof of Theorem
\ref{THE THEOREM}. We also record the inhomogeneous Strichartz estimates that will be useful in deriving the
frequency-localized interaction Morawetz inequality. Finally, we refine the bilinear Strichartz estimates of
\cite{ckstt:gwp} using a lemma of M. Christ and A. Kiselev. The main application of the bilinear Strichartz estimates is to
control the interaction between high frequencies, $P_{hi}u:=u_{hi}$, and low frequencies, $P_{lo}u:=u_{lo}$, when deriving
the frequency localization result. In fact, because of the small power of the nonlinearity in higher dimensions, we
have to control interactions between $u_{hi}$ and fractional powers of $u_{lo}$ (and vice versa); this is dealt
with via interpolation and the refined bilinear Strichartz estimates.


In Section 3, we record perturbation results from \cite{TV} that we will use repeatedly in the proof of
Theorem~\ref{THE THEOREM}.

In Section 4 we prove various localization and concentration results. We expect that a minimal energy blowup solution
should be localized in both physical and frequency space.  For if not, it could be decomposed into two essentially
separate solutions, each with strictly smaller energy than the original.
By Lemma \ref{lemma induct on energy} we can then extend these smaller energy solutions to all of $I_*$.  As each of the
separate evolutions exactly solves \eqref{schrodinger equation}, we expect their sum to solve \eqref{schrodinger equation}
approximately.  We could then use the perturbation theory results (specifically Lemma~\ref{lemma induct on energy})
to bound the $L^{\frac{2(n+2)}{n-2}}_{t,x}$-norm of $u$ in terms of $\eta_0, \eta_1, \eta_2, \eta_3, \text{ and } \eta_4$,
thus contradicting the fact that $\eta_5$ can be chosen arbitrarily small in \eqref{HUGE}.

The spatial concentration result follows in a similar manner, but is a bit more technical. To derive it, we use an
idea of Bourgain, \cite{borg:scatter}, and restrict our analysis to a subinterval $I_0 \subset I_*$. We need to use
both the frequency localization result and the fact that the potential energy of a minimal energy blowup solution
is bounded away from zero in order to prove spatial concentration.

In Section 5 we obtain the frequency-localized interaction Morawetz inequality \eqref{flim}, which will be used
to derive a contradiction to the frequency localization and spatial concentration results just described.

A typical example of a Morawetz inequality for \eqref{schrodinger equation} is the bound
$$
\int_I \int_{\R^n}\frac{|u(t,x)|^{\frac{2n}{n-2}}}{|x|}dxdt\lesssim \sup_{t\in I}\|u(t)\|_{\dot{H}^{1/2}(\R^n)}^2
$$
for all time intervals $I$ and all sufficiently regular solutions $u:I\times\R^n\to \C$.

This estimate is not particularly useful for the energy-critical problem since the $\dot{H}^{1/2}_x$-norm is
supercritical with respect to the scaling \eqref{scaling}. To get around this problem, J. Bourgain and M. Grillakis
introduced a spatial cutoff obtaining the variant
$$
\int_I \int_{|x|\leq A|I|^{1/2}}\frac{|u(t,x)|^{\frac{2n}{n-2}}}{|x|}dxdt\lesssim A|I|^{1/2} E(u)
$$
for all $A\geq 1$, where $|I|$ denotes the length of the time interval $I$. While this estimate is better suited for
the energy-critical NLS (it involves the energy on the right-hand side), it only prevents concentration of $u$ at the spatial
origin $x=0$. This is especially useful in the spherically-symmetric case $u(t,x)=u(t,|x|)$, since the spherical symmetry
combined with the bounded energy assumption can be used to show that $u$ cannot concentrate at any location other than the
spatial origin. However, it does not provide much information about the solution away from the origin.
Following \cite{ckstt:gwp}, we develop a frequency-localized interaction Morawetz inequality which is better suited to
handle nonradial solutions.

While the previously mentioned Morawetz inequalities were \textit{a priori} estimates, the frequency-localized
interaction Morawetz inequality we will develop is not; it only applies to minimal
energy blowup solutions. While our model in obtaining this estimate is \cite{ckstt:gwp}, there are two significant
differences. We manage to avoid localizing in space (which adds significantly to the complexity of \cite{ckstt:gwp});
however, the low power of the nonlinearity necessitates decomposing the high-frequency portions of the minimal energy
blowup solution into a `good' part, which is in $\dot{S}^0\bigcap\dot{S}^1$, and a `bad' part
which lives outside the Strichartz trapezoid.
While having slower decay in time than the `good' part, the `bad' part has better
spatial decay. This splitting of the high frequencies together with H\"older-type estimates (used as a substitute for
the standard fractional chain rule) enable us to control the error terms appearing in the frequency-localized interaction
Morawetz inequality. This machinery is employed to derive \eqref{flim} in dimensions $n\geq 6$. In dimension $n=5$,
the derivation of the frequency-localized interaction Morawetz inequality is somewhat simpler (for details see
\cite{my thesis}). One should mention that the method used to obtain this inequality in
dimension $n=5$ also works in dimensions $6$, $7$, and $8$; in dimensions $n\geq 9$ the small power of the nonlinearity
causes the argument to fail.

A corollary of \eqref{flim} (in all dimensions $n\geq 5$) is good $L^3_tL_x^{\frac{6n}{3n-4}}$ control over the
high-frequency part of a minimal energy blowup solution. One then has to use Sobolev embedding to bootstrap this
to $L^{\frac{2(n+2)}{n-2}}_{t,x}$ control. However, it is also necessary to
make sure that the solution is not shifting its energy from low to high frequencies causing the
$L^{\frac{2(n+2)}{n-2}}_{t,x}$-norm to blow up while the $L^3_tL_x^{\frac{6n}{3n-4}}$-norm stays bounded. This is done
in Section~6, where we prove a frequency-localized mass almost conservation law that prevents energy evacuation
to high modes.

We put all these pieces together in Section~7 where the contradiction argument is concluded.

\textbf{Acknowledgments}: I would like to thank my advisor, Terence Tao, for his existence and uniqueness.
I am also grateful to Rowan Killip for helpful comments.

\subsection{Notation}
We will often use the notation $X \lesssim Y$ whenever there exists some constant $C$, possibly depending on the
critical energy and the dimension $n$ but not on any other parameters, so that $X \leq CY$. Similarly we will write $X \sim Y$ if
$X \lesssim Y \lesssim X$.  We say $X \ll Y$ if $X \leq cY$ for some small constant $c$, again possibly depending on the
critical energy and the dimension $n$. We will use the abbreviation $\O(X)$ to denote a quantity that resembles $X$, that is, a finite linear
combination of terms that look like $X$, but possibly with some factors replaced by their complex conjugates. We also
use the notation $\langle x\rangle :=(1+|x|^2)^{1/2}.$ We will use the notation $X+ := X + \eps$, for some
$0<\eps \ll 1$; similarly $X- := X-\eps$. The derivative operator $\nabla$ refers to the space variable only.
We will occasionally write subscripts to denote spatial derivatives and will use the summation convention over
repeated indices.

We define the Fourier transform on $\R^n$ to be
$$
\hat f(\xi) := \int_{\R^n} e^{-2 \pi i x \cdot \xi} f(x) dx.
$$

We will make frequent use of the fractional differentiation operators $|\nabla|^s$ defined by
$$
\widehat{|\nabla|^sf}(\xi) := |\xi|^s \hat f (\xi).
$$
These define the homogeneous Sobolev norms
$$
\|f\|_{\dot H^s_x} := \| |\nabla|^s f \|_{L^2_x}.
$$

Let $e^{it\Delta}$ be the free Schr\"odinger propagator.  In physical space this is given by the formula
$$
e^{it\Delta}f(x) = \frac{1}{(4 \pi i t)^{n/2}} \int_{\R^n} e^{i|x-y|^2/4t} f(y) dy,
$$
while in frequency space one can write this as
\begin{equation}\label{fourier rep}
\widehat{e^{it\Delta}f}(\xi) = e^{-4 \pi^2 i t |\xi|^2}\hat f(\xi).
\end{equation}
In particular, the propagator preserves the above Sobolev norms and obeys the \emph{dispersive inequality}
\begin{equation}\label{dispersive ineq}
\|e^{it\Delta}f\|_{L^\infty_x} \lesssim
|t|^{-\frac{n}{2}}\|f\|_{L^1_x}
\end{equation}
for all times $t\neq 0$.  We also recall \emph{Duhamel's formula}
\begin{align}\label{duhamel}
u(t) = e^{i(t-t_0)\Delta}u(t_0) - i \int_{t_0}^t e^{i(t-s)\Delta}(iu_t + \Delta u)(s) ds.
\end{align}

We will also need some Littlewood-Paley theory.  Specifically, let $\varphi(\xi)$ be a smooth symmetric bump supported in the ball
$|\xi| \leq 2$ and equalling one on the ball $|\xi| \leq 1$.  For each dyadic number $N \in 2^\Z$ we define the
Littlewood-Paley operators
\begin{align*}
\widehat{P_{\leq N}f}(\xi) &:=  \varphi(\xi/N)\hat f (\xi),\\
\widehat{P_{> N}f}(\xi) &:=  (1-\varphi(\xi/N))\hat f (\xi),\\
\widehat{P_N f}(\xi) &:=  [\varphi(\xi/N) - \varphi (2 \xi /N)] \hat f (\xi).
\end{align*}
Similarly we can define $P_{<N}$, $P_{\geq N}$, and $P_{M < \cdot \leq N} := P_{\leq N} - P_{\leq M}$, whenever $M$ and
$N$ are dyadic numbers.  We will frequently write $f_{\leq N}$ for $P_{\leq N} f$ and similarly for the other operators.

The Littlewood-Paley operators commute with derivative operators, the free propagator, and complex conjugation.
They are self-adjoint and bounded on every $L^p_x$ and $\dot H^s_x$ space for $1 \leq p \leq \infty$ and $s\geq 0$.  They
also obey the following Sobolev and Bernstein estimates that we will use repeatedly:
\begin{align*}
\|P_{\geq N} f\|_{L^p_x} &\lesssim N^{-s} \| |\nabla|^s P_{\geq N} f \|_{L^p_x},\\
\| |\nabla|^s  P_{\leq N} f\|_{L^p_x} &\lesssim N^{s} \| P_{\leq N} f\|_{L^p_x},\\
\| |\nabla|^{\pm s} P_N f\|_{L^p_x} &\sim N^{\pm s} \| P_N f \|_{L^p_x},\\
\|P_{\leq N} f\|_{L^q_x} &\lesssim N^{\frac{n}{p}-\frac{n}{q}} \|P_{\leq N} f\|_{L^p_x},\\
\|P_N f\|_{L^q_x} &\lesssim N^{\frac{n}{p}-\frac{n}{q}} \| P_N f\|_{L^p_x},
\end{align*}
whenever $s \geq 0$ and $1 \leq p \leq q \leq \infty$.

For instance, we can use the above Bernstein estimates and the kinetic energy bound \eqref{kinetic energy bound} to
control the mass at high frequencies
\begin{equation}\label{mass high freq bound}
\|P_{>M}u\|_{L^2(\R^n)} \lesssim \frac{1}{M} \quad \text{ for all }M \in 2^{\Z}.
\end{equation}

For any dyadic frequency $N\in 2^\Z$, the kernel of the operator $P_{\leq N}$ is not positive.  To resolve this problem, we introduce an operator
$P_{\leq N}'$.  More precisely, if $K_{\leq N}$ is the kernel associated to $P_{\leq N}$, we let $P_{\leq N}'$ be the operator associated to
$N^{-n}(K_{\leq N})^2$. Please note that since $\varphi(\xi)$ is symmetric, $K_{\leq N}\in \R$ and thus $N^{-n}(K_{\leq N})^2\geq 0$.  Moreover, as
$$
[N^{-n}(K_{\leq N})^2]\widehat{\ }(\xi)= N^{-n} \varphi(\xi/N)*\varphi(\xi/N),
$$
the kernel of $P_{\leq N}'$ is bounded in $L_x^1$ independently of $N$.  Therefore, the operator $P_{\leq N}'$ is bounded on every $L_x^p$
for $1\leq p\leq \infty$. Furthermore, for $s \geq 0$ and $1 \leq p \leq q \leq \infty$, we have
\begin{align*}
\| |\nabla|^s  P_{\leq N}' f\|_{L^p_x} &\lesssim N^{s} \| P_{\leq N}' f\|_{L^p_x},\\
\|P_{\leq N}' f\|_{L^q_x} &\lesssim N^{\frac{n}{p}-\frac{n}{q}} \|P_{\leq N}' f\|_{L^p_x}.
\end{align*}

%
%
%
%

\section{Strichartz numerology}
In this section we recall the Strichartz estimates and develop bilinear Strichartz estimates in $\R^{1+n}$.

We use $\lqlr$ to denote the spacetime norm
$$
\|u\|_{\llr}=\|u\|_{q,r} :=\Bigl(\int_{\R}\Bigl(\int_{\R^n} |u(t,x)|^r dx \Bigr)^{q/r} dt\Bigr)^{1/q},
$$
with the usual modifications when $q$ or $r$ is infinity, or when the domain $\R \times \R^n$ is replaced by some smaller
spacetime region.  When $q=r$ we abbreviate $\lqlr$ by $L^q_{t,x}$.

\subsection{Linear and bilinear Strichartz estimates}\label{strichartz section}
We say that a pair of exponents $(q,r)$ is Schr\"odinger-\emph{admissible} if $\tfrac{2}{q} + \tfrac{n}{r} = \frac{n}{2}$
and $2 \leq q,r \leq \infty$. If $I \times \R^n$ is a spacetime slab, we define the $\szi$
\emph{Strichartz norm} by
\begin{equation}\label{s0}
\|u\|_{\szi} := \sup \Bigl(\sum_N \| P_N u \|^2_{\lli}\Bigr)^{1/2}
\end{equation}
where the supremum is taken over all admissible pairs $(q,r)$.  For $s>0$ we also define the $\dot S^s (\ir)$ \emph{Strichartz norm}
to be
$$
\|u\|_{\dot S^s (\ir)} := \| |\nabla|^s u \|_{\szi}.
$$

We observe the inequality
\begin{equation}\label{square sum}
\Bigl\|\Bigl(\sum_N |f_N|^2 \Bigr)^{1/2}\Bigr\|_{\lli} \leq \Bigl(\sum_N \|f_N\|^2_{\lli}\Bigr)^{1/2}
\end{equation}
for all $2 \leq q,r \leq \infty$ and arbitrary functions $f_N$, which one proves by interpolating between the trivial
cases $(2,2)$, $(2,\infty)$, $(\infty,2)$, and $(\infty,\infty)$. In particular, \eqref{square sum} holds for all
admissible exponents $(q,r)$.  Combining this with the Littlewood-Paley inequality, we find
\begin{align*}
\| u \|_{\lli}& \lesssim \Bigl\|\Bigl(\sum_N |P_N u|^2\Bigr)^{1/2}\Bigr\|_{\lli}\\
            & \lesssim \Bigl(\sum_N \|P_N u \|^2_{\lli}\Bigr)^{1/2}\\
            & \lesssim \| u \|_{\szi},
\end{align*}
which in particular implies
\begin{equation}\label{grad less s1}
\|\nabla u \| _{\lli} \lesssim \|u\|_{\soi}.
\end{equation}

In fact, by \eqref{grad less s1} and Sobolev embedding, the $\so$ norm controls the following spacetime norms:

\begin{lemma}\label{lemma strichartz norms}
For any $\dot{S^1}$ function $u$ on $\ir$, we have
\begin{align*}\label{strichartz norms}
&\|\nabla u\|_{\infty,2} + \|\nabla u\|_{3,\frac{6n}{3n-4}}+ \|\nabla u \|_{\frac{2(n+2)}{n-2},\frac{2n(n+2)}{n^2+4}} +
\|\nabla u\|_{\frac{2(n+2)}{n},\frac{2(n+2)}{n}} + \|\nabla u\|_{2,\frac{2n}{n-2}} \\
&\quad + \|u\|_{\infty,\frac{2n}{n-2}} + \|u\|_{3,\frac{6n}{3n-10}} + \|u\|_{\frac{2(n+2)}{n-2},\frac{2(n+2)}{n-2}} +\|u\|_{\frac{2(n+2)}{n},\frac{2n(n+2)}{n^2-2n-4}} + \|u\|_{2,\frac{2n}{n-4}}\\
&\phantom{\quad + \|u\|_{\infty,\frac{2n}{n-2}} +}\lesssim \|u\|_{\so}
\end{align*}
where all spacetime norms are on $\ir$.
\end{lemma}

\begin{figure}[ht]
\begin{center}
\setlength{\unitlength}{1.75mm}
\begin{picture}(75,65)(-4,-4)
\put(0,0){\vector(1,0){68}}\put(67,-2.5){$\frac1r$}
\put(0,0){\vector(0,1){60}}\put(-2.0,58){$\frac1q$}
\put(25,25){\circle*{0.7}}\put(20,24.5){$(4,4)$}
\put(50,0){\circle*{0.7}}\put(50,-2.2){$(\infty,2)$}
\put(43.333,50){\circle*{0.7}}\put(37,51){$(2,\tfrac{2n}{n-2})$}
\qbezier(50,0)(46.667,25)(43.333,50)
\put(36.667,50){\circle*{0.7}} \put(29,51){$(2,\tfrac{2n}{n-4})$}
\put(43.333,0){\circle*{0.7}}\put(35,-2.2){$(\infty,\tfrac{2n}{n-2})$} 
\qbezier(36.667,50)(40,25)(43.333,0)
\put(38.235,38.235){\circle*{0.7}} \put(23.5,37.735){$(\tfrac{2(n+2)}{n-2},\tfrac{2(n+2)}{n-2})$}
\put(46.667,50){\circle*{0.7}}\put(47,51){$(2,\tfrac{2n}{n-1})$} 
\put(45.641,50){\circle*{0.7}}\put(45.641,54.4){\vector(0,-1){3.4}}\put(42,56){$(2,\tfrac{2n(n-2)}{n^2-3n-2})$}
\qbezier[100](45.641,50)(47.821,25)(50,0)
\put(46.222,43.333){\circle*{0.7}}\put(48,42.833){$(\tfrac{2n}{n-2},\tfrac{2n^2}{(n+1)(n-2)})$}
\put(45.882,38.235){\circle*{0.7}}\put(48,37.735){$(\tfrac{2(n+2)}{(n-2)},\tfrac{2n(n+2)}{(n-2)(n+3)})$}
\qbezier(46.667,50)(45.882,38.235)(43.333,0)
\put(45.555,33.333){\circle*{0.7}}\put(48,32.833){$(3,\tfrac{6n}{3n-4})$}
\put(58.974,50){\circle*{0.7}}\put(58.974,51){$(2,\tfrac{2n(n-2)}{n^2+n-10})$}
\put(56.666,50){\circle*{0.7}}\put(56.666,54.4){\vector(0,-1){3.4}}\put(54.666,56){$(2,\tfrac{2n}{n+2})$}
\end{picture}
Figure 1: The Strichartz trapezoid $(n > 10)$.
\end{center}
\end{figure}

Next, we recall the Strichartz estimates:

\begin{lemma}\label{lemma linear strichartz}
Let $I$ be a compact time interval, $s\geq 0$, and let $u : \ir \rightarrow \C$ be a solution to the forced Schr\"odinger equation
\begin{equation*}
i u_t + \Delta u = \sum_{m=1}^M F_m
\end{equation*}
for some functions $F_1 ,\dots,F_M$.  Then,
\begin{equation}
\||\nabla|^s u\|_{\dot S^0(\ir)} \lesssim \|u(t_0)\|_{\dot H^s (\R^n)} + \sum_{m=1}^M \||\nabla|^s F_m \|_{L^{q'_m}_t L^{r'_m}_x (\ir)}
\end{equation}
for any time $t_0 \in I$ and any admissible exponents $(q_1,r_1),\dots,(q_m,r_m)$. As
usual, $p'$ denotes the dual exponent to $p$, that is, $1/p + 1/p' = 1$.
\end{lemma}

\begin{proof}
To prove Lemma \ref{lemma linear strichartz}, let us first make the following reductions. We note that it suffices to
take $M=1$, since the claim for general $M$ follows from Duhamel's formula and the triangle inequality. We can also take $s$
to be 0, since the estimate for $s>0$ follows by applying $|\nabla|^s$ to both sides of the equation and
noting that $|\nabla|^s$ commutes with $i\partial_t+\Delta$. As the Littlewood-Paley operators also commute with
$i\partial_t+\Delta$, we have
$$
(i\partial_t+\Delta)P_{N}u=P_N F_1
$$
for all dyadic $N$'s. Applying the standard Strichartz estimates (see \cite{tao:keel}), we get
\begin{align}\label{P N Strichartz}
\|P_N u\|_{L_t^qL_x^r(I\times \R^n)}\lesssim \|P_Nu(t_0)\|_{L^2_x}+\|P_N F_1\|_{L_t^{q'_1}L_x^{r'_1}(I\times\R^n)}
\end{align}
for all admissible exponents $(q,r)$ and $(q_1,r_1)$. Squaring \eqref{P N Strichartz}, summing in $N$, using the definition
of the $\dot{S}^0$-norm and the Littlewood-Paley inequality, together with the dual of \eqref{square sum}, we get the
claim.
\end{proof}

We recall next the inhomogeneous Strichartz estimates. We say that the pair $(q,r)$ is Schr\"odinger-\emph{acceptable}
if $1\leq q,r\leq \infty$ and $\frac{1}{q}<n(\frac{1}{2}-\frac{1}{r})$, or $(q,r)=(\infty ,2)$. We have the following result,
which is a special case of Theorem 1.4 from \cite{foschi}:
\begin{theorem}\label{foschi}
Let $I$ be a compact time interval. Let $(q,r)$ and
$(\tilde{q},\tilde{r})$ be two Schr\"odinger-\emph{acceptable} pairs satisfying the scaling condition
$\frac{1}{q}+\frac{1}{\tilde{q}}=\frac{n}{2}\bigl(1-\frac{1}{r}-\frac{1}{\tilde{r}})$ and either
\begin{align*}
\frac{1}{q}+\frac{1}{\tilde{q}}=1, \quad
\frac{n-2}{n}<\frac{r}{\tilde{r}}<\frac{n}{n-2},\quad
\frac{1}{r}\leq \frac{1}{q},\quad\text{and}\quad
\frac{1}{\tilde{r}}\leq \frac{1}{\tilde{q}},
\end{align*}
or
\begin{align*}
\frac{1}{q}+\frac{1}{\tilde{q}}<1 \quad \text{and}\quad \frac{n-2}{n}\leq\frac{r}{\tilde{r}}\leq\frac{n}{n-2}.
\end{align*}
Then,
$$
\Bigl\|\int_{s<t}e^{i(t-s)\Delta} F(s)ds\Bigr\|_{L_t^qL_x^r(\ir)}\lesssim \|F\|_{L_t^{\tilde{q}'}L_x^{\tilde{r}'}(\ir)}.
$$
\end{theorem}

In particular, let us record the following inhomogeneous Strichartz estimates that we will use to derive the
frequency-localized interaction Morawetz inequality:
\begin{align}\label{inhomS}
\Bigl\|\int_{s<t}e^{i(t-s)\Delta} F(s)ds\Bigr\|_{L_t^2L_x^{\frac{2n(n-2)}{n^2-3n-2}}(\ir)}
\lesssim \|F\|_{L_t^2L_x^{\frac{2n(n-2)}{n^2+n-10}}(\ir)}
\end{align}
and
\begin{align}\label{inhomS2}
\Bigl\|\int_{s<t}e^{i(t-s)\Delta} F(s)ds\Bigr\|_{L_t^2L_x^{\frac{2n(n-2)}{n^2-3n-2}}(\ir)}
\lesssim \|F\|_{L_{t,x}^{\frac{2(n-2)(n+2)}{n^2+3n-14}}(\ir)}.
\end{align}
We leave it to the reader to check that the hypotheses of Theorem \ref{foschi} are satisfied for
$(q,r)=(2, \frac{2n(n-2)}{n^2-3n-2})$ and $(\tilde{q},\tilde{r})=(2,\frac{2n(n-2)}{n^2-5n+10})$
or $\tilde{q}=\tilde{r}=\frac{2(n-2)(n+2)}{n^2-3n+6}$, as long as $n\geq 5$.

For the remainder of this section, we develop bilinear Strichartz estimates that we will use later, in particular,
in deriving the
frequency localization result. We will adapt the bilinear Strichartz estimate obtained in \cite{ckstt:gwp}, which is
in turn a refinement of a Strichartz estimate due to J. Bourgain (see \cite{BRefine}), to better suit our nonlinearity.
The reason for which we need to make this modification is that the power of the nonlinearity gets smaller as the
dimension increases and we have no hope of placing it in $L_t^1\dot{H}^{s}_x$ for $n>6$. In order to achieve our
goal, we will use a lemma due to M.~Christ and A.~Kiselev,~\cite{christkiss}. The following version is from H. Smith
and C. Sogge, \cite{smithsogge}:
\begin{lemma}\label{CK}
Let $X,Y$ be Banach spaces and let $k(t,s)$ be the kernel of an operator $T:L^p([0,T];X)\to L^q([0,T];Y)$.
Define the lower triangular operator $\tilde{T}:L^p([0,T];X)\to L^q([0,T];Y)$ by
$$
\tilde{T}f(t)=\int_0^t k(t,s)f(s)ds.
$$
Then, the operator $\tilde{T}$ is bounded from $L^p([0,T];X)$ to $L^q([0,T];Y)$ and $\|\tilde{T}\|\lesssim \|T\|$,
provided $p<q$.
\end{lemma}

We are now ready to state and prove the following
\begin{lemma}\label{lemma bilinear strichartz}
Fix $n \geq 2$.  For any spacetime slab $I \times \R^n$, any $t_0 \in I$, and any $\delta > 0$, we have
\begin{multline}
\|uv\|_{L^2_{t,x} (I \times \R^n)} \leq
C(\delta)\Bigl(\|u(t_0)\|_{\dot H_x^{-1/2 +\delta}} + \||\nabla|^{-\frac{1}{2}+\delta}(i \partial_t + \Delta)u\|_{L^{q'}_t L^{r'}_x(I \times \R^n)}\Bigr)\\
\times \Bigl(\|v(t_0)\|_{\dot H_x^{\frac{n-1}{2}-\delta}} + \||\nabla|^{\frac{n-1}{2}-\delta}(i
\partial_t + \Delta)v\|_{L^{\tilde{q}'}_t L^{\tilde{r}'}_x(I \times \R^n)}\Bigr),
\end{multline}
for any Schr\"odinger-\emph{admissible} pairs $(q,r)$ and $(\tilde{q}, \tilde{r})$ satisfying $q,\tilde{q}>2$.
\end{lemma}

\begin{proof}
Throughout the proof all spacetime norms will be on the slab $\ir$. We define
$$
\|w\|_{k,q,r}:=\|w(t_0)\|_{\dot{H}_x^k} + \||\nabla|^k(i\partial_t+\Delta)w\|_{q',r'}
$$
and
$$
F_{k,q,r}=\{w :\ \|w\|_{k,q,r}< \infty\}.
$$
With this notation our goal is to show
$$
\|uv\|_{2,2}
   \leq C(\delta)\|u\|_{-\frac{1}{2}+\delta, q, r} \|v\|_{\frac{n-1}{2}-\delta, \tilde{q}, \tilde{r}},
$$
for any $(q,r)$ and $(\tilde{q}, \tilde{r})$ Schr\"odinger admissible pairs satisfying $q,\tilde{q}>2$.

The bilinear Strichartz estimate derived in \cite{ckstt:gwp} (see their Lemma 3.4) reads
\begin{align}\label{Iteambilinear}
\|uv\|_{2,2}
   \leq C(\delta)\|u\|_{-\frac{1}{2}+\delta, \infty, 2} \|v\|_{\frac{n-1}{2}-\delta, \infty, 2},
\end{align}
which proves the case $q=q'=\infty$. We will combine this result with Lemma \ref{CK} to obtain the full set
of exponents. A particular case of \eqref{Iteambilinear} is
\begin{align}\label{Iteambilinearfree}
\|e^{i(t-t_0)\Delta}u(t_0)e^{i(t-t_0)\Delta}v(t_0)\|_{2,2}
\leq C(\delta)\|u(t_0)\|_{\dot{H}_x^{-\frac{1}{2}+\delta}}\|v(t_0)\|_{\dot{H}_x^{\frac{n-1}{2}-\delta}}.
\end{align}

Now fix $(q,r)$ a Schr\"odinger admissible pair with $q>2$ and fix $v\in\emph{F}_{\frac{n-1}{2}-\delta, \infty, 2}$.
Consider the operator $u\mapsto uv$; we claim that this operator is bounded from $\emph{F}_{-\frac{1}{2}+\delta,q,r}$ to
$L_{t,x}^2$. Indeed, using Duhamel's formula for $u$ we estimate
$$
\|uv\|_{2,2}\leq \|e^{i(t-t_0)\Delta}u(t_0)v\|_{2,2}+
    \Bigl\|\Bigl(\int_{t_0}^t e^{i(t-s)\Delta}(i\partial_t+\Delta)u(s)ds\Bigr) v\Bigr\|_{2,2}.
$$
Using Duhamel's formula for $v$ and \eqref{Iteambilinearfree}, we get
\begin{align*}
\|e^{i(t-t_0)\Delta}u(t_0)v\|_{2,2}
&\lesssim \|e^{i(t-t_0)\Delta}u(t_0)e^{i(t-t_0)\Delta}v(t_0)\|_{2,2}\\
&\quad+\Bigl\|e^{i(t-t_0)\Delta}u(t_0)\int_{t_0}^t e^{i(t-s)\Delta}(i\partial_t+\Delta)v(s)ds\Bigr\|_{2,2}\\
&\leq C(\delta) \|u(t_0)\|_{\dot{H}_x^{-\frac{1}{2}+\delta}}\|v(t_0)\|_{\dot{H}_x^{\frac{n-1}{2}-\delta}}\\
&\quad+C(\delta) \|u(t_0)\|_{\dot{H}_x^{-\frac{1}{2}+\delta}}\int_{\R}\|(i\partial_t+\Delta)v(s)\|_{\dot{H}_x^{\frac{n-1}{2}-\delta}}ds\\
&\leq C(\delta) \|u(t_0)\|_{\dot{H}_x^{-\frac{1}{2}+\delta}}\|v\|_{\frac{n-1}{2}-\delta, \infty, 2}\\
&\leq C(\delta) \|u\|_{-\frac{1}{2}+\delta, q, r} \|v\|_{\frac{n-1}{2}-\delta, \infty, 2}.
\end{align*}
In order to conclude our claim, it suffices to see that
\begin{align}\label{Christkiss}
\Bigl\|\Bigl(\int_{t_0}^t e^{i(t-s)\Delta}(i\partial_t+\Delta)u(s)ds\Bigr) v\Bigr\|_{2,2}
   \leq C(\delta)\|u\|_{-\frac{1}{2}+\delta, q, r} \|v\|_{\frac{n-1}{2}-\delta, \infty, 2}.
\end{align}
By Lemma \ref{CK}, for $q>2$, \eqref{Christkiss} is implied by
$$
\Bigl\|\Bigl(\int_{\R} e^{i(t-s)\Delta}(i\partial_t+\Delta)u(s)ds\Bigr) v\Bigr\|_{2,2}
   \leq C(\delta)\|u\|_{-\frac{1}{2}+\delta, q, r} \|v\|_{\frac{n-1}{2}-\delta, \infty, 2}.
$$
But now, using again a Duhamel expansion for $v$ and proceeding as before, we get
\begin{align*}
\Bigl\| e^{it\Delta}\Bigl(\int_{\R} e^{-is\Delta}(i\partial_t+\Delta)&u(s)ds \Bigr)v\Bigr\|_{2,2}\\
&\lesssim \Bigl\|\int_{\R} e^{-is\Delta}(i\partial_t+\Delta)u(s)ds \Bigr\|_{\dot{H}_x^{-\frac{1}{2}+\delta}}\|v\|_{\frac{n-1}{2}-\delta, \infty, 2}.
\end{align*}
By the standard linear Strichartz estimates,
\begin{align*}
\Bigl\|\int_{\R} e^{-is\Delta}(i\partial_t+\Delta)u(s)ds \Bigr\|_{\dot{H}_x^{-\frac{1}{2}+\delta}}
\lesssim \||\nabla|^{-\frac{1}{2}+\delta}(i\partial_t+\Delta)u\|_{q',r'}
&\lesssim \|u\|_{-\frac{1}{2}+\delta,q,r},
\end{align*}
and \eqref{Christkiss} follows.

To conclude the proof of Lemma \ref{lemma bilinear strichartz}, we run the same argument for
$v \in \emph{F}_{\frac{n-1}{2}-\delta,\tilde{q},\tilde{r}}$ with $u\in \emph{F}_{-\frac{1}{2}+\delta,q,r}$ fixed.

\end{proof}

%
%
%
%

\section{Perturbation Theory}
In this section we review the local theory for \eqref{schrodinger equation}. As mentioned in the introduction,
the Cauchy problem for \eqref{schrodinger equation} is locally well-posed in $\ho(\R^n)$.
A large part of the local theory for the energy-critical NLS is due to
Cazenave and Weissler, \cite{cw0}, \cite{cwI}, who showed the existence of local solutions for large energy data
and that of global solutions for small energy data. As is to be expected for a critical equation, the time
of existence of the local solutions depends on the profile of the initial data and not only on its energy. They also
proved uniqueness of these solutions in certain Strichartz spaces in which the solution was shown to depend
continuously\footnote{For the defocusing energy-critical NLS the continuity was established in $L^q_t \dot{H}^1_x$
for any $q\leq \infty$.  See \cite{cazenave:book} for details.} on the initial data in the energy space $\dot H^1(\R^n)$.
A later argument of Cazenave, \cite{cazenave:book}, also demonstrates that the uniqueness
is in fact unconditional in the category of strong solutions (see also \cite{katounique},
\cite{twounique}, \cite{ckstt:gwp} for some related arguments).

These preliminary results are not completely satisfactory as the arguments that establish continuous dependence
on the data do not establish uniformly continuous dependence on the data in energy-critical spaces. In \cite{ckstt:gwp}
and \cite{RV}, it was shown that the dependence on the data is Lipschitz in dimensions $n=3$, respectively $n=4$,
results which extend nicely to treat dimensions $n=5,6$; see \cite{TV}. However, in dimensions $n>6$, the low power
of the nonlinearity causes the same argument to fail; in this case, the dependence on the data was shown to be H\"older
continuous (rather than Lipschitz) in \cite{TV}.

Closely related to the continuous dependence on the data and an essential tool for induction on energy type arguments
is the stability theory for the equation \eqref{schrodinger equation}. More precisely, given an \emph{approximate}
solution
\begin{equation}\label{equation 1-approx}
\begin{cases}
i \tilde u_t +\Delta \tilde u &= |\tilde u|^{\frac{4}{n-2}}\tilde u + e\\
\tilde u(t_0,x) &= \tilde u_0(x) \in \dot H^1(\R^n)
\end{cases}
\end{equation}
to \eqref{schrodinger equation}, with $e$ small in a suitable space and $\tilde u_0$ and $u_0$ close in $\dot H^1_x$,
is it possible to show that the solution $u$ to \eqref{schrodinger equation} stays very close to $\tilde u$?
Note that the question of continuous dependence on the data corresponds to the case $e=0$.  In dimensions $n=3,4$,
an analysis based on Strichartz estimates yields a satisfactory theory; see \cite{ckstt:gwp}, \cite{RV}. In the general
case, particularly $n>6$, a more careful analysis is needed; the relevant results were obtained by T. Tao and M. Visan,
\cite{TV}, and we record them below.

\begin{lemma}[Long-time perturbations]\label{lemma long time}
Let $I$ be a compact time interval and let $\util$ be an approximate solution to \eqref{schrodinger equation}
on $I\times\R^n$ in the sense that
$$
(i\partial_t+\Delta)\util=|\util|^{\frac{4}{n-2}}\util+e
$$
for some function $e$. Assume that
\begin{align}
\|\util\|_{L_{t,x}^{\frac{2(n+2)}{n-2}}(\ir)}&\leq M \label{finite S norm} \\
\|\util\|_{L_t^{\infty}\dot{H}^1_x(\ir)}&\leq E \label{finite energy}
\end{align}
for some constants $M, E>0$. Let $t_0\in I$ and let $u(t_0)$ close to $\util(t_0)$ in the sense that
\begin{align}\label{close}
\|u(t_0)-\util(t_0)\|_{\dot{H}^1_x}\leq E'
\end{align}
for some $E'>0$. Assume also the smallness conditions
\begin{align}
\Bigl(\sum_N \|P_N \nabla e^{i(t-t_0)\Delta}\bigl(u(t_0)-\util(t_0)\bigr)\|^2_{L_t^{\frac{2(n+2)}{n-2}}L_x^{\frac{2n(n+2)}{n^2+4}}(\ir)}\Bigr)^{1/2} &\leq \eps \label{closer} \\
\|\nabla e\|_{L_t^2L_x^{\frac{2n}{n+2}}(\ir)}&\leq \eps \label{error small}
\end{align}
for some $0<\eps \leq \eps_1$, where $\eps_1=\eps_1(E, E', M)$ is a small constant. Then there exists a solution
$u$ to \eqref{schrodinger equation} on $\ir$ with the specified initial data $u(t_0)$ at time $t=t_0$ satisfying
\begin{align}
\|\nabla(u-\util)\|_{L_t^{\frac{2(n+2)}{n-2}}L_x^{\frac{2n(n+2)}{n^2+4}}(\ir)}&\leq C(E ,E', M)\bigl(\eps+\eps^{\frac{7}{(n-2)^2}}\bigr) \label{close in L^p}\\
\|u-\util\|_{\dot{S}^1(\ir)}&\leq C(E ,E', M)\bigl(E'+\eps+ \eps^{\frac{7}{(n-2)^2}}\bigr) \label{close in S^1}\\
\|u\|_{\dot{S}^1(\ir)}&\leq C(E, E', M). \label{u in S^1}
\end{align}
Here, $C(E,E',M) > 0$ is a non-decreasing function of $E,E',M$, and the dimension~$n$.
\end{lemma}

\begin{remark}\label{redundant}
By Strichartz estimates and Plancherel's theorem, we have
\begin{align*}
\Bigl(\sum_N \|P_N \nabla e^{i(t-t_0)\Delta}\bigl(u(t_0)-\util(t_0)&\bigr)\|^2_{L_t^{\frac{2(n+2)}{n-2}}L_x^{\frac{2n(n+2)}{n^2+4}}(\ir)}\Bigr)^{1/2}\\
&\lesssim \Bigl(\sum_N \|P_N \nabla(u(t_0)-\util(t_0)\bigr)\|^2_{\infty,2}\Bigr)^{1/2}\\
&\lesssim \|\nabla (u(t_0)-\util(t_0)\bigr)\|_{\infty,2} \\
&\lesssim E'
\end{align*}
on the slab $\ir$; hence, the hypothesis \eqref{closer} is redundant if $E'=O(\eps)$.
\end{remark}

We end this section with a few related results. The first asserts that if a solution cannot be continued strongly
beyond a time $T_*$, then the $L_{t,x}^{\frac{2(n+2)}{n-2}}$-norm must blow up at that time.

\begin{lemma}[Standard blowup criterion, \cite{cw0}, \cite{cwI}, \cite{TV}]\label{blow}
Let $u_0\in \dot{H}^1_x$ and let $u$ be a strong $\dot S^1$ solution to \eqref{schrodinger equation} on the slab
$[t_0, T_0]\times\R^n$ such that
\begin{align}\label{norm finite}
\|u\|_{L_{t,x}^{\frac{2(n+2)}{n-2}}([t_0, T_0]\times\R^n)} < \infty.
\end{align}
Then there exists $\delta=\delta(u_0)>0$ such that the solution $u$ extends to a strong $\dot S^1$ solution to
\eqref{schrodinger equation} on the slab $[t_0, T_0+\delta]\times\R^n$.
\end{lemma}

The last result we mention here was used in the proof of Lemma~\ref{lemma long time} above and shows that once
we have $L^{\frac{2(n+2)}{n-2}}_{t,x}$ control of a finite-energy solution, we control all Strichartz norms as well.
Details can be found in \cite{TV}.

\begin{lemma}[Persistence of regularity]\label{persistence of regularity}
Let $0\leq s<1+\frac{4}{n-2}$, $I$ a compact time interval, and $u$ a finite-energy solution to \eqref{schrodinger equation} obeying
$$
\|u\|_{L^{\frac{2(n+2)}{n-2}}_{t,x}(I\times\R^n)} \leq M.
$$
Then, if $t_0 \in I$ and $u(t_0) \in \dot H_x^s$, we have
\begin{equation}\label{persistence of regularity eq}
\|u\|_{\dot S^s (\ir)} \leq C(M,E(u))\|u(t_0)\|_{\dot H_x^s}.
\end{equation}
\end{lemma}

\begin{proof}
We first consider the case $0\leq s\leq 1$.  We subdivide the interval $I$ into $N\sim (1+\frac{M}{\eta})^{\frac{2(n+2)}{n-2}}$
subintervals $I_j=[t_j,t_{j+1}]$ such that on each slab $I_j\times\R^n$ we have
\begin{align*}
\|u\|_{L_{t,x}^{\frac{2(n+2)}{n-2}}(I_j\times\R^n)}\leq\eta,
\end{align*}
where $\eta$ is a small positive constant to be chosen momentarily. By Strichartz and the fractional chain rule
(see for instance \cite{ChW:fractional chain rule, KP:fractional chain rule}), on each slab $I_j\times\R^n$ we obtain
\begin{align*}
\|u\|_{\dot S^s(I_j\times\R^n)}
&\lesssim \|u(t_j)\|_{\dot H^s_x} + \||\nabla|^s\bigr(|u|^{\frac{4}{n-2}}u\bigl)\|_{L_{t,x}^{\frac{2(n+2)}{n+4}}(I_j\times\R^n)}\\
&\lesssim \|u(t_j)\|_{\dot H^s_x} + \||\nabla|^s u\|_{L_{t,x}^{\frac{2(n+2)}{n}}(I_j\times\R^n)}\|u\|_{L_{t,x}^{\frac{2(n+2)}{n-2}}(I_j\times\R^n)}^{\frac{4}{n-2}}\\
&\lesssim \|u(t_j)\|_{\dot H^s_x} + \eta^{\frac{4}{n-2}}\|u\|_{\dot S^s(I_j\times\R^n)}.
\end{align*}
Choosing $\eta$ sufficiently small, we obtain
\begin{align*}
\|u\|_{\dot S^s(I_j\times\R^n)}\lesssim \|u(t_j)\|_{\dot H^s_x}.
\end{align*}
For the range of $s$ under discussion, i.e., $0\leq s\leq 1$, the conclusion \eqref{persistence of regularity eq} follows
by adding these estimates over all subintervals $I_j$.

We now consider the case $1<s<1+\frac{4}{n-2}$.  As the solution $u$ has finite energy, from the previous case with $s=1$,
we deduce that
\begin{align*}
\|\nabla u\|_{L_t^{\frac{2(n+2)}{n-2}}L_x^{\frac{2n(n+2)}{n^2+4}}(\ir)}
\lesssim \|u\|_{\dot S^1 (\ir)} \leq C(M,E(u))\|u(t_0)\|_{\dot H_x^1}.
\end{align*}
We subdivide the interval $I$ into $N\sim (1+\frac{C(M, E(u))}{\eta})^{\frac{2(n+2)}{n-2}}$ subintervals $I_j=[t_j,t_{j+1}]$ such that
on each slab $I_j\times\R^n$ we have
\begin{align*}
\|\nabla u\|_{L_t^{\frac{2(n+2)}{n-2}}L_x^{\frac{2n(n+2)}{n^2+4}}(I_j\times\R^n)}\leq\eta,
\end{align*}
where $\eta$ is a small positive constant to be chosen later.  By Sobolev embedding, we also have
\begin{align*}
\|u\|_{L_{t,x}^{\frac{2(n+2)}{n-2}}(I_j\times\R^n)}\lesssim\eta.
\end{align*}
By Strichartz, on each slab $I_j\times\R^n$ we obtain
\begin{align}\label{comput1}
\|u\|_{\dot S^s(I_j\times\R^n)}
&\lesssim \|u(t_j)\|_{\dot H^s_x} + \||\nabla|^s\bigr(|u|^{\frac{4}{n-2}}u\bigl)\|_{L_{t,x}^{\frac{2(n+2)}{n+4}}(I_j\times\R^n)} \notag\\
&\lesssim \|u(t_j)\|_{\dot H^s_x} + \||\nabla|^{s-1} \bigl(\nabla u F_z(u)+\nabla \bar{u} F_{\bar z}(u)\bigr)\|_{L_{t,x}^{\frac{2(n+2)}{n+4}}(I_j\times\R^n)}.
\end{align}
By the fractional product rule (see \cite{ChW:fractional chain rule, KP:fractional chain rule}) and Sobolev embedding, we estimate
\begin{align*}
\||\nabla|^{s-1} \bigl(&\nabla u F_z(u)\bigr)\|_{L_{t,x}^{\frac{2(n+2)}{n+4}}(I_j\times\R^n)}\\
&\lesssim \||\nabla|^s u\|_{L_{t,x}^{\frac{2(n+2)}{n}}(I_j\times\R^n)}\|F_z(u)\|_{L_{t,x}^{\frac{n+2}{2}}(I_j\times\R^n)}\\
&\quad +\|\nabla u\|_{L_t^{\frac{2(n+2)}{n}}L_x^{\frac{2n(n+2)}{n^2-2(n+2)(s-1)}}(I_j\times\R^n)} \||\nabla|^{s-1} F_z(u)\|_{L_t^{\frac{n+2}{2}}L_x^{\frac{n(n+2)}{2n+(n+2)(s-1)}}(I_j\times\R^n)}\\
&\lesssim \|u\|_{\dot S^s(I_j\times\R^n)}\|u\|_{L_{t,x}^{\frac{2(n+2)}{n-2}}(I_j\times\R^n)}^{\frac{4}{n-2}}\\
&\quad + \||\nabla|^{s}u\|_{L_{t,x}^{\frac{2(n+2)}{n}}(I_j\times\R^n)}\||\nabla|^{s-1} F_z(u)\|_{L_t^{\frac{n+2}{2}}L_x^{\frac{n(n+2)}{2n+(n+2)(s-1)}}(I_j\times\R^n)}\\
&\lesssim \eta^{\frac{4}{n-2}}\|u\|_{\dot S^s(I_j\times\R^n)}+\|u\|_{\dot S^s(I_j\times\R^n)}\||\nabla|^{s-1} F_z(u)\|_{L_t^{\frac{n+2}{2}}L_x^{\frac{n(n+2)}{2n+(n+2)(s-1)}}(I_j\times\R^n)}.
\end{align*}
Similarly,
\begin{align*}
\||\nabla|^{s-1} \bigl(&\nabla \bar u F_{\bar z}(u)\bigr)\|_{L_{t,x}^{\frac{2(n+2)}{n+4}}(I_j\times\R^n)}\\
&\lesssim \eta^{\frac{4}{n-2}}\|u\|_{\dot S^s(I_j\times\R^n)}+\|u\|_{\dot S^s(I_j\times\R^n)}\||\nabla|^{s-1} F_{\bar z}(u)\|_{L_t^{\frac{n+2}{2}}L_x^{\frac{n(n+2)}{2n+(n+2)(s-1)}}(I_j\times\R^n)},
\end{align*}
so \eqref{comput1} becomes
\begin{align}\label{comput2}
\|u\|_{\dot S^s(I_j\times\R^n)}
&\lesssim \|u(t_j)\|_{\dot H^s_x} +\eta^{\frac{4}{n-2}}\|u\|_{\dot S^s(I_j\times\R^n)}\\
&\quad +\|u\|_{\dot S^s(I_j\times\R^n)}\||\nabla|^{s-1} F_z(u)\|_{L_t^{\frac{n+2}{2}}L_x^{\frac{n(n+2)}{2n+(n+2)(s-1)}}(I_j\times\R^n)}\notag\\
&\quad +\|u\|_{\dot S^s(I_j\times\R^n)}\||\nabla|^{s-1} F_{\bar z}(u)\|_{L_t^{\frac{n+2}{2}}L_x^{\frac{n(n+2)}{2n+(n+2)(s-1)}}(I_j\times\R^n)}.\notag
\end{align}
In order to estimate
\begin{align*}
\||\nabla|^{s-1} F_z(u)\|_{L_t^{\frac{n+2}{2}}L_x^{\frac{n(n+2)}{2n+(n+2)(s-1)}}(I_j\times\R^n)}
\end{align*}
and
\begin{align*}
\||\nabla|^{s-1} F_{\bar z}(u)\|_{L_t^{\frac{n+2}{2}}L_x^{\frac{n(n+2)}{2n+(n+2)(s-1)}}(I_j\times\R^n)},
\end{align*}
we will exploit the H\"older continuity of the functions $z\mapsto F_z(z)$ and $z\mapsto F_{\bar z}(z)$.
Using Proposition~\ref{fdfp} with $\alpha:=\tfrac{4}{n-2}$, $\sigma:=s-1$, $s:=r$, and $p=\tfrac{n(n+2)}{2n+(n+2)(s-1)}$, and applying Sobolev embedding,
we get
\begin{align*}
\||\nabla|^{s-1} F_{ z}(u)\|_{L_x^{\frac{n(n+2)}{2n+(n+2)(s-1)}}}
&\lesssim \|u\|^{\frac{4}{n-2}-\frac{s-1}{r}}_{L_x^{(\frac{4}{n-2}-\frac{s-1}{r})p_1}}\||\nabla|^r u\|^{\frac{s-1}{r}}_{L_x^{\frac{s-1}{r}p_2}}\\
&\lesssim \|u\|^{\frac{4}{n-2}-\frac{s-1}{r}}_{L_x^{(\frac{4}{n-2}-\frac{s-1}{r})p_1}}\|\nabla u\|^{\frac{s-1}{r}}_{L_x^{\frac{n(s-1)p_2}{nr+(1-r)(s-1)p_2}}}.
\end{align*}
Choosing $p_2$ such that $\tfrac{n(s-1)p_2}{nr+(1-r)(s-1)p_2}=\tfrac{2n(n+2)}{n^2+4}$ and applying H\"older's inequality with respect
to time, on the slab $I_j\times\R^n$ we get
\begin{align*}
\||\nabla|^{s-1} F_z(u)\|_{L_t^{\frac{n+2}{2}}L_x^{\frac{n(n+2)}{2n+(n+2)(s-1)}}}
&\lesssim \|u\|_{L_{t,x}^{\frac{2(n+2)}{n-2}}}^{\frac{4}{n-2}-\frac{s-1}{r}}
   \|\nabla u\|_{L_t^{\frac{2(n+2)}{n-2}}L_x^{\frac{2n(n+2)}{n^2+4}}}^{\frac{s-1}{r}}\\
&\lesssim \eta^{\frac{4}{n-2}}.
\end{align*}
Similarly,
\begin{align*}
\||\nabla|^{s-1} F_{\bar z}(u)\|_{L_t^{\frac{n+2}{2}}L_x^{\frac{n(n+2)}{2n+(n+2)(s-1)}}}\lesssim \eta^{\frac{4}{n-2}},
\end{align*}
and hence, returning to our previous computation, i.e., \eqref{comput2}, we obtain
\begin{align*}
\|u\|_{\dot S^s(I_j\times\R^n)}
\lesssim \|u(t_j)\|_{\dot H^s_x} +\eta^{\frac{4}{n-2}}\|u\|_{\dot S^s(I_j\times\R^n)}.
\end{align*}
Choosing $\eta$ sufficiently small, we get
\begin{align*}
\|u\|_{\dot S^s(I_j\times\R^n)} \lesssim \|u(t_j)\|_{\dot H^s_x}.
\end{align*}
The claim \eqref{persistence of regularity eq} follows by adding these bounds over all time intervals $I_j$.
\end{proof}

%
%
%
%

\section{Frequency Localization and Space Concentration}

Recall from the introduction that we expect a minimal energy blowup solution to be localized in both frequency and space.
In this section we will prove that this is indeed the case (we will not actually prove that the solution is localized in
space, just that it concentrates; see the discussion after the proof of Corollary \ref{lemma freq loc}). The first step
is the following proposition:

\begin{proposition}[Frequency delocalization $\Rightarrow$ spacetime bound]\label{lemma freq loc imp spacetime bds}
Let $\eta > 0$ and suppose there exist a dyadic frequency $N_{lo}>0$ and a time $t_0 \in I_*$ such that we have the
energy separation conditions
\begin{equation}\label{freq deloc lo freq}
\|P_{\leq N_{lo}} u(t_0)\|_{\ho_x} \geq \eta
\end{equation}
and
\begin{equation}\label{freq deloc hi freq}
\|P_{\geq K(\eta)N_{lo}} u(t_0)\|_{\ho_x} \geq \eta.
\end{equation}
If $K(\eta)$ is sufficiently large depending on $\eta$, we have
\begin{equation}\label{freq deloc ls}
\|u\|_{L_{t,x}^{\frac{2(n+2)}{n-2}} (I_* \times \R^n)}\leq C(\eta).
\end{equation}
\end{proposition}

\begin{proof}
Let $0 < \eps = \eps(\eta) \ll 1$ be a small number to be chosen later.  If $K(\eta)$ is sufficiently large depending on
$\eps$, then one can find $\eps^{-2}$ disjoint intervals $[\eps^2 N_j, \eps^{-2} N_j]$ contained in
$[N_{lo}, K(\eta)N_{lo}]$.  By \eqref{kinetic energy bound} and the pigeonhole principle, we may find an $N_j$ such that
the localization of $u(t_0)$ to the interval $[\eps^2 N_j, \eps^{-2}N_j]$ has very little energy:
\begin{align}\label{no middle}
\|P_{\eps^2 N_j \leq \cdot \leq \eps^{-2}N_j}u(t_0)\|_{\ho_x} \lesssim \eps.
\end{align}
As both the statement and conclusion of the proposition are invariant under the scaling \eqref{scaling}, we normalize
$N_j = 1$.

Define $\ulo(t_0) := P_{\leq \eps}u(t_0)$ and $\uhi(t_0) = P_{\geq \eps^{-1}}u(t_0)$.  We claim that $\uhi$ and $\ulo$
have smaller energy than $u$.

\begin{lemma}\label{lemma uhi ulo smaller energy}
If $\eps$ is sufficiently small depending on $\eta$, we have
$$
E(\ulo(t_0)), \ E(\uhi(t_0)) \leq E_{crit} - c\eta^C.
$$
\end{lemma}

\begin{proof}
We will prove this for $\ulo$; the proof for $\uhi$ is similar. Define $\uhip(t_0) := P_{> \eps}u(t_0)$ so that
$u(t_0) = \ulo (t_0) + \uhip(t_0)$ and consider the quantity
\begin{equation}\label{freq deloc energy 1}
|E(u(t_0)) - E(\ulo (t_0)) - E(\uhip(t_0))|.
\end{equation}
By the definition of energy, we can bound \eqref{freq deloc energy 1} by
\begin{equation}\label{freq deloc energy 2}
|\langle \nabla \ulo (t_0),\nabla \uhip (t_0) \rangle|
 + \Bigl|\int_{\R^n}\bigl(|u(t_0)|^{\frac{2n}{n-2}}-|\ulo(t_0)|^{\frac{2n}{n-2}}-|\uhip(t_0)|^{\frac{2n}{n-2}}\bigr) dx\Bigr|.
\end{equation}
We deal with the potential energy term first.  We have the pointwise estimate
$$
\bigl||u|^{\frac{2n}{n-2}}-|\ulo|^{\frac{2n}{n-2}}-|\uhip|^{\frac{2n}{n-2}}\bigr| \lesssim
\begin{cases}
|\uhip||\ulo|^{\frac{n+2}{n-2}}, & |\uhip|\leq |\ulo|\\
|\ulo ||\uhip |^{\frac{n+2}{n-2}}, & |\ulo|\leq |\uhip|.
\end{cases}
$$
Take the case $|\uhip(t_0)|\leq|\ulo(t_0)|$ and use H\"older to estimate
$$
\|\uhip(t_0)|\ulo(t_0)|^{\frac{n+2}{n-2}}\|_{L^1_x}\lesssim \|\uhip(t_0)\|_{L^2_x}\|\ulo(t_0)\|_{L^{\frac{2(n+2)}{n-2}}_x}^{\frac{n+2}{n-2}}.
$$
An application of Bernstein and Sobolev embedding yields
\begin{align*}
\|\ulo (t_0)\|_{L^{\frac{2(n+2)}{n-2}}_x}
\lesssim \eps^{\frac{n-2}{n+2}} \|\ulo (t_0)\|_{L^{\frac{2n}{n-2}}_x}
\lesssim \eps^{\frac{n-2}{n+2}} \|\ulo (t_0)\|_{\dot{H}^1_x}
   \lesssim \eps^{\frac{n-2}{n+2}}.
\end{align*}
Similarly, by Bernstein and \eqref{no middle},
\begin{align*}
\|\uhip (t_0)\|_{L^2_x}
   &\lesssim \sum_{N > \eps} \|P_N u(t_0)\|_{L^2_x}
     \lesssim \sum_{N > \eps}N^{-1} \|P_N u (t_0)\|_{\ho_x}\\
   &\lesssim \sum_{N > \eps^{-2}} N^{-1} + \sum_{\eps < N \leq \eps^{-2}}N^{-1}\eps
     \lesssim \eps^2 + \eps\eps^{-1}\\
   &\lesssim 1.
\end{align*}
Thus, for $|\uhip(t_0)|\leq|\ulo(t_0)|$,
$$
\|\uhip(t_0)|\ulo(t_0)|^{\frac{n+2}{n-2}}\|_{L^1_x}\lesssim \eps.
$$

Now take the case $|\ulo(t_0)|\leq|\uhip(t_0)|$ and use H\"older and the previous estimates on $\|\uhip (t_0)\|_{L^2_x}$
to get
\begin{align*}
\|\ulo(t_0)|\uhip(t_0)|^{\frac{n+2}{n-2}}\|_{L^1_x}
&\lesssim \||\ulo(t_0)|^{\frac{4}{n-2}}|\uhip(t_0)|^2\|_{L^1_x}
 \lesssim \|\ulo(t_0)\|_{L^{\infty}_x}^{\frac{4}{n-2}}\|\uhip(t_0)\|_{L^2_x}^2\\
&\lesssim \|\ulo(t_0)\|_{L^{\infty}_x}^{\frac{4}{n-2}}.
\end{align*}
Another application of Bernstein plus Sobolev embedding yields
\begin{align*}
\|\ulo (t_0)\|_{L^{\infty}_x}
\lesssim \eps^{\frac{n-2}{2}} \|\ulo (t_0)\|_{L^{\frac{2n}{n-2}}_x}
\lesssim \eps^{\frac{n-2}{2}} \|\ulo (t_0)\|_{\dot{H}^1_x}
    \lesssim \eps^{\frac{n-2}{2}}.
\end{align*}
Hence, if $|\ulo(t_0)|\leq|\uhip(t_0)|$, we have $\|\ulo(t_0)|\uhip(t_0)|^{\frac{n+2}{n-2}}\|_{L^1_x}\lesssim \eps^2.$
Combining the two cases, we get control over the potential energy term in \eqref{freq deloc energy 2}:
$$
\Bigl|\int_{\R^n}\bigl(|u(t_0)|^{\frac{2n}{n-2}}-|\ulo(t_0)|^{\frac{2n}{n-2}}-|\uhip(t_0)|^{\frac{2n}{n-2}}\bigr) dx\Bigr|\lesssim \eps.
$$

Next, we deal with the kinetic energy part of \eqref{freq deloc energy 2}.
We estimate
\begin{align*}
|\langle \nabla \ulo (t_0),\nabla \uhip (t_0) \rangle|
   &\lesssim |\langle \nabla P_{>\eps}P_{\leq\eps}u(t_{0}), \nabla u(t_0)\rangle| \\
   &\lesssim \|\nabla P_{>\eps}P_{\leq\eps}u(t_{0})\|_{L^2_x} \|\nabla u(t_0)\|_{L^2_x}.
\end{align*}
As
$$
\|\nabla u(t_0)\|_{L^2_x}\lesssim 1
$$
and
\begin{align*}
\|\nabla P_{>\eps} P_{\leq\eps}u(t_{0})\|_{L^2_x}
   &=\|(\nabla P_{>\eps}P_{\leq\eps}u(t_{0}))^{\wedge}\|_{L^2_x} \\
   &=\|\xi \varphi(\xi/\eps)(1-\varphi(\xi/\eps))\widehat{u(t_0)}(\xi)\|_{L^2_x} \\
   &\lesssim \eps \|\widehat{u_{hi'}(t_0)}\|_{L^2_x}
     \lesssim \eps,
\end{align*}
we obtain control over the kinetic energy term in \eqref{freq deloc energy 2},
$$
|\langle \nabla \ulo (t_0),\nabla \uhip (t_0) \rangle|\lesssim \eps.
$$

Therefore $\eqref{freq deloc energy 1} \lesssim \eps$. As
$$
E(u) \leq E_{crit}
$$
and, by hypothesis,
$$
E(\uhip(t_0))\gtrsim \|\uhip(t_0)\|^2_{\ho_x} \gtrsim \eta^2,
$$
the triangle inequality implies $E(\ulo(t_0)) \leq E_{crit} -c\eta^C$, provided we choose $\eps$ sufficiently small.

Similarly, one proves $E(\uhi(t_0)) \leq E_{crit} -c\eta^C$.
\end{proof}

Now, since $E(\ulo(t_0)),E(\uhi(t_0)) \leq \ecrit - c\eta^C < \ecrit$, we can apply Lemma~\ref{lemma induct on energy} to deduce that there exist $\dot{S}^1$ solutions $\ulo$ and $\uhi$ on the slab $I_*\times\R^n$ with initial data $\ulo (t_0)$ and
$\uhi (t_0)$ such that
\begin{align}
\|\ulo\|_{\sois} &\leq C(\eta) \label{uloso}
\end{align}
and
\begin{align}
\|\uhi\|_{\sois} &\leq C(\eta)\label{uhiso}.
\end{align}
From Lemma~\ref{persistence of regularity}, we also have
\begin{align}
\|\ulo\|_{\dot S^{1+s}(\is)} &\lesssim C(\eta)\|u_{lo}(t_0)\|_{\dot
H^{1+s}}\leq C(\eta)\eps^s, \ \ \ \forall \ 0<s<\tfrac{4}{n-2}
\label{ulos1+s}
\end{align}
and
\begin{align}
\|\uhi\|_{\szis} &\lesssim C(\eta)\|u_{hi}(t_0)\|_{L_x^2}\leq C(\eta)\eps \label{uhisz}.
\end{align}

Define $\tilde u := \ulo + \uhi$.  We claim that $\tilde u$ is a near-solution to \eqref{schrodinger equation}.

\begin{lemma}\label{lemma uhi plus ulo near soln}
We have
$$i \tilde u_t + \Delta \tilde u = |\tilde u|^{\frac{4}{n-2}} \tilde u - e$$
where the error $e$ obeys the bound
\begin{equation}\label{freq deloc error bound}
\|\nabla e\|_{L^2_t L^{\frac{2n}{n+2}}_x (\is)} \lesssim C(\eta)\eps^{\frac{3}{(n-2)(2n+1)}}.
\end{equation}
\end{lemma}

\begin{proof}
In order to estimate one derivative of the error term
$$
e=|\tilde{u}|^{\frac{4}{n-2}}\tilde{u}-|u_{lo}|^{\frac{4}{n-2}}u_{lo}-|u_{hi}|^{\frac{4}{n-2}}u_{hi}
 =F(u_{lo}+u_{hi})-F(u_{lo})-F(u_{hi}),
$$
we use \eqref{diff3} to obtain
\begin{align*}
\|\nabla e\|_{L^2_t L^{\frac{2n}{n+2}}_x (\is)}
\lesssim \|\nabla u_{lo}|u_{hi}|^{\frac{4}{n-2}}\|_{L^2_t L^{\frac{2n}{n+2}}_x (\is)} + \|\nabla u_{hi}|u_{lo}|^{\frac{4}{n-2}}\|_{L^2_t L^{\frac{2n}{n+2}}_x (\is)}
\end{align*}
in dimension $n\geq 6$ and
\begin{align*}
\|\nabla e\|_{L^2_t L^{\frac{10}{7}}_x (\is)}
&\lesssim \|\nabla u_{lo}|u_{hi}||u_{lo}|^{\frac{1}{3}}\|_{L^2_t L^{\frac{10}{7}}_x (\is)} + \|\nabla u_{hi}|u_{lo}||u_{hi}|^{\frac{1}{3}}\|_{L^2_t L^{\frac{10}{7}}_x (\is)}\\
&\quad  +\|\nabla u_{lo}|u_{hi}|^{\frac{4}{3}}\|_{L^2_t L^{\frac{10}{7}}_x (\is)} + \|\nabla u_{hi}|u_{lo}|^{\frac{4}{3}}\|_{L^2_t L^{\frac{10}{7}}_x (\is)}
\end{align*}
in dimension $n=5$.

Thus, proving \eqref{freq deloc error bound} amounts to showing
\begin{align}
\|\nabla u_{lo}|u_{hi}|^{\frac{4}{n-2}}\|_{L^2_t L^{\frac{2n}{n+2}}_x (\is)}
&\lesssim C(\eta)\eps^{\frac{3}{(n-2)(2n+1)}},  \ \ \ \ n\geq 5 \label{error1}\\
\|\nabla u_{hi}|u_{lo}|^{\frac{4}{n-2}}\|_{L^2_t L^{\frac{2n}{n+2}}_x (\is)}
&\lesssim C(\eta)\eps^{\frac{3}{(n-2)(2n+1)}},  \ \ \ \ n\geq 5  \label{error2}\\
\|\nabla u_{lo}|u_{hi}||u_{lo}|^{\frac{1}{3}}\|_{L^2_t L^{\frac{10}{7}}_x (\is)}
&\lesssim C(\eta)\eps^{\frac{1}{11}}, \ \ \ \ n=5 \label{error3}\\
\|\nabla u_{hi}|u_{lo}||u_{hi}|^{\frac{1}{3}}\|_{L^2_t L^{\frac{10}{7}}_x (\is)}
&\lesssim C(\eta)\eps^{\frac{1}{11}}, \ \ \ \ n=5. \label{error4}
\end{align}
To prove \eqref{error1}, we make use of H\"older, interpolation, \eqref{uloso}, \eqref{uhiso}, \eqref{ulos1+s}, and \eqref{uhisz}:
\begin{align*}
\|\nabla u_{lo}|u_{hi}|^{\frac{4}{n-2}}\|_{2,\frac{2n}{n+2} }
&\lesssim \||u_{hi}|^{\frac{4}{n-2}}\|_{\infty,\frac{n(n-2)}{2(n-1)}}\|\nabla u_{lo}\|_{2,\frac{2n(n-2)}{n^2-4n}} \\
&\lesssim \|u_{hi}\|_{\infty,\frac{2n}{n-1}}^{\frac{4}{n-2}}\||\nabla|^{1+\frac{2}{n-2}} u_{lo}\|_{2,\frac{2n}{n-2}} \\
&\lesssim \|u_{hi}\|_{\infty,2}^{\frac{2}{n-2}}\|u_{hi}\|_{\infty,\frac{2n}{n-2}}^{\frac{2}{n-2}}\||\nabla|^{1+\frac{2}{n-2}} u_{lo}\|_{\dot S^0(\is)} \\
&\lesssim \|u_{hi}\|_{\dot S^0(\is)}^{\frac{2}{n-2}}\|u_{hi}\|_{\dot S^1(\is)}^{\frac{2}{n-2}}\|u_{lo}\|_{\dot S^{1+\frac{2}{n-2}}(\is)}\\
&\lesssim C(\eta)\eps^{\frac{4}{n-2}}\\
&\lesssim C(\eta)\eps^{\frac{3}{(n-2)(2n+1)}},
\end{align*}
where all spacetime norms are on $\is$.
We turn now towards \eqref{error2}; on $\is$, we estimate
\begin{align*}
\|\nabla u_{hi} |u_{lo}|^{\frac{4}{n-2}}\|_{2,\frac{2n}{n+2}}
\lesssim \bigl\|\bigl|\nabla u_{hi} u_{lo} |u_{lo}|^{2n}\bigr|^{\frac{n}{(n-2)(n+1)}}\bigr\|_{\frac{2(n-2)(n+1)}{n},1}^{\frac{4(n+1)}{n(2n+1)}}\|\nabla u_{hi}\|_{2,\frac{2n}{n-2}}^{1-\frac{4}{(n-2)(2n+1)}}.
\end{align*}
Using \eqref{uhiso}, we get
$$
\|\nabla u_{hi}\|_{L^2_tL_x^{\frac{2n}{n-2}}(\is)}^{1-\frac{4}{(n-2)(2n+1)}}
  \lesssim \| u_{hi}\|_{\dot{S}^1(\is)}^{1-\frac{4}{(n-2)(2n+1)}}
  \lesssim C(\eta).
$$
Fixing $t\in I_*$, we use H\"older to estimate
\begin{align*}
\int_{\R^n}\bigl|\nabla &u_{hi}(t) u_{lo}(t) |u_{lo}|^{2n}(t)\bigr|^p dx\\
&\lesssim \sum_{N_1\leq\dots\leq N_{2n+1}}\int_{\R^n}\bigr|\nabla u_{hi}(t)P_{N_1}u_{lo}(t)P_{N_2}u_{lo}(t)\cdots P_{N_{2n+1}}\overline{u_{lo}}(t)\bigr|^p dx\\
&\lesssim \sum_{N_1\leq\dots\leq N_{2n+1}}\|\nabla u_{hi}(t)P_{N_1}u_{lo}(t)\|_{L^2_x}^p \|P_{N_2}u_{lo}(t)\|_{L_x^r}^p\cdots \|P_{N_{2n+1}}u_{lo}(t)\|_{L_x^r}^p,
\end{align*}
where we denoted $p:=\frac{n}{(n-2)(n+1)}$ and $r:=\frac{4n^2}{2n^2-3n-4}$.

Integrating with respect to time, on the slab $\is$ we get
\begin{align}\label{???}
\bigl\|\bigl|\nabla u_{hi}& u_{lo} |u_{lo}|^{2n}\bigr|^{\frac{n}{(n-2)(n+1)}}\bigr\|_{L_t^\frac{2(n-2)(n+1)}{n}L_x^1} \notag\\
&\lesssim \sum_{N_1\leq\dots\leq N_{2n+1}}\|\nabla u_{hi}P_{N_1}u_{lo}\|_{L_{t,x}^2}^p
\|P_{N_2}u_{lo}\|_{L_t^\infty L_x^r}^p\cdots \|P_{N_{2n+1}}u_{lo}\|_{L_t^\infty L_x^r}^p
\end{align}

By Bernstein,
\begin{align*}
\|P_N u_{lo}\|_{L_t^\infty L_x^r}\lesssim N^{\frac{3n+4}{4n}}\|P_N u_{lo}\|_{L_t^\infty L_x^2}.
\end{align*}
Thus, in view of \eqref{uloso} and \eqref{ulos1+s}, we obtain
\begin{align}\label{indiv project}
\|P_N u_{lo}\|_{L_t^\infty L_x^r}\lesssim C(\eta) \min(N^{-\frac{n-4}{4n}},N^{-\frac{n-4}{4n}-s}\eps^s)
\end{align}
for all $0<s<\tfrac{4}{n-2}$.

To bound $\|\nabla u_{hi}P_{N_1}u_{lo}\|_{L^2_{t,x}(\is)}$ we use the bilinear Strichartz estimates we have developed in Lemma \ref{lemma bilinear strichartz}. On $\is$ we estimate
\begin{align}\label{uhiun1}
\|\nabla u_{hi} P_{N_1}u_{lo}\|_{L^2_{t,x}}
\leq C(\delta)\Bigl(\|\nabla u_{hi}(t_0)\|_{\dot H_x^{-1/2 +\delta}} + \||\nabla|^{-\frac{1}{2}+\delta}(i \partial_t + \Delta)\nabla u_{hi}\|_{L^{\frac{2(n+2)}{n+4}}_{t,x}}\Bigr) \notag\\
\times \Bigl(\|P_{N_1}u_{lo}(t_0)\|_{\dot H_x^{\frac{n-1}{2}-\delta}} + \||\nabla|^{\frac{n-1}{2}-\delta}(i\partial_t + \Delta)P_{N_1}u_{lo}\|_{L^{\frac{2(n+2)}{n+4}}_{t,x}}\Bigr).
\end{align}
Interpolating between $ \|u_{hi}(t_0)\|_{L^2_x}\lesssim \eps$ and
$\|\nabla u_{hi}(t_0)\|_{L_x^2}\lesssim 1$,
we get
\begin{align}\label{!!}
\|\nabla u_{hi}(t_0)\|_{\dot H_x^{-1/2 +\delta}}\lesssim \eps^{\frac{1}{2}-\delta}.
\end{align}
Using \eqref{uhiso} and \eqref{uhisz}, we estimate
\begin{align*}
\|(i\partial_t+\Delta)u_{hi}\|_{L^{\frac{2(n+2)}{n+4}}_{t,x}(\is)}
&=\||u_{hi}|^{\frac{4}{n-2}}u_{hi}\|_{L^{\frac{2(n+2)}{n+4}}_{t,x}(\is)}\\
&\lesssim \|u_{hi}\|_{L^{\frac{2(n+2)}{n}}_{t,x}(\is)}\|u_{hi}\|_{L^{\frac{2(n+2)}{n-2}}_{t,x}(\is)}^{\frac{4}{n-2}}\\
&\lesssim \|u_{hi}\|_{\dot{S}^0(\is)}\|u_{hi}\|_{\dot{S}^1(\is)}^{\frac{4}{n-2}}\\
&\lesssim \eps C(\eta)
\end{align*}
and
\begin{align*}
\|\nabla (i\partial_t+\Delta)u_{hi}\|_{L^{\frac{2(n+2)}{n+4}}_{t,x}(\is)}
&=\|\nabla\bigl(|u_{hi}|^{\frac{4}{n-2}}u_{hi}\bigr)\|_{L^{\frac{2(n+2)}{n+4}}_{t,x}(\is)}\\
&\lesssim \|\nabla u_{hi}\|_{L^{\frac{2(n+2)}{n}}_{t,x}(\is)}\|u_{hi}\|_{L^{\frac{2(n+2)}{n-2}}_{t,x}(\is)}^{\frac{4}{n-2}}\\
&\lesssim \|u_{hi}\|_{\dot{S}^1(\is)}^{\frac{n+2}{n-2}} \lesssim C(\eta).
\end{align*}
Interpolating between the two estimate above, we obtain
$$
\||\nabla|^{-\frac{1}{2}+\delta}(i \partial_t + \Delta)\nabla u_{hi}\|_{L^{\frac{2(n+2)}{n+4}}_{t,x}(\is)}\lesssim C(\eta)\eps^{\frac{1}{2}-\delta}.
$$
Hence, combining this with \eqref{!!} gives
\begin{align}\label{uhiCK}
\|\nabla u_{hi}(t_0)\|_{\dot H_x^{-1/2 +\delta}} + \||\nabla|^{-\frac{1}{2}+\delta}(i \partial_t + \Delta)\nabla u_{hi}\|_{L^{\frac{2(n+2)}{n+4}}_{t,x}(\is)}
   \lesssim C(\eta)\eps^{\frac{1}{2}-\delta}.
\end{align}

We turn now to the factor in \eqref{uhiun1} containing $P_{N_1}u_{lo}$ and use Bernstein, \eqref{uloso}, and \eqref{ulos1+s}
to estimate
$$
\|P_{N_1}u_{lo}(t_0)\|_{\dot H_x^{\frac{n-1}{2}-\delta}}
  \lesssim C(\eta)\min(N_1^{\frac{n-3}{2}-\delta}, N_1^{\frac{n-3}{2}-\delta-s}\eps^s),
$$
for every $0<s<\tfrac{4}{n-2}$.  Similarly, by Bernstein and \eqref{uloso},
\begin{align*}
\||\nabla|^{\frac{n-1}{2}-\delta}(i\partial_t + \Delta)&P_{N_1}u_{lo}\|_{L^{\frac{2(n+2)}{n+4}}_{t,x}(\is)}\\
&\lesssim N_1^{\frac{n-3}{2}-\delta}\|\nabla(i\partial_t + \Delta)P_{N_1}u_{lo}\|_{L^{\frac{2(n+2)}{n+4}}_{t,x}(\is)}\\
&\lesssim N_1^{\frac{n-3}{2}-\delta}\|\nabla u_{lo}\|_{L^{\frac{2(n+2)}{n}}_{t,x}(\is)}\|u_{lo}\|_{L^{\frac{2(n+2)}{n-2}}_{t,x}(\is)}^{\frac{4}{n-2}}\\
&\lesssim N_1^{\frac{n-3}{2}-\delta}\|u_{lo}\|_{\dot{S}^1(\is)}^{\frac{n+2}{n-2}} \\
&\lesssim C(\eta) N_1^{\frac{n-3}{2}-\delta},
\end{align*}
while by Bernstein, \eqref{ulos1+s}, and the same arguments as in the proof of Lemma~\ref{persistence of regularity}, for any
$0<s<\tfrac{4}{n-2}$ we have
\begin{align*}
\||\nabla|^{\frac{n-1}{2}-\delta}(i\partial_t + \Delta)&P_{N_1}u_{lo}\|_{L^{\frac{2(n+2)}{n+4}}_{t,x}(\is)}\\
&\lesssim N_1^{\frac{n-3}{2}-\delta-s}\||\nabla|^{1+s}(i\partial_t + \Delta)u_{lo}\|_{L^{\frac{2(n+2)}{n+4}}_{t,x}(\is)}\\
&\lesssim N_1^{\frac{n-3}{2}-\delta-s} \|\nabla u_{lo}\|^{\frac{4}{n-2}}_{L_t^{\frac{2(n+2)}{n-2}}L_x^{\frac{2n(n+2)}{n^2+4}}(\is)} \|u_{lo}\|_{\dot{S}^{1+s}(\is)}\\
&\lesssim C(\eta) N_1^{\frac{n-3}{2}-\delta-s}\eps^s.
\end{align*}
Hence, for any $0<s<\tfrac{4}{n-2}$ we obtain
\begin{align}\label{uloCK}
\|P_{N_1}u_{lo}(t_0)\|_{\dot H_x^{\frac{n-1}{2}-\delta}}
+\||\nabla|^{\frac{n-1}{2}-\delta}(i\partial_t + \Delta)& P_{N_1}u_{lo}\|_{L^{\frac{2(n+2)}{n+4}}_{t,x}(\is)}\notag\\
&\lesssim C(\eta)\min(N_1^{\frac{n-3}{2}-\delta}, N_1^{\frac{n-3}{2}-\delta-s}\eps^s).
\end{align}

Thus, putting together \eqref{uhiun1}, \eqref{uhiCK}, and \eqref{uloCK}, we get
\begin{align}\label{?!}
\|\nabla u_{hi} P_{N_1}u_{lo}\|_{L^2_{t,x} (\is)}\lesssim C(\eta)\eps^{\frac{1}{2}-\delta}\min(N_1^{\frac{n-3}{2}-\delta}, N_1^{\frac{n-3}{2}-\delta-s}\eps^s).
\end{align}

Returning to our earlier computation, \eqref{???}, and using \eqref{indiv project} and \eqref{?!}, we conclude
\begin{align*}
\bigl\|\bigl|\nabla u_{hi} u_{lo}& |u_{lo}|^{\frac{n-4}{2}}\bigr|^{\frac{4}{n-2}}\bigr\|_{L_t^{\frac{n-2}{2}}L_x^1(\is)}\\
&\lesssim C(\eta)\eps^{(\frac{1}{2}-\delta)p}\sum_{N_1\leq\dots\leq N_{2n+1}}
\bigl[\min(N_1^{\frac{n-3}{2}-\delta}, N_1^{\frac{n-3}{2}-\delta-s}\eps^s)\bigr]^p\cdots\\
&\phantom{\lesssim C(\eta)\eps^{(\frac{1}{2}-\delta)p}\sum_{N_1\leq\dots\leq N_{2n+1}}}
\cdots \bigl[\min(N_{2n+1}^{-\frac{n-4}{4n}},N_{2n+1}^{-\frac{n-4}{4n}-s}\eps^s)\bigr]^p.
\end{align*}
In order to estimate the sum, we split it into three parts as follows:
\begin{align*}
\sum_{N_1\leq\cdots\leq N_{2n+1}}
&=\sum_{N_1\leq\cdots\leq N_{2n+1}\leq \eps}+\sum_{j=2}^{2n}\mathop{\sum_{N_1\leq\cdots \leq N_j\leq \eps}}_{\eps\leq N_{j+1}\leq \cdots\leq N_{2n+1}}
 +\sum_{\eps\leq N_1\leq\cdots\leq N_{2n+1}}\\
&=I+II+III.
\end{align*}
We have
\begin{align*}
I&\lesssim \sum_{N_1\leq\cdots\leq N_{2n+1}\leq \eps}N_1^{p(\frac{n-3}{2}-\delta)}N_2^{-p\frac{n-4}{4n}}\cdots N_{2n+1}^{-p\frac{n-4}{4n}}\\
 &\lesssim \sum_{N_1\leq \eps} N_1^{p(\frac{n-3}{2}-\delta)}N_1^{-p\frac{n-4}{2}}\\
 &\lesssim \eps^{(\frac{1}{2}-\delta)p}\\
II&\lesssim \sum_{j=2}^{2n}\mathop{\sum_{N_1\leq\cdots N_j\leq \eps}}_{\eps\leq N_{j+1}\leq \cdots\leq N_{2n+1}}
   N_1^{p(\frac{n-3}{2}-\delta)}N_2^{-p\frac{n-4}{4n}}\cdots N_j^{-p\frac{n-4}{4n}}\times\\
  &\phantom{\sum_{j=2}^{2n}\mathop{\sum_{N_1\leq\cdots\leq N_j\leq \eps}}_{\eps\leq N_{j+1}\leq \cdots\leq N_{2n+1}}N_1^{p(\frac{n-3}{2}-\delta)}}
   \qquad \times N_{j+1}^{-p(\frac{n-4}{4n}+s)}\eps^{sp}\cdots N_{2n+1}^{-p(\frac{n-4}{4n}+s)}\eps^{sp}\\
  &\lesssim \sum_{j=2}^{2n}\eps^{sp(2n+1-j)}\sum_{N_1\leq \eps\leq N_{j+1}}N_1^{p(\frac{n-3}{2}-\delta)-p\frac{n-4}{4n}(j-1)}N_{j+1}^{-p(\frac{n-4}{4n}+s)(2n+1-j)}\\
  &\lesssim \eps^{(\frac{1}{2}-\delta)p}\\
III&\lesssim \sum_{\eps\leq N_1\leq\cdots\leq N_{2n+1}}N_1^{p(\frac{n-3}{2}-\delta-s)}\eps^{sp}N_{2}^{-p(\frac{n-4}{4n}+s)}\eps^{sp}\cdots N_{2n+1}^{-p(\frac{n-4}{4n}+s)}\eps^{sp}\\
&\lesssim \eps^{sp(2n+1)}\sum_{\eps\leq N_1} N_1^{p(\frac{n-3}{2}-\delta-s)-p(\frac{n-4}{4n}+s)2n}\\
&\lesssim \eps^{sp(2n+1)}\sum_{\eps\leq N_1} N_1^{p[\frac{1}{2}-\delta-(2n+1)s]}\\
&\lesssim \eps^{(\frac{1}{2}-\delta)p},
\end{align*}
where the last inequality follows as soon as we choose $0<s<\tfrac{4}{n-2}$ such that $\frac{1}{2}-\delta-(2n+1)s<0$; in particular,
it holds for $s=\tfrac{1}{2(2n+1)}$.

Putting everything together, we obtain
$$
\|\nabla u_{hi}|u_{lo}|^{\frac{4}{n-2}}\|_{L^2_t L^{\frac{2n}{n+2}}_x (\is)}
\lesssim C(\eta)\eps^{(1-2\delta)\frac{4}{(n-2)(2n+1)}}
$$
and \eqref{error2} follows by choosing $\delta$ sufficiently small.

We examine next \eqref{error3}.  Using H\"older, interpolation, \eqref{uloso}, \eqref{uhiso}, \eqref{ulos1+s}, and \eqref{uhisz} we estimate
\begin{align*}
\|\nabla u_{lo} |u_{hi}||u_{lo}|^{\frac{1}{3}}\|_{2,\frac{10}{7}}
&\lesssim \|\nabla u_{lo}\|_{2,\frac{10}{3}} \|u_{hi}\|_{\infty, 3}\|u_{lo}\|_{\infty,5}^{\frac{1}{3}}\\
&\lesssim \|u_{lo}\|_{\dot S^1} \|u_{hi}\|_{\infty,2}^{\frac{1}{6}}\|u_{hi}\|_{\infty,\frac{10}{3}}^{\frac{5}{6}}\||\nabla|^{\frac{3}{2}}u_{lo}\|_{\infty,2}^{\frac{1}{3}}\\
&\lesssim C(\eta) \eps^{\frac{1}{6}}\eps^{\frac{1}{6}}\\
&\lesssim C(\eta)\eps^{\frac{1}{11}}.
\end{align*}

Finally, we consider \eqref{error4}.  By H\"older and conservation of energy, we estimate
\begin{align*}
\|\nabla u_{hi} |u_{lo}||u_{hi}|^{\frac{1}{3}}\|_{2,\frac{10}{7}}
&\lesssim \|\nabla u_{hi}\|_{\infty,\frac{10}{3}}^{\frac{1}{3}}\|\nabla u_{hi}|u_{lo}|\|_{2,\frac{5}{3}}
\lesssim \|\nabla u_{hi}|u_{lo}|\|_{2,\frac{5}{3}}.
\end{align*}
By interpolation,
\begin{align*}
\|\nabla u_{hi}|u_{lo}|\|_{2,\frac{5}{3}}
\lesssim \|\nabla u_{hi}\|_{2,\frac{10}{3}}^{\frac{1}{2}}\|\nabla u_{hi}|u_{lo}|^2\|_{2,\frac{10}{9}}^{\frac{1}{2}}
\lesssim C(\eta)\|\nabla u_{hi}|u_{lo}|^2\|_{2,\frac{10}{9}}^{\frac{1}{2}}.
\end{align*}
Using H\"older, Bernstein, \eqref{uloso}, \eqref{uhiso}, \eqref{ulos1+s}, \eqref{uhisz}, and \eqref{?!}, we estimate
\begin{align*}
\|\nabla u_{hi}|u_{lo}|^2\|_{2,\frac{10}{9}}
&\lesssim \sum_{N_1\leq N_2} \|\nabla u_{hi} P_{N_1} u_{lo}\|_{2,2}\|P_{N_2}u_{lo}\|_{\infty,\frac{5}{2}}\\
&\lesssim \sum_{N_1\leq N_2}C(\eta)\eps^{\frac{1}{2}-\delta}\min(N_1^{1-\delta}, N_1^{1-\delta-s}\eps^s)N_2^{-\frac{1}{2}}\|\nabla P_{N_2}u_{lo}\|_{\infty,2}\\
&\lesssim C(\eta)\eps^{\frac{1}{2}-\delta}\sum_{N_1\leq N_2}\min(N_1^{1-\delta}, N_1^{1-\delta-s}\eps^s)\min(N_2^{-\frac{1}{2}},N_2^{-\frac{1}{2}-s}\eps^s).
\end{align*}
Decomposing the sum into three sums, i.e., $\sum_{\eps\leq N_1\leq N_2}$, $\sum_{ N_1\leq \eps\leq N_2}$, and $\sum_{ N_1\leq N_2\leq \eps}$,
and taking $s>\frac{1}{4}$, we get
\begin{align*}
\|\nabla u_{hi}|u_{lo}|^2\|_{2,\frac{10}{9}}\lesssim C(\eta)\eps^{1-2\delta}.
\end{align*}
\eqref{error4} follows by taking $\delta$ sufficiently small.
\end{proof}

Next, we derive estimates on $u$ from those on $\tilde u$ via perturbation theory.  More precisely, we know
from \eqref{no middle} that
$$
\|u(t_0)-\tilde u(t_0)\|_{\ho_x} \lesssim \eps
$$
and hence, by Remark \ref{redundant},
$$
\Bigl(\sum_N \|P_N \nabla e^{i(t-t_0)\Delta}\bigl(u(t_0)-\util(t_0)\bigr)\|^2_{L_t^{\frac{2(n+2)}{n-2}}L_x^{\frac{2n(n+2)}{n^2+4}}(I_*\times\R^n)}\Bigr)^{1/2}\\
 \lesssim \eps.
$$
By Strichartz, we also have that
$$
\|\tilde{u}\|_{L_t^\infty \dot{H}^1_x(\is)}
 \lesssim \|\tilde{u}\|_{\dot{S}^1(\is)}\lesssim \|\ulo\|_{\so(I_* \times\R^n) } + \|\uhi\|_{\so(I_* \times\R^n)}
  \lesssim C(\eta)
$$
and hence,
$$
\|\tilde u \|_{L_{t,x}^{\frac{2(n+2)}{n-2}}(I_* \times\R^n)}
  \lesssim \|\tilde{u}\|_{\dot{S}^1(\is)}\lesssim C(\eta).
$$
So in view of \eqref{freq deloc error bound}, if $\eps$ is sufficiently small depending on $\eta$, we can apply
Lemma~\ref{lemma long time} to deduce the bound \eqref{freq deloc ls}.  This concludes the proof of Proposition
\ref{lemma freq loc imp spacetime bds}.
\end{proof}

Comparing \eqref{freq deloc ls} with \eqref{HUGE} gives the desired contradiction if $u$ satisfies the hypotheses of
Proposition \ref{lemma freq loc imp spacetime bds}.  We therefore expect $u$ to be localized in frequency for each $t$.
Indeed we have:

\begin{corollary}[Frequency localization of energy at each time]\label{lemma freq loc}
Let $u$ be a minimal energy blowup solution of \eqref{schrodinger
equation}.  Then, for each time $t \in I_*$ there exists a dyadic
frequency $N(t)\in 2^{\Z}$ such that for every $\eta_4 \leq \eta
\leq \eta_0$ we have small energy at frequencies $\ll N(t)$
\begin{equation}
\|P_{\leq c(\eta)N(t)}u(t)\|_{\ho_x} \leq \eta,
\end{equation}
small energy at frequencies $\gg N(t)$
\begin{equation}\label{energy hifreq freq loc}
\|P_{\geq C(\eta)N(t)}u(t)\|_{\ho_x} \leq \eta,
\end{equation}
and large energy at frequencies $\sim N(t)$
\begin{equation}
\|P_{c(\eta)N(t)<\cdot<C(\eta)N(t)}u(t)\|_{\ho_x} \sim 1,
\end{equation}
where the values of $0 < c(\eta) \ll 1 \ll C(\eta) < \infty$ depend on $\eta$.
\end{corollary}

\begin{proof}
For $t \in I_*$ define
$$
N(t) := \sup \{N \in 2^{\Z} : \|P_{\leq N} u(t) \|_{\ho_x} \leq \eta_0 \},
$$
which is clearly positive.  As $\|u\|_{L^\infty_t \ho_x} \sim 1$, $N(t)$ is also finite.  From the definition of $N(t)$
we have that
$$
\|P_{\leq 2N(t)} u(t)\|_{\ho_x} > \eta_0.
$$
Let $\eta_4 \leq \eta \leq \eta_0$. If $C(\eta) \gg 1$ then we must
have \eqref{energy hifreq freq loc}, since otherwise Proposition
\ref{lemma freq loc imp spacetime bds} would imply
$\|u\|_{L^{\frac{2(n+2)}{n-2}}_{t,x}(I_*\times\R^n)} \lesssim
C(\eta)$, which would contradict $u$ being a minimal energy blowup
solution.

Also, as by the definition of $N(t)$, $\|P_{\leq N(t)}u(t)\|_{\dot{H}^1_x}\leq \eta_0$,
$\|u\|_{L^\infty_t \ho_x} \sim 1$ and \eqref{energy hifreq freq loc} imply that
\begin{align}\label{energy medfreq freq loc eta0}
\|P_{N(t) < \cdot < C(\eta_0)N(t)} u(t) \|_{\ho_x} \sim 1
\end{align}
and therefore,
$$
\|P_{c(\eta)N(t)< \cdot < C(\eta)N(t)} u(t)\|_{\ho_x} \sim 1
$$
for all $\eta_4 \leq \eta \leq \eta_0$.  Thus, if $c(\eta) \ll 1$ then $\|P_{\leq c(\eta)N(t)} u(t)\|_{\ho_x} \leq \eta$
for all $\eta_4 \leq \eta \leq \eta_0$, since otherwise \eqref{energy medfreq freq loc eta0} and
Proposition~\ref{lemma freq loc imp spacetime bds} would again imply
$\|u\|_{L^{\frac{2(n+2)}{n-2}}_{t,x}(I_*\times\R^n)} \lesssim C(\eta)$.
\end{proof}

Having shown that a minimal energy blowup solution must be localized in frequency, we turn our attention to space.
In physical space, we will not need the full strength of a localization result. We will settle instead for a weaker
property concerning the spatial concentration of a minimal energy blowup solution.  Roughly, \emph{concentration} will mean
large at some point, while we reserve \emph{localization} to mean simultaneously concentrated and small at points far
from the concentration point.  To obtain the concentration result,
we use an idea of Bourgain (see \cite{borg:scatter}). We divide the interval $I_*$ into three consecutive subintervals
$I_* = I_{-} \cup I_0 \cup I_{+}$, each containing a third of the $L^{\frac{2(n+2)}{n-2}}_{t,x}$ mass of $u$:
$$
\int_I \int_{\R^n} |u(t,x)|^{\frac{2(n+2)}{n-2}}dx dt = \frac{1}{3} \int_{I_*}\int_{\R^n} |u(t,x)|^{\frac{2(n+2)}{n-2}} dx dt \ \ \text{ for } I = I_-,I_0,I_+.
$$
It is on the middle interval $I_0$ that we will show space concentration.  The first step is:

\begin{proposition}[Potential energy bounded from below]\label{lemma potential bdd below}
For any minimal energy blowup solution to \eqref{schrodinger equation} and all $t \in I_0$ we have
\begin{equation}\label{potential bdd below}
\|u(t)\|_{L^{\frac{2n}{n-2}}_x} \geq \eta_1.
\end{equation}
\end{proposition}

\begin{proof}
If the linear evolution of the solution does not concentrate at
some point in spacetime, then we can use the small data theory and iterate. So say the linear evolution concentrates at
some point $(t_1, x_1)$.  If the solution is small in $L^{\frac{2n}{n-2}}_x$ at time $t=t_0$, we show that $t_0$ must be
far from $t_1$.  We then remove the energy concentrating at $(t_1, x_1)$ and use induction on energy.

More formally, we will argue by contradiction.  Suppose there exists some time $t_0 \in I_0$ such that
\begin{equation}\label{contra for pot bdd below}
\|u(t_0)\|_{L^{\frac{2n}{n-2}}_x} < \eta_1.
\end{equation}
Using \eqref{scaling} we scale $N(t_0)=1$.  If the linear evolution $e^{i(t-t_0)\Delta}u(t_0)$ had small
$L^{\frac{2(n+2)}{n-2}}_{t,x}$-norm then, by perturbation theory (see Lemma~\ref{lemma long time}),
the nonlinear solution would have small $L^{\frac{2(n+2)}{n-2}}_{t,x}$-norm as well.  Hence, we may assume
$$\|e^{i(t-t_0)\Delta}u(t_0)\|_{L^{\frac{2(n+2)}{n-2}}_{t,x}(\rr)} \gtrsim 1.$$

On the other hand, Corollary \ref{lemma freq loc} implies that
$$
\| \pl u(t_0) \|_{\ho_x} + \| \ph u(t_0)\|_{\ho_x} \lesssim \eta_0,
$$
where $\pl = P_{<c(\eta_0)}$ and $\ph = P_{>C(\eta_0)}$.  Strichartz estimates yield
$$
\|e^{i(t-t_0)\Delta} \pl u(t_0)\|_{L^{\frac{2(n+2)}{n-2}}_{t,x}(\R\times\R^n)} + \|e^{i(t-t_0)\Delta} \ph
u(t_0)\|_{L^{\frac{2(n+2)}{n-2}}_{t,x}(\R\times\R^n)} \lesssim \eta_0.
$$
Thus,
$$
\|e^{i(t-t_0)\Delta} \pmed u(t_0)\|_{L^{\frac{2(n+2)}{n-2}}_{t,x}(\R\times\R^n)} \sim 1
$$
where $\pmed = 1 - \pl - \ph$.  However, $\pmed u(t_0)$ has bounded energy (by \eqref{kinetic energy bound}) and Fourier
support in $c(\eta_0) \lesssim |\xi| \lesssim C(\eta_0)$; an application of Strichartz  and
\eqref{mass high freq bound} yields
$$
\|e^{i(t-t_0)\Delta} \pmed u(t_0)\|_{L^{\frac{2(n+2)}{n}}_{t,x}} \lesssim \|\pmed u(t_0)\|_{L^2_x} \lesssim C(\eta_0).
$$
Combining these estimates with H\"older we get
$$
\|e^{i(t-t_0)\Delta} \pmed u(t_0)\|_{L^\infty_{t,x}} \gtrsim c(\eta_0).
$$
In particular, there exist a time $t_1 \in \R$ and a point $x_1 \in \R^n$ so that
\begin{align}\label{concentrated}
|e^{i(t_1-t_0)\Delta} (\pmed u(t_0))(x_1)| \gtrsim c(\eta_0).
\end{align}
We may perturb $t_1$ so that $t_1 \neq t_0$ and, by time reversal symmetry, we may take $t_1 < t_0$.  Let
$\delta_{x_1}$ be the Dirac mass at $x_1$.  Define $f(t_1) := \pmed \delta_{x_1}$ and for $t > t_1$ define
$f(t) := e^{i(t-t_1)\Delta} f(t_1)$.  One should think of $f(t_1)$ as basically $u$ at $(t_1,x_1)$. The point is then to
compare $u(t_0)$ to the linear evolution of $f(t_1)$ at time $t_0$. We will show that $f(t)$ is fast decaying in any
$L^p_x$-norm for $1\leq p\leq \infty$.

\begin{lemma}\label{lemma potential bdd below lemma}
For any $t \in \R$ and any $1 \leq p \leq \infty$ we have
$$\|f(t)\|_{L^p_x} \lesssim C(\eta_0) \langle t-t_1\rangle^{\frac{n}{p}-\frac{n}{2}}.$$
\end{lemma}

\begin{proof}
We may translate so that $t_1 = x_1 = 0$. By Bernstein and the unitarity of $e^{it\Delta}$, we get
$$
\|f(t)\|_{L^\infty_x} \lesssim C(\eta_0) \|f(t)\|_{L^2_x}
   = C(\eta_0) \|\pmed \delta_{x_1}\|_{L^2_x} \lesssim C(\eta_0).
$$
By \eqref{dispersive ineq} we also have
$$
\|f(t)\|_{L^\infty_x} \lesssim |t|^{-\frac{n}{2}}\|\pmed \delta_{x_1}\|_{L^1_x} \lesssim C(\eta_0) |t|^{-\frac{n}{2}}.
$$
Combining these two estimates, we obtain
$$
\|f(t)\|_{L^\infty_x} \lesssim C(\eta_0) \langle t\rangle^{-\frac{n}{2}}.
$$
This proves the lemma in the case $p = \infty$.

For other $p$'s we use \eqref{fourier rep} to write
$$
f(t,x) = \int_{\R^n} e^{2 \pi i (x \cdot \xi - 2 \pi t |\xi|^2)} \phi_{med} (\xi) d\xi
$$
where $\phi_{med}$ is the Fourier multiplier corresponding to $\pmed$.  For $|x| \gg 1 + |t|$, repeated integration
by parts shows $|f(t,x)| \lesssim |x|^N$ for any $N\leq 0$.  On $|x| \lesssim 1 + |t|$, one integrates using the above
$L_x^\infty$-bound.
\end{proof}

{}From \eqref{contra for pot bdd below} and H\"older we have
$$
| \langle u(t_0),f(t_0)\rangle| \lesssim \|f(t_0)\|_{L^{\frac{2n}{n+2}}_x}\|u(t_0)\|_{L^{\frac{2n}{n-2}}_x}
                                \lesssim \eta_1 C(\eta_0) \langle t_1 - t_0\rangle.
                                $$
On the other hand, by \eqref{concentrated} we get
$$
| \langle u(t_0),f(t_0) \rangle | = | \langle e^{i(t_1 - t_0)\Delta}\pmed u(t_0),\delta_{x_1}\rangle| \gtrsim c(\eta_0).
$$
So $\langle t_1 - t_0\rangle \gtrsim c(\eta_0)/\eta_1$, i.e.,
$t_1$ is far from $t_0$.  In particular, the time of concentration
must be far from where the $L^{\frac{2n}{n-2}}_x$-norm is small.

Also, from Lemma \ref{lemma potential bdd below lemma} we see
that $f$ has small $L^{\frac{2(n+2)}{n-2}}_tL_x^{\frac{2n(n+2)}{n^2+4}}$-norm to the
future of $t_0$ (recall $t_1 < t_0$):
\begin{align}\label{small future lp}
\|f\|_{L^{\frac{2(n+2)}{n-2}}_tL_x^{\frac{2n(n+2)}{n^2+4}}([t_0,\infty)\times \R^n)}
&\lesssim C(\eta_0) \|\langle \cdot - t_1\rangle^{-\frac{n-2}{n+2}}\|_{L_t^{\frac{2(n+2)}{n-2}}([t_0,\infty))}\notag\\
&\lesssim C(\eta_0)|t_1 - t_0|^{-\frac{n-2}{2(n+2)}}
\lesssim C(\eta_0)\eta_1^{\frac{n-2}{2(n+2)}}.
\end{align}

Now we use the induction hypothesis.  Split $u(t_0) = v(t_0) + w(t_0)$ where
$w(t_0) = \delta e^{i \theta} \Delta^{-1} f(t_0)$ for some small $\delta = \delta(\eta_0) > 0$ and phase $\theta$ to
be chosen later.\footnote{The presence of $\Delta^{-1}$ in the definition of $w(t_0)$ is due to the fact
that inner products are taken in $\dot{H}^1_x$.} One should think of $w(t_0)$ as the contribution coming from the point
$(t_1,x_1)$ where the solution concentrates. We will show that for an appropriate choice of $\delta$ and $\theta$,
$v(t_0)$ has slightly smaller energy than $u$.  By the definition of $f$ and an integration by parts we have
\begin{align*}
\frac{1}{2} \int_{\R^n} |\nabla v(t_0)|^2 dx &= \frac{1}{2} \int_{\R^n} |\nabla u(t_0) - \nabla w(t_0)|^2 dx\\
&=\frac{1}{2} \int_{\R^n} |\nabla u(t_0)|^2 dx - \delta \text{Re} \int_{\R^n} e^{-i \theta} \overline{\nabla \Delta^{-1}f(t_0)} \cdot \nabla u(t_0) dx\\
&\quad +O(\delta^2\|\Delta^{-1}f(t_0)\|_{\ho_x}^2)\\
& \leq \ecrit + \delta \text{Re } e^{-i\theta}\langle u(t_0),f(t_0) \rangle + O(\delta^2 C(\eta_0)).
\end{align*}
Choosing $\delta$ and $\theta$ appropriately we get
$$
\frac{1}{2} \int_{\R^n} |\nabla v(t_0)|^2 dx \leq \ecrit - c(\eta_0).
$$
Also, by Lemma \ref{lemma potential bdd below lemma} we have
$$
\|w(t_0)\|_{L^{\frac{2n}{n-2}}_x} \lesssim C(\eta_0)\|f(t_0)\|_{L^{\frac{2n}{n-2}}_x}
   \lesssim C(\eta_0)\langle t_1 - t_0\rangle^{-1}
   \lesssim C(\eta_0)\eta_1.
$$
So, by \eqref{contra for pot bdd below} and the triangle inequality we obtain
$$
\int_{\R^n} |v(t_0)|^{\frac{2n}{n-2}} dx \lesssim C(\eta_0)\eta_1^{\frac{2n}{n-2}}.
$$
Combining the above two energy estimates and taking $\eta_1$ sufficiently small depending on $\eta_0$, we obtain
$$
E(v(t_0))\leq \ecrit - c(\eta_0).
$$
Lemma \ref{lemma induct on energy} implies that there exists a global solution $v$ to \eqref{schrodinger equation}
with initial data $v(t_0)$ at time $t_0$ satisfying
$$
\|v\|_{\dot{S}^1(\R\times\R^n)}\lesssim C(\eta_0).
$$
In particular,
$$
\|v\|_{L_t^\infty \dot{H}^1_x([t_0, \infty) \times \R^n)} \lesssim C(\eta_0)
$$
and
$$
\|v\|_{L^{\frac{2(n+2)}{n-2}}_{t,x}([t_0, \infty) \times \R^n)} \lesssim C(\eta_0).
$$

Moreover, by Bernstein,
$$
\|w(t_0)\|_{\dot{H}^1_x}
\lesssim \delta \|\nabla \Delta^{-1} f(t_0)\|_{L_x^2} \lesssim C(\eta_0).
$$
By \eqref{small future lp} and frequency localization, we estimate
\begin{align*}
\sum_{N}\|P_N\nabla & e^{i(t-t_0)\Delta} w(t_0)\|_{L^{\frac{2(n+2)}{n-2}}_tL_x^{\frac{2n(n+2)}{n^2+4}}([t_0,\infty)\times \R^n)}^2\\
&\lesssim \sum_{N\leq C(\eta_0)}\|P_N\nabla e^{i(t-t_0)\Delta} w(t_0)\|_{L^{\frac{2(n+2)}{n-2}}_tL_x^{\frac{2n(n+2)}{n^2+4}}([t_0,\infty)\times \R^n)}^2\\
&\quad + \sum_{C(\eta_0)<N}\|P_N\nabla e^{i(t-t_0)\Delta} w(t_0)\|_{L^{\frac{2(n+2)}{n-2}}_tL_x^{\frac{2n(n+2)}{n^2+4}}([t_0,\infty)\times \R^n)}^2\\
&\lesssim \sum_{N\leq C(\eta_0)}N^2C(\eta_0)\|f\|_{L^{\frac{2(n+2)}{n-2}}_tL_x^{\frac{2n(n+2)}{n^2+4}}([t_0,\infty)\times \R^n)}^2\\
&\quad + \sum_{C(\eta_0)<N}N^{-2} \|f\|_{L^{\frac{2(n+2)}{n-2}}_tL_x^{\frac{2n(n+2)}{n^2+4}}([t_0,\infty)\times \R^n)}^2\\
&\lesssim C(\eta_0)\eta_1^{\frac{n-2}{n+2}}
\end{align*}
and hence
\begin{align*}
\Bigl(\sum_{N}\|P_N\nabla & e^{i(t-t_0)\Delta} w(t_0)\|_{L^{\frac{2(n+2)}{n-2}}_tL_x^{\frac{2n(n+2)}{n^2+4}}([t_0,\infty)\times \R^n)}^2\Bigr)^{\frac{1}{2}}
\lesssim C(\eta_0)\eta_1^{\frac{n-2}{2(n+2)}}.
\end{align*}

So, if $\eta_1$ is sufficiently small depending on $\eta_0$, we can apply Lemma \ref{lemma long time} with
$\tilde u = v$ and $e = 0$ to conclude that $u$ extends to all of $[t_0, \infty)$ and obeys
$$
\|u\|_{L^{\frac{2(n+2)}{n-2}}_{t,x}([t_0,\infty)\times\R^n)} \lesssim C(\eta_0,\eta_1).
$$
As $[t_0,\infty)$ contains $I_+$, the above estimate contradicts \eqref{HUGE} if $\eta_5$ is chosen sufficiently small.
This concludes the proof of Proposition \ref{lemma potential bdd below}.
\end{proof}

Using \eqref{potential bdd below} we can deduce the desired concentration result:

\begin{proposition}[Spatial concentration of energy at each time]\label{lemma physical concentration}
For any minimal energy blowup solution to \eqref{schrodinger equation} and for each $t \in I_0$, there exists
$x(t) \in \R^n$ such that
\begin{equation}\label{physical conc kinetic}
\int_{|x-x(t)| \leq C(\eta_1)/N(t)} |\nabla u(t,x)|^2 dx \gtrsim c(\eta_1)
\end{equation}
and
\begin{equation}\label{physical conc lp}
\int_{|x-x(t)| \leq C(\eta_1)/N(t)} |u(t,x)|^p dx \gtrsim c(\eta_1) N(t)^{\frac{n-2}{2}p-n}
\end{equation}
for all $1 < p < \infty$, where the implicit constants depend on $p$.
\end{proposition}

\begin{proof}
Fix $t$ and normalize $N(t)=1$.  By Corollary \ref{lemma freq loc} we have
$$
\|P_{<c(\eta_1)} u(t)\|_{\ho_x} + \|P_{>C(\eta_1)} u(t)\|_{\ho_x} \lesssim \eta_1^{100}  .
$$
Sobolev embedding implies
$$
\|P_{<c(\eta_1)} u(t)\|_{L^{\frac{2n}{n-2}}_x} +\|P_{>C(\eta_1)} u(t)\|_{L^{\frac{2n}{n-2}}_x}\lesssim \eta_1^{100}
$$
and so, by \eqref{potential bdd below},
\begin{align}\label{space med1}
\|\pmed u(t)\|_{L^{\frac{2n}{n-2}}_x} \gtrsim \eta_1,
\end{align}
where $\pmed = P_{c(\eta_1) \leq \cdot \leq C(\eta_1)}$.  On the other hand, by \eqref{mass high freq bound} we have
\begin{align}\label{space med2}
\|\pmed u(t)\|_{L^2_x} \lesssim C(\eta_1).
\end{align}
Thus, by H\"older, \eqref{space med1}, and \eqref{space med2}, we get
$$
\|\pmed u(t)\|_{L^\infty_x} \gtrsim c(\eta_1).
$$
In particular, there exists a point $x(t) \in \R^n$ so that
\begin{equation}\label{pmed bdd below}
c(\eta_1) \lesssim |\pmed u(t,x(t))|.
\end{equation}
As our function is now localized both in frequency and in space, all the Sobolev norms are practically equivalent.
So let's consider the operator $\pmed \nabla \Delta^{-1}$ and let $K_{med}$ denote its kernel.  Then,
\begin{align*}
c(\eta_1)
&\lesssim |\pmed u(t,x(t))| \lesssim |K_{med} * \nabla u(t,x(t))| \lesssim \int_{\R^n} |K_{med}(x(t)-x)||\nabla u(t,x)| dx\\
& \sim \int_{|x-x(t)|<C(\eta_1)}|K_{med}(x(t)-x)||\nabla u(t,x)| dx\\
&\phantom{\sim \int_{|x-x(t)|<C(\eta_1)}|K_{med}} + \int_{|x-x(t)|\geq C(\eta_1)}|K_{med}(x(t)-x)||\nabla u(t,x)| dx\\
& \lesssim C(\eta_1) \Bigl( \int_{|x-x(t)|< C(\eta_1)}|\nabla u(t,x)|^2 dx \Bigr)^{1/2} + \int_{|x-x(t)| \geq C(\eta_1)}
         \frac{|\nabla u(t,x)|}{|x-x(t)|^{100n}} dx ,
\end{align*}
where in order to obtain the last inequality we used Cauchy-Schwarz and that $K_{med}$ is a Schwartz function.
Therefore, by \eqref{kinetic energy bound} we have
$$
c(\eta_1) \lesssim \Bigl( \int_{|x-x(t)|< C(\eta_1)}|\nabla u(t,x)|^2 dx \Bigr)^{1/2} + C(\eta_1)^{-\alpha}
$$
for some $\alpha > 0$, proving \eqref{physical conc kinetic}.

Now let $\tilde K_{med}$ be the kernel associated to $\pmed$ and let $1 < p < \infty$.  As above, we get
\begin{align*}
c(\eta_1)
&\lesssim \int_{\R^n} |\tilde K_{med}(x(t)-x)|| u(t,x)| dx\\
& \sim \int_{|x-x(t)|<C(\eta_1)}|\tilde K_{med}(x(t)-x)|| u(t,x)| dx\\
&\phantom{\lesssim C(\eta_1) \int_{|x-x(t)|< C(\eta_1)}} + \int_{|x-x(t)|\geq C(\eta_1)}|\tilde K_{med}(x(t)-x)||u(t,x)| dx\\
& \lesssim C(\eta_1) \Bigl( \int_{|x-x(t)|< C(\eta_1)}| u(t,x)|^p dx \Bigr)^{1/p} \\
&\phantom{\lesssim C(\eta_1) \int_{|x-x(t)|< C(\eta_1)}} +\|u(t)\|_{L^{\frac{2n}{n-2}}_x} \Bigl( \int_{|x-x(t)| \geq C(\eta_1)} \frac{1}{|x-x(t)|^{100n \cdot \frac{2n}{n+2}}} dx \Bigr)^{\frac{n+2}{2n}}\\
& \lesssim C(\eta_1) \Bigl( \int_{|x-x(t)|< C(\eta_1)}| u(t,x)|^p dx \Bigr)^{1/p} + C(\eta_1)^{-\alpha}
\end{align*}
for some $\alpha > 0$, which, after scaling, proves \eqref{physical conc lp}.
\end{proof}

%
%
%
%

\section{Frequency-Localized Interaction Morawetz Inequality}

The goal of this section is to prove

\begin{proposition}[Frequency-localized interaction Morawetz estimate (FLIM)]\label{prop flim}\leavevmode\\
Assuming u is a minimal energy blowup solution to \eqref{schrodinger equation} and $N_{*}<c(\eta_{2})N_{min}$, we have
\begin{align}
\int_{I_{0}}\int_{\R^{n}}\int_{\R^{n}} &\frac{|P_{\geq N_{*}}u(t,y)|^{2}  |P_{\geq N_{*}}u(t,x)|^{2}}{|x-y|^{3}}dxdydt  \label{flim} \\
+ &\quad \int_{I_{0}}\int_{\R^{n}}\int_{\R^{n}}\frac{|P_{\geq N_{*}}u(t,y)|^{2}|P_{\geq N_{*}}u(t,x)|^{\frac{2n}{n-2}}}{|x-y|}dxdydt
     \lesssim\eta_{1}N_{*}^{-3}. \notag
\end{align}
Here, $N_{min}:=\inf_{t \in I_0} N(t)$.
\end{proposition}

\begin{remark}\label{rem Nmin}
$N_{min}>0$. Indeed, if $N_{min}=\inf_{t \in I_0} N(t)= 0$, there would exist a sequence $\{t_j\}_{j\in \N}\subset I_0$
such that $N(t_j)\to 0$ as $j\to \infty$. By passing, if necessary, to a subsequence, we may assume
$\{t_j\}_{j\in \N}$ converges to $t_\infty\in I_0$. By definition (see Corollary~\ref{lemma freq loc}),
$$
\|P_{\leq 2N(t_j)}u(t_j)\|_{\dot{H}^1_x}>\eta_0.
$$
{}From the triangle inequality, we get
\begin{align*}
\eta_0<\|P_{\leq 2N(t_j)}u(t_j)\|_{\dot{H}^1_x}
&\leq \|P_{\leq 2N(t_j)}(u(t_j)-u(t_\infty))\|_{\dot{H}^1_x}+\|P_{\leq 2N(t_j)}u(t_\infty)\|_{\dot{H}^1_x}\\
&\lesssim \|u(t_j)-u(t_\infty)\|_{\dot{H}^1_x}+\|P_{\leq 2N(t_j)}u(t_\infty)\|_{\dot{H}^1_x}.
\end{align*}
As $u\in C_t^0\dot{H}^1_x(I_0\times\R^n)$, a limiting argument combined with the dominated convergence theorem leads
to a contradiction.
\end{remark}
As the right-hand side in \eqref{flim} does not depend on $|I_0|$, this estimate excludes the formation of solitons,
at least for frequencies `close' to $N_{min}$ and provided $|I_0|$ is taken sufficiently large. Frequencies much
larger than $N_{min}$ will be dealt with in Section~6.

In what follows, we will only prove Proposition~\ref{prop flim} in dimensions $n\geq 6$. More precisely, the subsections
5.3 and 5.4 below are adapted to treat the case $n\geq 6$ only. For a proof of Proposition~\ref{prop flim} in dimension
$n=5$ see \cite{my thesis}.

\subsection{An interaction virial identity and a general interaction Morawetz estimate}
The calculations in this subsection are difficult to justify without additional assumptions on the solution. This
obstacle can be dealt with in the standard manner: mollify the initial data and the nonlinearity to make the interim
calculations valid and observe that the mollifications can be removed at the end. For expository reasons, we skip
the details and keep all computations on a formal level.

We start by recalling the standard Morawetz action centered at a point. Let $a$  be a function on the slab
$I\times \R^{n}$ and $\phi$ satisfying
\begin{align}\label{phi equation}
i\phi_{t}+\Delta \phi=\mathcal{N}
\end{align}
on $I\times \R^{n}$. We define the Morawetz action centered at zero to be
$$
M_{a}^0(t)=2\int_{\R^{n}}a_{j}(x)\Im(\overline{\phi(x)}\phi_{j}(x))dx.
$$

A calculation establishes
\begin{lemma}
$$
\partial_{t}M_{a}^0=\int_{\R^n} (-\Delta\Delta a)|\phi|^{2}
  +4\int_{\R^n} a_{jk}\Re(\overline{\phi_{j}}\phi_{k})+2\int_{\R^n} a_{j}\{\mathcal{N},\phi\}_{p}^{j},
$$
where we define the momentum bracket to be $\{f, g\}_p=\Re(f\nabla \bar{g}- g\nabla \bar{f})$ and repeated indices are
implicitly summed.
\end{lemma}

Note that when $\mathcal{N}$ is the energy-critical nonlinearity in dimension $n$, we have
$\{\mathcal{N},\phi\}_p=-\frac{2}{n}\nabla(|\phi|^{\frac{2n}{n-2}})$.

Now let $a(x)=|x|$. For this choice of the function $a$, one should interpret $M_a^0$ as a spatial average of the
radial component of the $L^2_x$-mass current. Easy computations show that in dimension $n\geq 4$ we have the
following identities:
\begin{align*}
a_{j}(x)=&\frac{x_{j}}{|x|} \\
a_{jk}(x)=&\frac{\delta_{jk}}{|x|}-\frac{x_{j}x_{k}}{|x|^{3}} \\
\Delta a(x)=&\frac{n-1}{|x|} \\
-\Delta \Delta a(x)=&\frac{(n-1)(n-3)}{|x|^{3}}
\end{align*}
and hence,
\begin{align*}
\partial_{t}M_{a}^0
&=(n-1)(n-3)\int_{\R^n} \frac{|\phi(x)|^{2}}{|x|^{3}}dx
   +4\int_{\R^n} \Bigl(\frac{\delta_{jk}}{|x|}-\frac{x_{j}x_{k}}{|x|^{3}}\Bigr) \Re(\overline{\phi_{j}}\phi_{k})(x)dx\\
&\phantom{=((n-1)(n-3)\int_{\R^n} \frac{|\phi(x)|^{2}}{|x|^{3}}dx} +2\int_{\R^n} \frac{x_{j}}{|x|}\{\mathcal{N},\phi\}_{p}^{j}(x)dx \\
&=(n-1)(n-3)\int_{\R^n} \frac{|\phi(x)|^{2}}{|x|^{3}}dx
   +4\int_{\R^n} \frac{1}{|x|} |\nabla_{0}\phi(x)|^{2}dx\\
&\phantom{=((n-1)(n-3)\int_{\R^n} \frac{|\phi(x)|^{2}}{|x|^{3}}dx} +2\int_{\R^n} \frac{x}{|x|} \{\mathcal{N},\phi\}_{p}(x)dx,
\end{align*}
where we use $\nabla_{0}$ to denote the complement of the radial portion of the gradient, that is,
$\nabla_{0}=\nabla-\frac{x}{|x|}\bigl(\frac{x}{|x|}\cdot\nabla\bigr)$.

We may center the above argument at any other point $y\in \R^{n}$. Choosing $a(x)=|x-y|$,
we define the Morawetz action centered at $y$ to be
$$
M_{a}^y(t)=2\int_{\R^{n}}\frac{x-y}{|x-y|}\Im(\overline{\phi(x)}\nabla
\phi(x))dx.
$$
The same considerations now yield
\begin{align*}
\partial_{t}M_a^{y}
&=(n-1)(n-3)\int_{\R^n} \frac{|\phi(x)|^{2}}{|x-y|^{3}}dx
   +4\int_{\R^n} \frac{1}{|x-y|} |\nabla_{y}\phi(x)|^{2}dx\\
&\phantom{=((n-1)(n-3)\int_{\R^n} \frac{|\phi(x)|^{2}}{|x-y|^{3}}dx}
+2\int_{\R^n} \frac{x-y}{|x-y|} \{\mathcal{N},\phi\}_{p}(x)dx.
\end{align*}

We are now ready to define the interaction Morawetz potential, which is a way of quantifying how mass is interacting
with (moving away from) itself:
\begin{align*}
M^{interact}(t)
&=\int_{\R^n} |\phi(t,y)|^{2}M_a^{y}(t)dy\\
&=2\Im \int_{\R^n} \int_{\R^n} |\phi(t,y)|^{2}\frac{x-y}{|x-y|}\nabla\phi(t,x)\overline{\phi(t,x)}dxdy.
\end{align*}
One gets immediately the easy estimate
$$
|M^{interact}(t)| \leq 2 \|\phi(t)\|_{L_{x}^{2}}^{3} \|\phi(t)\|_{\dot{H}_{x}^{1}}.
$$

Calculating the time derivative of the interaction Morawetz potential, we get the following virial-type identity:
\begin{align}
\partial_{t}M^{interact}
=&(n-1)(n-3)\int_{\R^n} \int_{\R^n} \frac{|\phi(y)|^{2}|\phi(x)|^{2}}{|x-y|^{3}}dxdy  \label{vi1} \\
 &+4\int_{\R^n} \int_{\R^n} \frac{|\phi(y)|^{2}|\nabla_{y}\phi(x)|^{2}}{|x-y|}dxdy  \label{vi2} \\
 &+2\int_{\R^n} \int_{\R^n} |\phi(y)|^{2}\frac{x-y}{|x-y|} \{\mathcal{N},\phi\}_{p}(x)dxdy  \label{vi3} \\
 &+2 \int_{\R^n} \partial_{y_{k}} \Im(\phi\overline{\phi_{k}})(y)M_a^{y}dy  \label{vi4} \\
 &+4\Im \int_{\R^n} \int_{\R^n} \{\mathcal{N},\phi\}_{m}(y)\frac{x-y}{|x-y|}\nabla\phi(x)\overline{\phi(x)}dxdy, \label{vi5}
\end{align}
where the mass bracket is defined to be $\{f,g\}_m=\Im(f\bar{g})$.\\

As far as the terms in the above identity are concerned, at the end of the subsection we will establish
\begin{lemma}\label{lemmaviterms}
\eqref{vi4} ${} + {}$ \eqref{vi2} ${} \geq {} 0$.
\end{lemma}

Thus, integrating over the compact interval $I_{0}$ we get:

\begin{proposition}[Interaction Morawetz inequality] \label{intmorineq}
\begin{align*}
&(n-1)(n-3)\int_{I_{0}} \int_{\R^n} \int_{\R^n} \frac{|\phi(t,y)|^{2}|\phi(t,x)|^{2}}{|x-y|^{3}}dxdydt \\
&\phantom{(n-1)(n-3)\int_{I_{0}}} +2\int_{I_{0}} \int_{\R^n} \int_{\R^n} |\phi(t,y)|^{2}\frac{x-y}{|x-y|}\{\mathcal{N},\phi\}_{p}(t,x)dxdydt \\
&\phantom{(n-1)(n)} \leq 4\|\phi\|_{L_{t}^{\infty}L_{x}^{2}(I_{0}\times\R^{n})}^{3} \|\phi\|_{L_{t}^{\infty}\dot{H}_{x}^{1}(I_{0}\times\R^{n})} \\
&\phantom{(n-1)(n-3)\int_{I_{0}}}+4\int_{I_{0}}\int_{\R^n}\int_{\R^n}|\{\mathcal{N},\phi\}_{m}(t,y)||\nabla\phi(t,x)||\phi(t,x)|dxdydt.
\end{align*}
\end{proposition}

Note that in the particular case $\mathcal{N}=|u|^{\frac{4}{n-2}}u$, after performing an integration by parts in the
momentum bracket term, the inequality becomes
\begin{align}\label{int mor est}
(n-1)(n-3)\int_{I_{0}} \int_{\R^n} &\int_{\R^n}\frac{|u(t,y)|^{2}|u(t,x)|^{2}}{|x-y|^{3}}dxdydt \notag \\
&\quad+\frac{4(n-1)}{n}\int_{I_{0}} \int_{\R^n} \int_{\R^n} \frac{|u(t,y)|^{2}|u(t,x)|^{\frac{2n}{n-2}}}{|x-y|}dxdydt\\
\leq &4\|u\|_{L_{t}^{\infty}L_{x}^{2}(I_{0}\times\R^{n})}^{3}\|u\|_{L_{t}^{\infty}\dot{H}_{x}^{1}(I_{0}\times\R^{n})}.\notag
\end{align}

Assuming $u$ has finite mass, this estimate is an expression of dispersion (as the interaction between the masses of
two particles is weak) and local smoothing (as it implies $|u|^2\in L_t^2\dot{H}^{-\frac{n-3}{2}}_x$). However,
we have made no assumptions regarding the finiteness of the $L_x^2$-norm of the initial data $u_0$ and thus
\eqref{int mor est} cannot be used directly.

We turn now to the proof of Lemma \ref{lemmaviterms}. We write
$$
\eqref{vi4}=4\int_{\R^n} \int_{\R^n} \partial_{y_{k}}
\Im(\phi(y)\overline{\phi_{k}(y)})\frac{x_j-y_j}{|x-y|}\Im(\overline{\phi(x)}\phi_j(x))dxdy,
$$
where repeated indices are implicitly summed. We integrate by parts moving
$\partial_{y_k}$ to the unit vector $\frac{x-y}{|x-y|}$. Using the
identity
$$
\partial_{y_k}\Bigl(\frac{x_j-y_j}{|x-y|}\Bigr)=-\frac{\delta_{kj}}{|x-y|}+\frac{(x_k-y_k)(x_j-y_j)}{|x-y|^3}
$$
and the notation $p(x)=2\Im(\overline{\phi(x)}\nabla \phi(x))$ for
the momentum density, we rewrite \eqref{vi4} as
$$
-\int_{\R^n} \int_{\R^n} \Bigl[ p(y)p(x)-\Bigl(
p(y)\frac{x-y}{|x-y|} \Bigr) \Bigl( p(x)\frac{x-y}{|x-y|} \Bigr)
\Bigr] \frac{dxdy}{|x-y|}.
$$
In the quantity between the square brackets we recognize the inner product between the projections of the
momentum densities $p(x)$ and $p(y)$ onto the orthogonal complement of $(x-y)$. As
\begin{align*}
|\pi_{(x-y)^\perp}p(y)|
   &=\Bigl|p(y)-\frac{x-y}{|x-y|}\Bigl(\frac{x-y}{|x-y|}p(y)\Bigr)\Bigr|
    =2|\Im(\overline{\phi(y)}\nabla_x \phi(y))| \\
   &\leq 2|\phi(y)||\nabla_x \phi(y))|
\end{align*}
and the same estimate holds when we switch $y$ and $x$, we get
\begin{align*}
\eqref{vi4}
    &\geq -4\int_{\R^n} \int_{\R^n} |\phi(y)||\nabla_x \phi(y))| |\phi(x)| |\nabla_y \phi(x))|\frac{dxdy}{|x-y|} \\
    &\geq -2 \int_{\R^n} \int_{\R^n} \frac{|\phi(y)|^{2}|\nabla_{y}\phi(x)|^{2}}{|x-y|}dxdy
       -2 \int_{\R^n} \int_{\R^n} \frac{|\phi(x)|^{2}|\nabla_{x}\phi(y)|^{2}}{|x-y|}dxdy \\
    &\geq -\eqref{vi2}.
\end{align*}

\subsection{FLIM: the setup} We are now ready to start the proof of Proposition \ref{prop flim}. As the
statement is invariant under scaling, we normalize $N_{*}=1$ and define $u_{hi}=P_{> 1}u$ and $u_{lo}=P_{\leq 1}u$. As
we assume $1=N_{*}<c(\eta_{2})N_{min}$, we have $1<c(\eta_{2})N(t)$, $\forall t\in I_{0}$. Choosing
$c(\eta_{2})$ sufficiently small (smaller than $\eta_2 \tilde{c}(\eta_2)$ where $\tilde{c}(\eta_2)$ is the constant
appearing in Corollary~\ref{lemma freq loc}), the frequency localization result and Sobolev embedding yield
\begin{equation}\label{sijz}
\|u_{<\eta_2^{-1}}\|_{L_{t}^{\infty}\dot{H}_{x}^{1}(I_{0}\times\R^{n})}
  +\|u_{<\eta_2^{-1}}\|_{L_{t}^{\infty}L_{x}^{\frac{2n}{n-2}}(I_{0}\times\R^{n})}
  \lesssim \eta_{2}.
\end{equation}
In particular, this implies that $u_{lo}$ has small energy
\begin{align}
\|u_{lo}\|_{L_{t}^{\infty}\dot{H}_{x}^{1}(I_{0}\times\R^{n})}
  +\|u_{lo}\|_{L_{t}^{\infty}L_{x}^{\frac{2n}{n-2}}(I_{0}\times\R^{n})}\lesssim
\eta_{2}. \label{lowfreqsmall}
\end{align}
Using \eqref{mass high freq bound} and \eqref{sijz}, one also sees that $u_{hi}$ has small mass
\begin{align}
\|u_{hi}\|_{L_{t}^{\infty}L_{x}^{2}(I_{0}\times\R^{n})}\lesssim \eta_{2}. \label{smallmass}
\end{align}

Our goal is to prove \eqref{flim}, which, in particular, implies
\begin{align}\label{show}
\int_{I_{0}}\int_{\R^{n}}\int_{\R^{n}}\frac{|u_{hi}(t,x)|^{2}|u_{hi}(t,y)|^{2}}{|x-y|^{3}}dx dy dt\lesssim\eta_{1}.
\end{align}
As in dimension $n$ convolution with $1/|x|^3$ is basically the same as the
fractional integration operator $|\nabla|^{-(n-3)}$, the above estimate translates into
\begin{align}
\| |u_{hi}|^2 \|_{L^2_t \dot H^{-\frac{n-3}{2}}_x(I_{0} \times \R^n)} \lesssim \eta^{1/2}_1. \label{hifreqasmp}
\end{align}

By a standard continuity argument, it suffices to prove \eqref{hifreqasmp} under the bootstrap hypothesis
\begin{align}
\| |u_{hi}|^2 \|_{L^2_t \dot H^{-\frac{n-3}{2}}_x(I_{0} \times \R^n)} \leq (C_0 \eta_1)^{\frac{1}{2}}, \label{bootstraphyp1}
\end{align}
for a large constant $C_0$ depending on energy but not on any of the $\eta$'s. In fact, we need to prove that
\eqref{bootstraphyp1} implies \eqref{show} whenever $I_0$ is replaced by a subinterval of $I_0$ in order to run the
continuity argument correctly. However, it will become clear to the reader that the argument below works not only
for $I_0$, but also for any of its subintervals.

First, let us note that \eqref{bootstraphyp1} implies
\begin{align}
\||\nabla|^{-\frac{n-3}{4}}u_{hi} \|_{L^4_{t,x}(I_0\times \R^n)}\lesssim (C_0\eta_1)^{\frac{1}{4}}, \label{bootstraphyp2}
\end{align}
as can be seen by taking $f=u_{hi}$ in the following

\begin{lemma}
\begin{align}\label{uhi in L^4}
\||\nabla|^{-\frac{n-3}{4}}f\|_4\lesssim \||\nabla|^{-\frac{n-3}{2}}|f|^2\|_2^{1/2}.
\end{align}
\end{lemma}

\begin{proof}
As $|\nabla|^{-\frac{n-3}{4}}$ and $|\nabla|^{-\frac{n-3}{2}}$
correspond to convolutions with positive kernels, it suffices to
prove \eqref{uhi in L^4} for a positive Schwartz function $f$. For
such an $f$, we will show the pointwise inequality
\begin{align}\label{square function}
S(|\nabla|^{-\frac{n-3}{4}}f)(x)\lesssim
[(|\nabla|^{-\frac{n-3}{2}}|f|^2)(x)]^{1/2},
\end{align}
where $S$ denotes the Littlewood-Paley square function $Sf :=
(\sum_N |P_N f|^2)^{1/2}$. Clearly \eqref{square function} implies
\eqref{uhi in L^4}:
$$
\||\nabla|^{-\frac{n-3}{4}}f\|_4\lesssim
\|S(|\nabla|^{-\frac{n-3}{4}}f)\|_4
    \lesssim\|(|\nabla|^{-\frac{n-3}{2}}|f|^2)^{1/2}\|_4\lesssim\||\nabla|^{-\frac{n-3}{2}}|f|^2\|_2^{1/2}.
$$

In order to prove \eqref{square function} we will estimate each of
the dyadic pieces,
$$
P_N(|\nabla|^{-\frac{n-3}{4}}f)(x)=\int e^{2\pi i
x\xi}\hat{f}(\xi)|\xi|^{-\frac{n-3}{4}}m(\xi/N)d\xi,
$$
where $m(\xi):=\phi(\xi)-\phi(2\xi)$ in the notation introduced in
Section~2. As $|\xi|^{-\frac{n-3}{4}}m(\xi/N)\sim
N^{-\frac{n-3}{4}}\tilde{m}(\xi/N)$ for $\tilde{m}$ a multiplier
with the same properties as $m$, we have
\begin{align*}
P_N(|\nabla|^{-\frac{n-3}{4}}f)(x)
 &\sim f* \bigl( N^{-\frac{n-3}{4}} [\tilde{m}(\xi/N)]\hbox{\Large$\check{\ }$}(x) \bigr)
    = N^{\frac{3(n+1)}{4}} f*\check{\tilde{m}}(Nx)  \\
 &= N^{\frac{3(n+1)}{4}} \int f(x-y) \check{\tilde{m}}(Ny) dy.
\end{align*}
An application of Cauchy-Schwartz yields
\begin{align*}
S(|\nabla|^{-\frac{n-3}{4}}f)(x)
  &=\Bigl(\sum_{N}|P_N(|\nabla|^{-\frac{n-3}{4}}f)(x)|^2\Bigr)^{1/2} \\
  &\lesssim \Bigl(\sum_{N} N^{\frac{3(n+1)}{2}} \Bigl|\int f(x-y)\check{\tilde{m}}(Ny) dy\Bigr|^2\Bigr)^{1/2} \\
  &\lesssim \Bigl(\sum_{N} N^{\frac{3(n+1)}{2}} \int |\check{\tilde{m}}(Ny)| dy \int |f(x-y)|^2|\check{\tilde{m}}(Ny)|dy\Bigr)^{1/2}\\
  &\lesssim \Bigl(\sum_{N} N^{\frac{n+3}{2}}\int |f(x-y)|^2|\check{\tilde{m}}(Ny)| dy \Bigr)^{1/2}.
\end{align*}
As $\check{\tilde{m}}$ is rapidly decreasing,
$$
\sum_{N} N^{\frac{n+3}{2}}|\check{\tilde{m}}(Ny)|
   \lesssim \sum_{N} N^{\frac{n+3}{2}}\min\{1, |Ny|^{-100n}\}
   \lesssim |y|^{-\frac{n+3}{2}}.
$$
In this way we get
$$
S(|\nabla|^{-\frac{n-3}{4}}f)(x)\lesssim \Bigl( \int
\frac{|f(x-y)|^2}{|y|^{\frac{n+3}{2}}} dy\Bigr)^{1/2}
    \sim [(|\nabla|^{-\frac{n-3}{2}}|f|^2)(x)]^{1/2},
$$
and the claim follows.
\end{proof}

We now use Proposition \ref{intmorineq} to derive an interaction Morawetz estimate for $\phi:=u_{hi}$.

\begin{proposition} \label{etprop}
With the notation and assumptions above we have
\begin{align*}
\int_{I_{0}}\int_{\R^n}\int_{\R^n} \frac{|u_{hi}(t,y)|^{2}|u_{hi}(t,x)|^{2}}{|x-y|^{3}}&dxdydt\\
   &+\int_{I_{0}}\int_{\R^n}\int_{\R^n}\frac{|u_{hi}(t,y)|^{2}|u_{hi}(t,x)|^{\frac{2n}{n-2}}}{|x-y|}dxdydt
\end{align*}
\begin{align}
&\lesssim  \eta_2^3 \label{et0}\\
&+\eta_{2} \int_{I_0} \int_{\R^n}|u_{hi}(t,x)||P_{hi}\bigl(|u|^{\frac{4}{n-2}}u-|u_{hi}|^{\frac{4}{n-2}}u_{hi}-|u_{lo}|^{\frac{4}{n-2}}u_{lo}\bigr)(t,x)|dxdt \label{et1}\\
&+\eta_{2} \int_{I_0} \int_{\R^n}|u_{hi}(t,x)||P_{lo}\bigl(|u_{hi}|^{\frac{4}{n-2}}u_{hi}\bigr)(t,x)| dxdt \label{et2}\\
&+\eta_{2} \int_{I_0} \int_{\R^n}|u_{hi}(t,x)||P_{hi}\bigl(|u_{lo}|^{\frac{4}{n-2}}u_{lo}\bigr)(t,x)| dxdt \label{et3}\\
&+\eta_{2}^2 \int_{I_0} \int_{\R^n}|\nabla u_{lo}(t,x)||u_{lo}(t,x)|^{\frac{4}{n-2}}|u_{hi}(t,x)|dxdt \label{et4}\\
&+\eta_{2}^2 \int_{I_0} \int_{\R^n}|\nabla u_{lo}(t,x)||u_{hi}(t,x)|^{\frac{n+2}{n-2}} dxdt \label{et5}\\
&+\eta_{2}^2 \int_{I_0} \int_{\R^n}|\nabla P_{lo}\bigl(|u|^{\frac{4}{n-2}}u\bigr)(t,x)||u_{hi}(t,x)|dxdt  \label{et6}\\
&+\int_{I_0}\int_{\R^n}\int_{\R^n}\frac{|u_{hi}(t,y)|^{2}|u_{lo}(t,x)|^{\frac{n+2}{n-2}}|u_{hi}(t,x)|}{|x-y|}dxdydt \label{et7}\\
&+\int_{I_0}\int_{\R^n}\int_{\R^n}\frac{|u_{hi}(t,y)|^{2}|u_{lo}(t,x)||u_{hi}(t,x)|^{\frac{n+2}{n-2}}}{|x-y|}dxdydt \label{et8}\\
&+\int_{I_0}\int_{\R^n}\int_{\R^n}\frac{|u_{hi}(t,y)|^{2}|P_{lo}\bigl(|u_{hi}|^{\frac{4}{n-2}}u_{hi}\bigr)(t,x)||u_{hi}(t,x)|}{|x-y|}dxdydt \label{et9}.
\end{align}
\end{proposition}

\begin{proof}
Applying Proposition \ref{intmorineq} with $\phi=u_{hi}$ and $\mathcal{N}=P_{hi}(|u|^{\frac{4}{n-2}}u)$, we find
\begin{align*}
&(n-1)(n-3)\int_{I_{0}} \int_{\R^n} \int_{\R^n}  \frac{|u_{hi}(t,y)|^{2}|u_{hi}(t,x)|^{2}}{|x-y|^{3}}dxdydt \\
&\phantom{(n-1)(n}+2\int_{I_{0}} \int_{\R^n} \int_{\R^n} |u_{hi}(t,y)|^{2}\frac{x-y}{|x-y|}\{P_{hi}(|u|^{\frac{4}{n-2}}u),u_{hi}\}_{p}(t,x)dxdydt \\
&\phantom{(n)}\leq 4\|u_{hi}\|_{L_{t}^{\infty}L_{x}^{2}(I_{0}\times\R^{n})}^{3} \|u_{hi}\|_{L_{t}^{\infty}\dot{H}_{x}^{1}(I_{0}\times\R^{n})} \\
&\phantom{(n-1)(n}+4 \int_{I_{0}} \int_{\R^n} \int_{\R^n} |\{P_{hi}(|u|^{\frac{4}{n-2}}u),u_{hi}\}_{m}(t,y)| |\nabla u_{hi}(t,x)||u_{hi}(t,x)|dxdydt.
\end{align*}

Observe that \eqref{smallmass} plus conservation of energy dictates
$$
\|u_{hi}\|_{L_{t}^{\infty}L_{x}^{2}(I_{0}\times\R^{n})}^{3}\|u_{hi}\|_{L_{t}^{\infty}\dot{H}_{x}^{1}(I_{0}\times\R^{n})}
   \lesssim \eta^3_2,
$$
which is the error term \eqref{et0}.

We consider the mass bracket term first. Exploiting cancellation, we write
\begin{align*}
\{P_{hi}\bigl(|u|^{\frac{4}{n-2}}u\bigr),u_{hi}\}_{m}
  &=\{P_{hi}\bigl(|u|^{\frac{4}{n-2}}u\bigr)-|u_{hi}|^{\frac{4}{n-2}}u_{hi},u_{hi}\}_{m} \\
  &=\{P_{hi}\bigl(|u|^{\frac{4}{n-2}}u-|u_{hi}|^{\frac{4}{n-2}}u_{hi}-|u_{lo}|^{\frac{4}{n-2}}u_{lo}\bigr),u_{hi}\}_{m} \\
  &\quad -\{P_{lo}\bigl(|u_{hi}|^{\frac{4}{n-2}}u_{hi}\bigr),u_{hi}\}_{m}+\{P_{hi}\bigl(|u_{lo}|^{\frac{4}{n-2}}u_{lo}\bigr),u_{hi}\}_{m}.
\end{align*}
Estimating
$$
\int_{\R^n} |u_{hi}(t,x)||\nabla u_{hi}(t,x)|dx\lesssim \|u_{hi}\|_{L_{t}^{\infty}L_{x}^{2}(I_{0}\times\R^{n})}
   \|\nabla u_{hi}\|_{L_{t}^{\infty}L_{x}^{2}(I_{0}\times\R^{n})} \lesssim \eta_2,
$$
we see that we can bound the contribution of the mass bracket term by the sum of \eqref{et1}, \eqref{et2}, and \eqref{et3}.

We turn now towards the momentum bracket term and write
\begin{align*}
\{P_{hi}\bigl(|u|^{\frac{4}{n-2}}u\bigr),u_{hi}\}_{p}
 &=\{|u|^{\frac{4}{n-2}}u,u_{hi}\}_{p}-\{P_{lo}\bigl(|u|^{\frac{4}{n-2}}u\bigr),u_{hi}\}_{p} \\
 &=\{|u|^{\frac{4}{n-2}}u,u\}_{p}-\{|u|^{\frac{4}{n-2}}u,u_{lo}\}_{p}-\{P_{lo}\bigl(|u|^{\frac{4}{n-2}}u\bigr),u_{hi}\}_{p} \\
 &=\{|u|^{\frac{4}{n-2}}u,u\}_{p}-\{|u_{lo}|^{\frac{4}{n-2}}u_{lo},u_{lo}\}_{p}- \{P_{lo}\bigl(|u|^{\frac{4}{n-2}}u\bigr),u_{hi}\}_{p}\\
 &\phantom{\{|u|^{\frac{4}{n-2}}u,u\}_{p}}-\{|u|^{\frac{4}{n-2}}u-|u_{lo}|^{\frac{4}{n-2}}u_{lo},u_{lo}\}_{p} \\
 &=-\frac{2}{n}\nabla(|u|^{\frac{2n}{n-2}}-|u_{lo}|^{\frac{2n}{n-2}})-\{|u|^{\frac{4}{n-2}}u-|u_{lo}|^{\frac{4}{n-2}}u_{lo},u_{lo}\}_{p} \\
 &\phantom{\{|u|^{\frac{4}{n-2}}u,u\}_{p}}-\{P_{lo}\bigl(|u|^{\frac{4}{n-2}}u\bigr),u_{hi}\}_{p} \\
 &=I + II + III.
\end{align*}

To estimate the contribution coming from $I$, we integrate by parts in the momentum bracket term; we obtain, up to
a constant,
\begin{align*}
&\int_{I_{0}} \int_{\R^n} \int_{\R^n} \frac{|u_{hi}(t,y)|^{2}}{|x-y|}\bigl(|u(t,x)|^{\frac{2n}{n-2}}-|u_{lo}(t,x)|^{\frac{2n}{n-2}}\bigr)dxdydt \\
&\quad =\int_{I_{0}} \int_{\R^n} \int_{\R^n} \frac{|u_{hi}(t,y)|^{2}|u_{hi}(t,x)|^{\frac{2n}{n-2}}}{|x-y|}dxdydt \\
&\quad \quad +\int_{I_{0}} \int_{\R^n} \int_{\R^n} \frac{|u_{hi}(t,y)|^{2}\bigl(|u(t,x)|^{\frac{2n}{n-2}}-|u_{lo}(t,x)|^{\frac{2n}{n-2}}-|u_{hi}(t,x)|^{\frac{2n}{n-2}}\bigr)}{|x-y|}dxdydt.
\end{align*}
In the above expression we recognize the left-hand side term in Proposition \ref{etprop} and an error that we estimate by
the sum of \eqref{et7} and \eqref{et8}.

In order to estimate the contribution of $II$, we write $\{f,g\}_{p}=\nabla \O(fg)+\O(f\nabla g)$ and hence,
\begin{align}
&\{|u|^{\frac{4}{n-2}}u-|u_{lo}|^{\frac{4}{n-2}}u_{lo},u_{lo}\}_p  \notag\\
&\quad = \nabla\O\bigl[\bigl(|u|^{\frac{4}{n-2}}u-|u_{lo}|^{\frac{4}{n-2}}u_{lo}\bigr)u_{lo}\bigr]\label{II1} \\
&\quad \quad + \O\bigl[\bigl(|u|^{\frac{4}{n-2}}u-|u_{lo}|^{\frac{4}{n-2}}u_{lo}\bigr)\nabla u_{lo}\bigr]\label{II2}.
\end{align}

Integrating by parts, we estimate the error coming form \eqref{II1} by a scalar multiple of
\begin{align*}
\int_{I_{0}}\int_{\R^n}\int_{\R^n}&\frac{|u_{hi}(t,y)|^{2}\bigl||u|^{\frac{4}{n-2}}u-|u_{lo}|^{\frac{4}{n-2}}u_{lo}\bigr|(t,x)|u_{lo}(t,x)|}{|x-y|}dxdydt\\
& \lesssim \eqref{et7} + \eqref{et8},
\end{align*}
where in order to obtain the last inequality we used \eqref{diff1}.

We turn now to the contribution of \eqref{II2}. Let us first note that
\begin{align}\label{uhi^2}
\|u_{hi}^2\|_{L^{\infty}_tL_x^1(I_{0} \times \R^n)}\lesssim \|u_{hi}\|_{L^{\infty}_tL_x^2(I_{0} \times \R^n)}^2
  \lesssim \eta_2^2.
\end{align}
Taking the absolute values inside the integrals and using \eqref{diff1} and \eqref{uhi^2}, we estimate the error
coming from \eqref{II2} by
\begin{align*}
\int_{I_{0}}\int_{\R^n}\int_{\R^n} & |u_{hi}(t,y)|^{2}\bigl||u|^{\frac{4}{n-2}}u-|u_{lo}|^{\frac{4}{n-2}}u_{lo}\bigr|(t,x)|\nabla u_{lo}(t,x)|dxdydt \\
& \lesssim \eqref{et4} + \eqref{et5}.
\end{align*}

We consider next the contribution of $III$ to the momentum bracket term. When the derivative (from the definition of
the momentum bracket) falls on $P_{lo}(|u|^{\frac{4}{n-2}}u)$, we take the absolute values inside the integrals and use
\eqref{uhi^2} to estimate this contribution by
\begin{align*}
\int_{I_{0}} \int_{\R^n} \int_{\R^n} &|u_{hi}(t,y)|^{2}|\nabla P_{lo}\bigl(|u|^{\frac{4}{n-2}}u\bigr)(t,x)||u_{hi}(t,x)|dxdydt
\lesssim \eqref{et6}.
\end{align*}
When the derivative falls on $u_{hi}$, we first integrate by parts and then take the absolute values inside the integrals
to obtain, as an error, a scalar multiple of
\begin{align*}
&\int_{I_{0}} \int_{\R^n} \int_{\R^n} |u_{hi}(t,y)|^2|\nabla P_{lo}\bigl(|u|^{\frac{4}{n-2}}u\bigr)(t,x)||u_{hi}(t,x)|dxdydt\\
&\quad + \int_{I_{0}} \int_{\R^n} \int_{\R^n} \frac{|u_{hi}(t,y)|^2|P_{lo}\bigl(|u|^{\frac{4}{n-2}}u\bigr)(t,x)||u_{hi}(t,x)|}{|x-y|}dxdydt.
\end{align*}
The first term on the right-hand side of the above inequality is controlled by \eqref{et6}. The second term we estimate
via \eqref{diff1} by
\begin{align*}
&\int_{I_{0}} \int_{\R^n} \int_{\R^n} \frac{|u_{hi}(t,y)|^2|P_{lo}\bigl(|u_{hi}|^{\frac{4}{n-2}}u_{hi}\bigr)(t,x)||u_{hi}(t,x)|}{|x-y|}dxdydt\\
&\quad +\int_{I_{0}} \int_{\R^n} \int_{\R^n} \frac{|u_{hi}(t,y)|^2|P_{lo}\bigl(|u|^{\frac{4}{n-2}}u-|u_{hi}|^{\frac{4}{n-2}}u_{hi}\bigr)(t,x)||u_{hi}(t,x)|}{|x-y|}dxdydt\\
&\quad \phantom{+\int_{I_{0}} \int_{\R^n} \int_{\R^n}} \lesssim \eqref{et7}+\eqref{et8}+\eqref{et9}.
\end{align*}
\end{proof}

\subsection{Strichartz control on low and high frequencies} The purpose of this section is to obtain estimates on the
low and high-frequency parts of $u$, which we will use to bound the error terms in Proposition \ref{etprop}.
Throughout this section we take $n\geq 6$.

\begin{proposition}[Strichartz control on low and high frequencies]\label{Scontrol}
There exists a constant $C_1$ possibly depending on the energy, but not on any of the $\eta$'s, such that we have the
following estimates. The low frequencies satisfy
\begin{align}
\|u_{lo}\|_{\dot{S}^{1}(I_0\times \R^n)}\leq C_1 \eta_2^{\frac{4}{(n-2)^2}}. \label{slf}
\end{align}
The high frequencies of $u$ can be split into a `good' and a `bad' part, $u_{hi}=g+b$, such that
\begin{align}
&\|g\|_{\dot{S}^{0}(I_0\times \R^n)}\leq C_1 \eta_{2}^{\frac{2}{n-2}}, \label{gS0} \\
&\|g\|_{\dot{S}^{1}(I_0\times \R^n)}\leq C_1, \label{gS1} \\
&\||\nabla|^{-\frac{2}{n-2}}b\|_{L_t^2L_x^{\frac{2n(n-2)}{n^2-3n-2}}(I_0\times \R^n)}\leq C_1\eta_1^{\frac{1}{4}}. \label{bS}
\end{align}
\end{proposition}

\begin{proof}
We define the two functions, $g$ and $b$, to be the unique solutions to the initial value problems
\begin{align}\label{ecuatie g}
\begin{cases}
(i\partial_t+\Delta)g=G + P_{hi}F(u_{lo}) + P_{hi}\bigl(F(u_{lo}+g)-F(g)-F(u_{lo})\bigr)\\
g(t_0)=u_{hi}(t_0)
\end{cases}
\end{align}
and
\begin{align*}
\begin{cases}
(i\partial_t+\Delta)b= B+ P_{hi}\bigl(bF_z(u_{lo}+g)+\bar{b}F_{\bar z}(u_{lo}+g)\bigr)\\
\phantom{(i\partial_t+\Delta)b=}+P_{hi}\bigl(F(u_{lo}+g+b)-F(u_{lo}+g)-bF_z(u_{lo}+g)-\bar{b}F_{\bar z}(u_{lo}+g)\bigr)\\
b(t_0)=0,
\end{cases}
\end{align*}
where $F:\C\to \C$ is the function given by $F(z)=|z|^{\frac{4}{n-2}}z$ and $G$ and $B$ are such that
$P_{hi}\bigl(|g|^{\frac{4}{n-2}}g\bigl)=G+B$, as we will explain momentarily.  Note that $b\equiv u_{hi}-g$.
In Appendix~B, we prove the existence and uniqueness of local solutions to \eqref{ecuatie g}, which implies
the existence and uniqueness of $b$.

In order to prove Proposition \ref{Scontrol} we will use a bootstrap argument. Fix $t_0:=\inf I_0$ and let $\Omega_1$
be the set of all times $T\in I_0$ such that \eqref{slf} through \eqref{bS} hold on $[t_0, T]$ with $g$ and $b$ defined
above.

Define also $\Omega_2$ to be the set of all times $T\in I_0$ such that \eqref{slf} through \eqref{bS}
hold on $[t_0, T]$ with $C_1$ replaced by $2C_1$ and $g$ and $b$ defined above.
More precisely, for $T\in \Omega_2$ we have
\begin{align}\label{uloass}
\|u_{lo}\|_{\dot{S}^{1}([t_0, T]\times\R^n)}\leq 2C_1\eta_2^{\frac{4}{(n-2)^2}}
\end{align}
and
\begin{align}
&\|g\|_{\dot{S}^{0}([t_0, T]\times\R^n)}\leq 2C_1\eta_{2}^{\frac{2}{n-2}},  \label{gS0ass} \\
&\|g\|_{\dot{S}^{1}([t_0, T]\times\R^n)}\leq 2C_1, \label{gS1ass} \\
&\||\nabla|^{-\frac{2}{n-2}}b\|_{L_t^2L_x^{\frac{2n(n-2)}{n^2-3n-2}}([t_0, T]\times\R^n)}\leq 2C_1\eta_1^{\frac{1}{4}} \label{bSass}.
\end{align}

In order to run a bootstrap argument successfully, we need to check four things:\\
$\bullet$ First, we need to see that $t_0\in \Omega_1$; this follows immediately from the definition of $u_{lo}$, $g$,
and $b$ at the time $t=t_0$, provided $C_1$ is sufficiently large. \\
$\bullet$ Secondly, we need $\Omega_1$ to be closed; this follows from the definition of $\Omega_1$ and Fatou's lemma.\\
$\bullet$ Next, we need to prove that if $T\in\Omega_1$, then there exists a small neighborhood of $T$ contained in $\Omega_2$.
This property follows for $u_{lo}$ from the dominated convergence theorem and the fact that $u_{lo}$ is not only
in $\dot{S}^1([t_0,T]\times\R^n)$, but also in $C_t^0\dot{H}^1_x([t_0,T]\times\R^n)$ because of the smoothing effect
of the free propagator. As far as the high frequencies are concerned, it suffices to check it for $g$ since
$b\equiv u_{hi}-g$ and thus the claim for $b$ follows again from the dominated convergence theorem. To prove this property for
the function $g$ basically amounts to proving existence of $g$ on a tiny interval, since the dominated convergence
theorem and the smoothing effect of the free propagator can be used, as before, to conclude our claim. The existence
of $g$ is proved in the usual way: showing convergence of the iterates; the proof is standard and we will defer it
to Appendix~B as to not disrupt the flow of the presentation.\\
$\bullet$ The last thing one needs to check in order to complete the bootstrap argument is that
$\Omega_2\subset\Omega_1$ and this is what we will focus on for the rest of the proof of Proposition~\ref{Scontrol}.
Fix therefore $T\in \Omega_2$. Throughout the rest of the proof all spacetime norms will be on $[t_0, T]\times\R^n$.

Before we move on with our proof, let us make a few observations. First, note that by \eqref{gS1ass} and the conservation
of energy, we get
\begin{align}\label{benergy}
\|\nabla b\|_{\infty,2}
\leq \|\nabla u_{hi}\|_{\infty,2}+\|\nabla g\|_{\infty,2}
\leq 3 C_1,
\end{align}
provided $C_1$ is sufficiently large. Also, from \eqref{smallmass} and \eqref{gS0ass} and by taking $C_1$ sufficiently
large, one easily sees that the mass of $b$ is small:
\begin{align}\label{bmass}
\|b\|_{\infty,2}
\leq \|u_{hi}\|_{\infty,2}+\| g\|_{\infty,2}
\leq 3C_1 \eta_2^{\frac{2}{n-2}}.
\end{align}
Interpolating between \eqref{bSass} and \eqref{benergy}, we obtain the following estimate which we will repeatedly use
in what follows:
\begin{align}\label{b in 2n/n-2}
\|b\|_{\frac{2n}{n-2},\frac{2n^2}{(n+1)(n-2)}}
\leq C\||\nabla|^{-\frac{2}{n-2}}b\|_{2,\frac{2n(n-2)}{n^2-3n-2}}^{\frac{n-2}{n}}\|\nabla b\|_{\infty,2}^{\frac{2}{n}}
\leq 3C C_1 \eta_1^{\frac{n-2}{4n}},
\end{align}
where $C$ is a positive real constant.
Also, by interpolation, \eqref{gS0ass}, and \eqref{gS1ass}, we estimate
\begin{align*}
\|\nabla^{\frac{1}{4}}g\|_{4,\frac{2n}{n-1}}
\lesssim\|\nabla g\|^{\frac{1}{4}}_{4,\frac{2n}{n-1}}\|g\|^{\frac{3}{4}}_{4,\frac{2n}{n-1}}
\lesssim \|g\|_{\dot{S}^1}^{\frac{1}{4}}\|g\|_{\dot{S}^0}^{\frac{3}{4}}
\leq \eta_2^{\frac{1}{n-2}}.
\end{align*}
Sobolev embedding dictates
\begin{align}\label{g in L4}
\||\nabla|^{-\frac{n-3}{4}}g\|_{4,4}
 \lesssim\|\nabla^{\frac{1}{4}}g\|_{4,\frac{2n}{n-1}}
 \leq \eta_2^{\frac{1}{2(n-2)}}
\end{align}
 and hence, by the triangle inequality, \eqref{bootstraphyp2} and \eqref{g in L4} yield
\begin{align}\label{b in L4}
\||\nabla|^{-\frac{n-3}{4}}b\|_{4,4}\lesssim (C_0\eta_1)^{\frac{1}{4}}.
\end{align}

We are now ready to resume the proof. We consider the low frequencies first. Strichartz's inequality yields the bound
\begin{align}\label{ulo1}
\|u_{lo}\|_{\dot{S}^1}\lesssim \|e^{i(t-t_0)\Delta}u_{lo}(t_0)\|_{\dot{S}^1}
    +\|\nabla P_{lo}F(u)\|_{2,\frac{2n}{n+2}}.
\end{align}
By Strichartz and \eqref{lowfreqsmall},
$$
\|e^{i(t-t_0)\Delta}u_{lo}(t_0)\|_{\dot{S}^1}\lesssim \|\nabla u_{lo}\|_{\infty,2}
  \lesssim  \eta_2 \leq \frac{C_1}{100}\eta_2^{\frac{4}{(n-2)^2}}.
$$
To estimate the second term on the right-hand side of \eqref{ulo1}, we write
$$
\nabla P_{lo}F(u)= \nabla P_{lo}F(u_{lo})+\nabla P_{lo}\bigl(F(u)-F(u_{lo})\bigr).
$$
By \eqref{lowfreqsmall} and \eqref{uloass}, we estimate
\begin{align*}
\|\nabla P_{lo}F(u_{lo})\|_{2,\frac{2n}{n+2}}
\lesssim \|\nabla u_{lo}\|_{2,\frac{2n}{n-2}}\|u_{lo}\|_{\infty,\frac{2n}{n-2}}^{\frac{4}{n-2}}
\lesssim \eta_2^{\frac{4}{n-2}}\|u_{lo}\|_{\dot{S}^1}\leq \frac{C_1}{100}\eta_2^{\frac{4}{(n-2)^2}}.
\end{align*}
Using Bernstein to drop the derivative in front of $P_{lo}$ and then replacing the projection $P_{lo}$ by a positive-kernel operator $P_{lo}'$
having the same $L_x^p$-mapping and Bernstein properties as $P_{lo}$ (see Subsection~1.2 for the definition of $P_{lo}'$), we estimate
\begin{align*}
\|\nabla P_{lo}\bigl(F(u)-F(u_{lo})\bigr)\|_{2,\frac{2n}{n+2}}
\lesssim \|P_{lo}'\bigl(|u_{lo}|^{\frac{4}{n-2}}|u_{hi}|\bigr)\|_{2,\frac{2n}{n+2}}
+\|P_{lo}'\bigl(|u_{hi}|^{\frac{n+2}{n-2}}\bigr)\|_{2,\frac{2n}{n+2}}.
\end{align*}
Decomposing $u_{hi}=g+b$ and exploiting the positivity of the operator $P_{lo}'$, we get
\begin{align*}
\|\nabla P_{lo}\bigl(F(u)-F(u_{lo})\bigr)\|_{2,\frac{2n}{n+2}}
&\lesssim \|P_{lo}'\bigl(|u_{lo}|^{\frac{4}{n-2}}|g|\bigr)\|_{2,\frac{2n}{n+2}}
+\|P_{lo}'\bigl(|u_{lo}|^{\frac{4}{n-2}}|b|\bigr)\|_{2,\frac{2n}{n+2}}\\
&\quad +\|P_{lo}'\bigl(|g|^{\frac{n+2}{n-2}}\bigr)\|_{2,\frac{2n}{n+2}}
+\|P_{lo}'\bigl(|b|^{\frac{n+2}{n-2}}\bigr)\|_{2,\frac{2n}{n+2}}.
\end{align*}
Using Bernstein to lower the spatial exponent when necessary, \eqref{lowfreqsmall}, as well as our assumptions
\eqref{uloass} through \eqref{bSass}, \eqref{bmass}, and \eqref{b in 2n/n-2}, we estimate
\begin{align*}
\|P_{lo}'\bigl(|u_{lo}|^{\frac{4}{n-2}}|g|\bigr)\|_{2,\frac{2n}{n+2}}
&\lesssim \|u_{lo}\|_{\infty,\frac{2n}{n-2}}^{\frac{4}{n-2}}\|g\|_{2,\frac{2n}{n-2}}
 \lesssim \eta_2^{\frac{4}{n-2}}\|g\|_{\dot{S}^0}
 \leq \frac{C_1}{100}\eta_2^{\frac{4}{(n-2)^2}},\\
\|P_{lo}'\bigl(|u_{lo}|^{\frac{4}{n-2}}|b|\bigr)\|_{2,\frac{2n}{n+2}}
&\lesssim \|P_{lo}'\bigl(|u_{lo}|^{\frac{4}{n-2}}|b|\bigr)\|_{2,\frac{2n(n+2)}{n^2+5n-2}}\\
&\lesssim \|b\|_{\frac{2(n+2)}{n-2},\frac{2n}{n-1}}\|u_{lo}\|_{\frac{2(n+2)}{n-2},\frac{2(n+2)}{n-2}}^{\frac{4}{n-2}}\\
&\lesssim  \|b\|_{\frac{2n}{n-2},\frac{2n^2}{(n+1)(n-2)}}^{\frac{n}{n+2}}\|b\|_{\infty,2}^{\frac{2}{n+2}}\|u_{lo}\|_{\dot{S}^1}^{\frac{4}{n-2}}
\leq \frac{C_1}{100}\eta_2^{\frac{4}{(n-2)^2}},\\
\|P_{lo}'\bigl(|g|^{\frac{n+2}{n-2}}\bigr)\|_{2,\frac{2n}{n+2}}
&\lesssim \|P_{lo}'\bigl(|g|^{\frac{n+2}{n-2}}\bigr)\|_{2,\frac{2n(n-2)}{n^2+4}}
 \lesssim \| g\|_{\frac{2(n+2)}{n-2},\frac{2n(n+2)}{n^2+4}}^{\frac{n+2}{n-2}}
 \lesssim \|g\|_{\dot{S}^0}^{\frac{n+2}{n-2}}\\
&\leq \frac{C_1}{100}\eta_2^{\frac{4}{(n-2)^2}},\\
\|P_{lo}'\bigl(|b|^{\frac{n+2}{n-2}}\bigr)\|_{2,\frac{2n}{n+2}}
&\lesssim \| P_{lo}'\bigl(|b|^{\frac{n+2}{n-2}}\bigr)\|_{2,\frac{2n(n-2)}{(n+2)(n-1)}}
 \lesssim \|b\|_{\frac{2(n+2)}{n-2},\frac{2n}{n-1}}^{\frac{n+2}{n-2}}\\
&\lesssim \|b\|_{\frac{2n}{n-2},\frac{2n^2}{(n+1)(n-2)}}^{\frac{n}{n-2}}\|b\|_{\infty,2}^{\frac{2}{n-2}}
 \leq \frac{C_1}{100}\eta_2^{\frac{4}{(n-2)^2}}.
\end{align*}
Therefore, putting everything together we obtain control over the low frequencies,
$$
\|u_{lo}\|_{\dot{S}^1}\leq C_1\eta_2^{\frac{4}{(n-2)^2}}.
$$

We turn now to the high frequencies. We will first clarify what $G$ and $B$ are. The reason we need to split
$P_{hi}F(g)$ into $G+B$ is that $P_{hi}F(g)$ is neither `good enough' to be part of $g$ (as one cannot close the
bootstrap for the $\dot{S}^1$ bound on $g$) nor `sufficiently bad' to be part of $b$ (as it's not sufficiently fast
decaying to belong to the appropriate $L_x^p$ spaces, unless $6\leq n\leq 14$). We thus use an interpolation trick to
split $P_{hi}F(g)$ into a part which is small and has high spatial integrability, $G$, and a part which has low spatial
integrability, $B$. Indeed, we have
\begin{lemma}\label{GB}
There exist two functions $G$ and $B$ such that $P_{hi}F(g)=G+B$ and moreover,
\begin{gather}
\|G\|_{\frac{2(n+2)}{n+4},\frac{2(n+2)}{n+4}}\leq \eta_2^c \eta_2^{\frac{2}{n-2}}\label{GN0}\\
\|\nabla G\|_{\frac{2(n+2)}{n+4},\frac{2(n+2)}{n+4}}\ll \eta_1^{100}\label{GN1}\\
\||\nabla|^{-\frac{2}{n-2}}B\|_{\frac{2(n-2)(n+2)}{n^2+3n-14},\frac{2(n-2)(n+2)}{n^2+3n-14}}\ll\eta_1^{100},\label{B}
\end{gather}
where $c>0$ is a small constant depending on the dimension $n$.
\end{lemma}

\begin{proof}
Let us first note that for $6\leq n\leq 14$, we can choose $G:=0$ and $B:=P_{hi}F(g)$, since by Bernstein and interpolation,
\begin{align*}
\||\nabla|^{-\frac{2}{n-2}}P_{hi}F(g)&\|_{\frac{2(n-2)(n+2)}{n^2+3n-14},\frac{2(n-2)(n+2)}{n^2+3n-14}}\\
&\lesssim \|P_{hi}F(g)\|_{\frac{2(n-2)(n+2)}{n^2+3n-14},\frac{2(n-2)(n+2)}{n^2+3n-14}}\\
&\lesssim \|g\|_{\frac{2(n+2)^2}{n^2+3n-14},\frac{2(n+2)^2}{n^2+3n-14}}^{\frac{n+2}{n-2}}\\
&\lesssim \|g\|_{\frac{2(n+2)^2}{n^2+3n-14},\frac{2n(n+2)^2}{n^3+2n^2-2n+28}}^{\frac{3n^2-4n-20}{2(n+2)(n-2)}} \|g\|_{\frac{2(n+2)^2}{n^2+3n-14},\frac{2n(n+2)^2}{n^3-10n+20}}^{\frac{(n+2)(14-n)}{2(n+2)(n-2)}}\\
&\lesssim \|g\|_{\dot{S}^0}^{\frac{3n^2-4n-20}{2(n+2)(n-2)}}  \|g\|_{\dot{S}^1}^{\frac{(n+2)(14-n)}{2(n+2)(n-2)}}\\
&\lesssim \eta_2^c \ll\eta_1^{100},
\end{align*}
where the last line follows from \eqref{gS0ass} and \eqref{gS1ass}. Here, $c$ is a small positive constant\footnote{Throughout
the proof, the constant $c$ may vary from line to line; however, it will always remain positive and will depend only
on the dimension.} depending only on the dimension $n$.

We consider next the case $n> 14$. To decompose $P_{hi}F(g)$ into a part with high spatial integrability and
a part with low spatial integrability, we first need $P_{hi}F(g)$ to belong to an intermediate $L_x^p$ space.
We choose the space $L_{t,x}^{\frac{2(n-2)}{n}}$ and use \eqref{gS0ass} and \eqref{gS1ass} to estimate
\begin{align}\label{F(g)}
\|F(g)\|_{\frac{2(n-2)}{n},\frac{2(n-2)}{n}}
\lesssim \|g\|_{\frac{2(n+2)}{n},\frac{2(n+2)}{n}}^{\frac{n+2}{n-2}}
\lesssim \|g\|_{\dot{S}^0}^{\frac{n+2}{n-2}}\leq \eta_2^c \|g\|_{\dot{S}^0}
\end{align}
and, by the boundedness of the Riesz potentials on $L_x^p$, $1<p<\infty$,
\begin{align}\label{nabla F(g)}
\||\nabla| F(g)\|_{\frac{2(n-2)}{n},\frac{2(n-2)}{n}}
&\lesssim \|\nabla F(g)\|_{\frac{2(n-2)}{n},\frac{2(n-2)}{n}}
\lesssim \|\nabla g\|_{\frac{2(n+2)}{n},\frac{2(n+2)}{n}}\|g\|_{\frac{2(n+2)}{n},\frac{2(n+2)}{n}}^{\frac{4}{n-2}}\notag\\
&\lesssim \|g\|_{\dot{S}^1}\|g\|_{\dot{S}^0}^{\frac{4}{n-2}}\leq \eta_2^c.
\end{align}

We now decompose
\begin{align}\label{decomp1}
P_{hi}F(g)=P_{1<\cdot<\eta_2^{-100}}F(g) + P_{\geq\eta_2^{-100}}F(g).
\end{align}

Consider first the term in \eqref{decomp1} involving very high frequencies. Writing
$$
P_{\geq\eta_2^{-100}}F(g)=|\nabla|^{-1}P_{\geq\eta_2^{-100}}\bigl(|\nabla| F(g)\bigr),
$$
we define
$$
G_{vhi}:=|\nabla|^{-1}P_{\geq\eta_2^{-100}}\bigl(\chi_{\{||\nabla| F(g)|\leq 1\}}|\nabla| F(g) \bigr)
$$
and
$$
B_{vhi}:=|\nabla|^{-1}P_{\geq\eta_2^{-100}}\bigl((1-\chi_{\{||\nabla| F(g)|\leq 1\}})|\nabla| F(g) \bigr),
$$
where $\chi_{\{||\nabla| F(g)|\leq 1\}}$ is a smooth cutoff.

By Bernstein, H\"older, and \eqref{nabla F(g)}, we estimate
\begin{align}\label{Gvhi N0}
\|G_{vhi}\|_{\frac{2(n+2)}{n+4},\frac{2(n+2)}{n+4}}
&\lesssim \eta_2^{100}\|P_{\geq\eta_2^{-100}}\bigl(\chi_{\{||\nabla| F(g)|\leq 1\}}|\nabla| F(g) \bigr)\|_{\frac{2(n+2)}{n+4},\frac{2(n+2)}{n+4}}\notag\\
&\lesssim \eta_2^{100}\||\nabla| F(g)\|_{\frac{2(n-2)}{n},\frac{2(n-2)}{n}}^{\frac{(n-2)(n+4)}{n(n+2)}}
\leq \eta_2^{100}.
\end{align}
By the boundedness of the Riesz transforms on $L_x^p$ with $1<p<\infty$ and \eqref{nabla F(g)}, we estimate
\begin{align}\label{Gvhi N1}
\|\nabla G_{vhi}\|_{\frac{2(n+2)}{n+4},\frac{2(n+2)}{n+4}}
&\lesssim \|P_{\geq\eta_2^{-100}}\bigl(\chi_{\{||\nabla| F(g)|\leq 1\}}|\nabla| F(g) \bigr)\|_{\frac{2(n+2)}{n+4},\frac{2(n+2)}{n+4}}\notag\\
&\lesssim \||\nabla| F(g)\|_{\frac{2(n-2)}{n},\frac{2(n-2)}{n}}^{\frac{(n-2)(n+4)}{n(n+2)}}
\leq \eta_2^c \ll \eta_1^{100}.
\end{align}
By Bernstein, H\"older, and \eqref{nabla F(g)}, we get
\begin{align}\label{Bvhi}
\||\nabla|^{-\frac{2}{n-2}}B_{vhi}\|_{\frac{2(n-2)(n+2)}{n^2+3n-14},\frac{2(n-2)(n+2)}{n^2+3n-14}}
&\lesssim\eta_2^{\frac{200}{n-2}} \|B_{vhi}\|_{\frac{2(n-2)(n+2)}{n^2+3n-14},\frac{2(n-2)(n+2)}{n^2+3n-14}}\notag\\
&\lesssim \eta_2^{\frac{100n}{n-2}}\||\nabla| F(g)\|_{\frac{2(n-2)}{n},\frac{2(n-2)}{n}}^{\frac{n^2+3n-14}{n(n+2)}}
\leq \eta_2^{100}.
\end{align}

We consider next the medium frequency term in \eqref{decomp1} and write
$$
P_{1<\cdot<\eta_2^{-100}}F(g)\sim \sum_{1<N<\eta_2^{-100}}\tilde{P_N} P_N F(g),
$$
where $\tilde{P_N}$ is an operator having the same properties as $P_N$ and double support on the Fourier side.
For dyadic $N$'s between $1$ and $\eta_2^{-100}$, we define
\begin{gather*}
G_N:=\chi_{\{|P_NF(g)|\leq 1/N\}} P_NF(g),\\
B_N:= (1-\chi_{\{|P_NF(g)|\leq 1/N\}}) P_N F(g),
\end{gather*}
where $\chi_{\{|P_NF(g)|\leq 1/N\}}$ are again smooth cutoffs. We define
\begin{align*}
G:=G_{med}+G_{vhi} \ \ \text{and} \ \ B:=B_{med}+B_{vhi},
\end{align*}
where
\begin{gather*}
G_{med}:=\sum_{1<N<\eta_2^{-100}}\tilde{P_N} G_N\ \ \text{and} \ \
B_{med}:=\sum_{1<N<\eta_2^{-100}}\tilde{P_N} B_N.
\end{gather*}
Using H\"older and \eqref{F(g)}, we estimate
\begin{align}\label{GN N0}
\|\tilde{P_N} G_N\|_{\frac{2(n+2)}{n+4},\frac{2(n+2)}{n+4}}
\lesssim \| F(g)\|_{\frac{2(n-2)}{n},\frac{2(n-2)}{n}}^{\frac{(n-2)(n+4)}{n(n+2)}}N^{-\frac{8}{n(n+2)}}
\leq \|g\|_{\dot{S}^0}^{\frac{(n+4)}{n}}N^{-\frac{8}{n(n+2)}},
\end{align}
which implies together with \eqref{gS0ass} that
\begin{align}\label{Gmed N0}
\|G_{med}\|_{\frac{2(n+2)}{n+4},\frac{2(n+2)}{n+4}}
&\lesssim \sum_{1<N<\eta_2^{-100}}\|\tilde{P_N} G_N\|_{\frac{2(n+2)}{n+4},\frac{2(n+2)}{n+4}}\notag\\
&\lesssim \|g\|_{\dot{S}^0}^{\frac{(n+4)}{n}}\sum_{1<N<\eta_2^{-100}}N^{-\frac{8}{n(n+2)}}\notag\\
&\lesssim \|g\|_{\dot{S}^0}^{\frac{(n+4)}{n}}\leq \eta_2^c \|g\|_{\dot{S}^0}.
\end{align}

Now, by Bernstein, H\"older, and \eqref{nabla F(g)}, we get
\begin{align}\label{GN N1}
\|\nabla \tilde{P_N} G_N\|_{\frac{2(n+2)}{n+4},\frac{2(n+2)}{n+4}}
&\lesssim N \| P_N F(g)\|_{\frac{2(n-2)}{n},\frac{2(n-2)}{n}}^{\frac{(n-2)(n+4)}{n(n+2)}}N^{-\frac{8}{n(n+2)}}\notag\\
&\lesssim N N^{-\frac{(n-2)(n+4)}{n(n+2)}}\| P_N \nabla F(g)\|_{\frac{2(n-2)}{n},\frac{2(n-2)}{n}}^{\frac{(n-2)(n+4)}{n(n+2)}}N^{-\frac{8}{n(n+2)}}\notag\\
&\lesssim \eta_2^c.
\end{align}
As there are about $\log(\eta_2^{-1})$ dyadic numbers $N$ between $1$ and $\eta_2^{-100}$, by \eqref{nabla F(g)}
and \eqref{GN N1}, we get
\begin{align}\label{Gmed N1}
\|\nabla G_{med}\|_{\frac{2(n+2)}{n+4},\frac{2(n+2)}{n+4}}
&\lesssim \sum_{1<N<\eta_2^{-100}}\|\nabla \tilde{P_N} G_N\|_{\frac{2(n+2)}{n+4},\frac{2(n+2)}{n+4}}\notag\\
&\lesssim \log(\eta_2^{-1}) \eta_2^c
\ll\eta_1^{100}.
\end{align}

By H\"older, Bernstein, and \eqref{nabla F(g)}, we estimate
\begin{align}\label{BN}
\||\nabla|^{-\frac{2}{n-2}}\tilde{P_N}B_N&\|_{\frac{2(n-2)(n+2)}{n^2+3n-14},\frac{2(n-2)(n+2)}{n^2+3n-14}}\notag\\
&\lesssim N^{-\frac{2}{n-2}}\|P_N F(g)\|_{\frac{2(n-2)}{n},\frac{2(n-2)}{n}}^{\frac{n^2+3n-14}{n(n+2)}}N^{\frac{n-14}{n(n+2)}}\notag\\
&\lesssim N^{-\frac{2}{n-2}}N^{-\frac{n^2+3n-14}{n(n+2)}}\| P_N \nabla F(g)\|_{\frac{2(n-2)}{n},\frac{2(n-2)}{n}}^{\frac{n^2+3n-14}{n(n+2)}}N^{\frac{n-14}{n(n+2)}}\notag\\
&\lesssim N^{-\frac{n}{n-2}}\eta_2^c.
\end{align}
Hence,
\begin{align}\label{Bmed}
\||\nabla|^{-\frac{2}{n-2}}B_{med}&\|_{\frac{2(n-2)(n+2)}{n^2+3n-14},\frac{2(n-2)(n+2)}{n^2+3n-14}}\notag\\
&\lesssim \sum_{1<N<\eta_2^{-100}}\||\nabla|^{-\frac{2}{n-2}}\tilde{P_N}B_N\|_{\frac{2(n-2)(n+2)}{n^2+3n-14},\frac{2(n-2)(n+2)}{n^2+3n-14}}\notag\\
&\lesssim \eta_2^c \sum_{1<N<\eta_2^{-100}} N^{-\frac{n}{n-2}}
 \lesssim \eta_2^c.
\end{align}

Thus, by \eqref{Gvhi N0} and \eqref{Gmed N0} we get \eqref{GN0}, by \eqref{Gvhi N1} and \eqref{Gmed N1} we get
\eqref{GN1}, and by \eqref{Bvhi} and \eqref{Bmed} we get \eqref{B}.
\end{proof}

We are now ready to resume the bootstrap for $g$ and $b$. Consider first the `good' part, $g$. By Strichartz, Bernstein,
\eqref{diff2}, \eqref{lowfreqsmall}, \eqref{smallmass}, \eqref{uloass}, \eqref{gS0ass}, and \eqref{GN0}, we estimate
\begin{align*}
\|g\|_{\dot{S}^0}
&\lesssim \|u_{hi}\|_{\infty,2}+ \|G\|_{\frac{2(n+2)}{n+4},\frac{2(n+2)}{n+4}} + \|P_{hi}F(u_{lo})\|_{2,\frac{2n}{n+2}} \\
&\quad+\|P_{hi}\bigl(F(u_{lo}+g)-F(u_{lo})-F(g)\bigr)\|_{2, \frac{2n}{n+2}}\\
&\lesssim \eta_2+ \eta_2^c \eta_2^{\frac{2}{n-2}} + \|\nabla P_{hi}F(u_{lo})\|_{2,\frac{2n}{n+2}} + \|g|u_{lo}|^{\frac{4}{n-2}}\chi_{\{|g|\leq |u_{lo}|\}}\|_{2,\frac{2n}{n+2}}\\
&\quad +\|u_{lo}|g|^{\frac{4}{n-2}}\chi_{\{|u_{lo}|< |g|\}}\|_{2,\frac{2n}{n+2}}\\
&\lesssim \eta_2^c\eta_2^{\frac{2}{n-2}}+ \|u_{lo}\|_{\dot{S}^1}\|u_{lo}\|_{\infty,\frac{2n}{n-2}}^{\frac{4}{n-2}}+ \|g|u_{lo}|^{\frac{4}{n-2}}\|_{2,\frac{2n}{n+2}}\\
&\lesssim \eta_2^c\eta_2^{\frac{2}{n-2}}+\|g\|_{\dot{S}^0}\|u_{lo}\|_{\infty,\frac{2n}{n-2}}^{\frac{4}{n-2}}\\
&\leq C_1 \eta_2^{\frac{2}{n-2}}.
\end{align*}
Similarly, by Strichartz, \eqref{diff3}, \eqref{uloass}, \eqref{gS1ass}, and \eqref{GN1}, we estimate
\begin{align*}
\|g\|_{\dot{S}^1}
&\lesssim \|\nabla u_{hi}\|_{\infty, 2} + \|\nabla G\|_{\frac{2(n+2)}{n+4},\frac{2(n+2)}{n+4}}+ \|\nabla P_{hi}F(u_{lo})\|_{2,\frac{2n}{n+2}}\\
&\quad +\|\nabla P_{hi}\bigl(F(u_{lo}+g)-F(u_{lo})-F(g)\bigr)\|_{2, \frac{2n}{n+2}}\\
&\lesssim 1 + \eta_1^{100}+\|u_{lo}\|_{\dot{S}^1}\|u_{lo}\|_{\infty,\frac{2n}{n-2}}^{\frac{4}{n-2}}+\|\nabla g|u_{lo}|^{\frac{4}{n-2}}\|_{2,\frac{2n}{n+2}}+\|\nabla u_{lo}|g|^{\frac{4}{n-2}}\|_{2,\frac{2n}{n+2}}\\
&\lesssim 1 + \|u_{lo}\|_{\infty,\frac{2n}{n-2}}^{\frac{4}{n-2}}\|g\|_{\dot{S}^1}+\|u_{lo}\|_{\dot{S}^1}\|g\|_{\dot{S}^1}^{\frac{4}{n-2}}\\
&\leq C_1,
\end{align*}
provided $C_1$ is sufficiently large.

We turn now to $b$. Using the triangle inequality and the inhomogeneous Strichartz estimates \eqref{inhomS} and
\eqref{inhomS2}, we estimate
\begin{align*}
\||\nabla|^{-\frac{2}{n-2}}b\|_{2,\frac{2n(n-2)}{n^2-3n-2}}
&\lesssim \||\nabla|^{-\frac{2}{n-2}}B\|_{\frac{2(n-2)(n+2)}{n^2+3n-14},\frac{2(n-2)(n+2)}{n^2+3n-14}}\\
&\quad+ \||\nabla|^{-\frac{2}{n-2}} P_{hi}\bigl(bF_z(u_{lo}+g)+\bar{b}F_{\bar z}(u_{lo}+g)\bigr)\|_{2,\frac{2n(n-2)}{n^2+n-10}}\\
&\quad +\||\nabla|^{-\frac{2}{n-2}}P_{hi}\Bigl(F(u_{lo}+g+b)-F(u_{lo}+g)\\
&\qquad \quad -bF_z(u_{lo}+g)-\bar{b}F_{\bar z}(u_{lo}+g)\Bigr)\|_{2,\frac{2n(n-2)}{n^2+n-10}}.
\end{align*}
By \eqref{B},
$$
\||\nabla|^{-\frac{2}{n-2}}B\|_{\frac{2(n-2)(n+2)}{n^2+3n-14},\frac{2(n-2)(n+2)}{n^2+3n-14}}\ll\eta_1^{100}\leq \frac{C_1}{100}\eta_1^{\frac{1}{4}}.
$$
By the Fundamental Theorem of Calculus, we have
\begin{align*}
F(z+w)-F(z)-wF_z(z)-\bar{w}F_{\bar{z}}(z)
&=w\int_0^1[F_z(z+tw)-F_z(z)]dt\\
&\quad +\bar{w}\int_0^1[F_{\bar{z}}(z+tw)-F_{\bar{z}}(z)]dt.
\end{align*}
As $z\to F_z(z)$ and $z\to F_{\bar{z}}(z)$ are H\"older continuous of order $\frac{4}{n-2}$, we see that
$$
|F(z+w)-F(z)-wF_z(z)-\bar{w}F_{\bar{z}}(z)|\lesssim |w|\int_0^1|tw|^{\frac{4}{n-2}}dt\lesssim |w|^{\frac{n+2}{n-2}}.
$$
Therefore, by Sobolev embedding and the above considerations (with $z=u_{lo}+g$ and $w=b$), we get
\begin{align*}
\||\nabla|^{-\frac{2}{n-2}}P_{hi}\bigl(F(u_{lo}+g+b)-F(u_{lo}+g)-bF_z(u_{lo}+g)-\bar{b}F_{\bar z}(u_{lo}+g)\bigr)\|_{2,\frac{2n(n-2)}{n^2+n-10}}\\
\lesssim \||b|^{\frac{n+2}{n-2}}\|_{2,\frac{2n}{n+3}}\lesssim \|b\|^{\frac{n+2}{n-2}}_{\frac{2(n+2)}{n-2},\frac{2n(n+2)}{(n-2)(n+3)}}.
\end{align*}
Interpolating between \eqref{bSass}, \eqref{benergy}, and \eqref{b in L4}, we obtain
\begin{align}
\|b\|^{\frac{n+2}{n-2}}_{\frac{2(n+2)}{n-2},\frac{2n(n+2)}{(n-2)(n+3)}}
&\lesssim \||\nabla|^{-\frac{n-3}{4}}b\|_{4,4}^{\frac{8}{n^2-3n-2}}
  \||\nabla|^{-\frac{2}{n-2}}b\|_{2,\frac{2n(n-2)}{n^2-3n-2}}^{\frac{n^2-3n-6}{n^2-3n-2}}
  \|\nabla b\|_{\infty,2}^{\frac{4n(n-4)}{(n-2)(n^2-3n-2)}}\notag\\
&\lesssim \eta_1^{\frac{1}{4}+}.\label{b in 2(n+2)/(n-2)}
\end{align}
Hence,
\begin{align*}
\||\nabla|^{-\frac{2}{n-2}}P_{hi}\bigl(F(u_{lo}+g+b)-F(u_{lo}+g)-bF_z(u_{lo}+g)-&\bar{b}F_{\bar z}(u_{lo}+g)\bigr)\|_{2,\frac{2n(n-2)}{n^2+n-10}}\\
&\leq \frac{C_1}{100}\eta_1^{\frac{1}{4}}.
\end{align*}
We turn now towards the remaining two terms,
$$
\||\nabla|^{-\frac{2}{n-2}} P_{hi}\bigl(b|u_{lo}+g|^{\frac{4}{n-2}}\bigr)\|_{2,\frac{2n(n-2)}{n^2+n-10}}
$$
and
$$
\||\nabla|^{-\frac{2}{n-2}}P_{hi}\bigl(\bar{b}|u_{lo}+g|^{\frac{4}{n-2}}\frac{(u_{lo}+g)^2}{|u_{lo}+g|^2}\bigr)\|_{2,\frac{2n(n-2)}{n^2+n-10}}.
$$
As the method of treating them is the same, in particular it appeals to the fact that the maps $z\mapsto |z|^{\frac{4}{n-2}}$
and $z\mapsto |z|^{\frac{4}{n-2}}\frac{z^2}{|z|^2}$ are H\"older continuous of order $\frac{4}{n-2}$, let us pick,
for the sake of the exposition, the first one.
By the triangle inequality, we estimate
\begin{align*}
\||\nabla|^{-\frac{2}{n-2}} P_{hi}\bigl(b|u_{lo}+g|^{\frac{4}{n-2}}&\bigr)\|_{2,\frac{2n(n-2)}{n^2+n-10}}
\lesssim \||\nabla|^{-\frac{2}{n-2}} P_{hi}\bigl(b|u_{lo}|^{\frac{4}{n-2}}\bigr)\|_{2,\frac{2n(n-2)}{n^2+n-10}}\\
&+\||\nabla|^{-\frac{2}{n-2}} P_{hi}\bigl(b|u_{lo}+g|^{\frac{4}{n-2}}-b|u_{lo}|^{\frac{4}{n-2}}\bigr)\|_{2,\frac{2n(n-2)}{n^2+n-10}}.
\end{align*}
Using the H\"older continuity and Sobolev embedding, we bound
$$
\||\nabla|^{-\frac{2}{n-2}} P_{hi}\bigl(b|u_{lo}+g|^{\frac{4}{n-2}}-b|u_{lo}|^{\frac{4}{n-2}}\bigr)\|_{2,\frac{2n(n-2)}{n^2+n-10}}
\lesssim \|b|g|^{\frac{4}{n-2}}\|_{2,\frac{2n}{n+3}}.
$$
Now, by interpolation and our assumptions,
\begin{align*}
\|b|g|^{\frac{4}{n-2}}\|_{2,\frac{2n}{n+3}}
&\lesssim \|b\|_{\frac{2n}{n-2},\frac{2n^2}{(n+1)(n-2)}}\|g\|_{\frac{4n}{n-2},\frac{4n^2}{(2n+1)(n-2)}}^{\frac{4}{n-2}}\\
&\lesssim \|b\|_{\frac{2n}{n-2},\frac{2n^2}{(n+1)(n-2)}}\|g\|_{\frac{4n}{n-2},\frac{2n^2}{n^2-n+2}}^{\frac{3}{n}}\|g\|_{\frac{4n}{n-2},\frac{2n^2}{n^2-3n+2}}^{\frac{n+6}{n(n-2)}}\\
&\lesssim \|b\|_{\frac{2n}{n-2},\frac{2n^2}{(n+1)(n-2)}}\|g\|_{\dot{S}^0}^{\frac{3}{n}}\|g\|_{\dot{S}^1}^{\frac{n+6}{n(n-2)}}\\
&\leq \frac{C_1}{100}\eta_1^{\frac{1}{4}}.
\end{align*}
In order to estimate
$$
\||\nabla|^{-\frac{2}{n-2}} P_{hi}\bigl(b|u_{lo}|^{\frac{4}{n-2}}\bigr)\|_{2,\frac{2n(n-2)}{n^2+n-10}},
$$
we drop the projection onto the high frequencies, $P_{hi}$, but we split $|u_{lo}|^{\frac{4}{n-2}}$ into high and low
frequencies. By the triangle inequality, we get
\begin{align*}
\||\nabla|^{-\frac{2}{n-2}}(b|u_{lo}|^{\frac{4}{n-2}})\|_{2,\frac{2n(n-2)}{n^2+n-10}}
&\lesssim \||\nabla|^{-\frac{2}{n-2}}(bP_{\leq 1/4}|u_{lo}|^{\frac{4}{n-2}})\|_{2,\frac{2n(n-2)}{n^2+n-10}}\\
&\quad+\||\nabla|^{-\frac{2}{n-2}}(bP_{> 1/4}|u_{lo}|^{\frac{4}{n-2}})\|_{2,\frac{2n(n-2)}{n^2+n-10}}.
\end{align*}
As $b$ is high frequency, we see that
\begin{align*}
\||\nabla|^{-\frac{2}{n-2}}(bP_{\leq 1/4}|u_{lo}|^{\frac{4}{n-2}})\|_{2,\frac{2n(n-2)}{n^2+n-10}}
&\lesssim \|(|\nabla|^{-\frac{2}{n-2}}b)P_{\leq 1/4}|u_{lo}|^{\frac{4}{n-2}}\|_{2,\frac{2n(n-2)}{n^2+n-10}}\\
&\lesssim \||\nabla|^{-\frac{2}{n-2}}b\|_{2,\frac{2n(n-2)}{n^2-3n-2}}\|P_{\leq 1/4}|u_{lo}|^{\frac{4}{n-2}}\|_{\infty,\frac{n}{2}}\\
&\lesssim 2C_1\eta_1^{\frac{1}{4}}\|u_{lo}\|_{\infty,\frac{2n}{n-2}}^{\frac{4}{n-2}}\\
&\leq \frac{C_1}{100}\eta_1^{\frac{1}{4}}.
\end{align*}
As far as the term
$$
\||\nabla|^{-\frac{2}{n-2}}(bP_{> 1/4}|u_{lo}|^{\frac{4}{n-2}})\|_{2,\frac{2n(n-2)}{n^2+n-10}}
$$
is concerned, let us note that
\begin{align*}
\||\nabla|^{-\frac{2}{n-2}}&(bP_{> 1/4}|u_{lo}|^{\frac{4}{n-2}})\|_{2,\frac{2n(n-2)}{n^2+n-10}}\\
&\lesssim \||\nabla|^{-\frac{2}{n-2}}b\|_{2,\frac{2n(n-2)}{n^2-3n-2}}\||\nabla|^{\frac{2}{n-2}}P_{> 1/4}|u_{lo}|^{\frac{4}{n-2}}\|_{\infty, \frac{n(n-2)}{2(n-1)}}\\
&\lesssim 2C_1\eta_1^{\frac{1}{4}}\||\nabla|^{\frac{2}{n-2}}P_{> 1/4}|u_{lo}|^{\frac{4}{n-2}}\|_{\infty,\frac{n(n-2)}{2(n-1)}},
\end{align*}
as can easily be seen by taking $j=|\nabla|^{-\frac{2}{n-2}}b$ and
$k=|\nabla|^{\frac{2}{n-2}}P_{> 1/4}|u_{lo}|^{\frac{4}{n-2}}$ in the following
\begin{lemma}\label{multilin op lemma}
\begin{align}\label{multilin op ineq}
\||\nabla|^{-\frac{2}{n-2}}\{(|\nabla|^{\frac{2}{n-2}}j)(|\nabla|^{-\frac{2}{n-2}}k)\}\|_{L_x^{\frac{2n(n-2)}{n^2+n-10}}}
\lesssim \|j\|_{L_x^{\frac{2n(n-2)}{n^2-3n-2}}}\|k\|_{ L_x^{\frac{n(n-2)}{2(n-1)}}}.
\end{align}
\end{lemma}

\begin{proof}
In order to prove Lemma \ref{multilin op lemma}, we decompose the left-hand side into $\pi_{h,h}$, $\pi_{l,h}$, and
$\pi_{h,l}$ which represent the projections onto high-high, low-high, and high-low frequency interactions.

The high-high and low-high frequency interactions are going to be treated in the same manner. Let's consider for example
the first one. A simple application of Sobolev embedding yields
\begin{align*}
\||\nabla|^{-\frac{2}{n-2}}\pi_{h,h}\{(|\nabla|^{\frac{2}{n-2}} j)(|\nabla|^{-\frac{2}{n-2}}k)\}&\|_{L_x^{\frac{2n(n-2)}{n^2+n-10}}}\\
   &\lesssim \|\pi_{h,h}\{(|\nabla|^{\frac{2}{n-2}}j)(|\nabla|^{-\frac{2}{n-2}} k)\}\|_{L_x^{\frac{2n}{n+3}}}.
\end{align*}
Now we only have to notice that the multiplier associated to the operator
$T(j,k)=\pi_{h,h}\{(|\nabla|^{\frac{2}{n-2}} j)(|\nabla|^{-\frac{2}{n-2}}k)\}$, i.e.
$$
\sum_{N\sim M} |\xi_1|^{\frac{2}{n-2}}\widehat{P_N j}(\xi_1) |\xi_2|^{-\frac{2}{n-2}}\widehat{P_M k}(\xi_2),
$$
is a symbol of order one with $\xi = (\xi_1,\xi_2)$,
since then a theorem of R. R. Coifman and Y. Meyer (\cite{coifmey:1}, \cite{coifmey:2}) yields the claim.

To deal with the $\pi_{h,l}$ term, we first notice that the multiplier associated to the operator
$T(j,\tilde k)=|\nabla|^{-\frac{2}{n-2}} \pi_{h,l}\{(|\nabla|^{\frac{2}{n-2}} j)\tilde k\}$, i.e.
$$
\sum_{N\gtrsim M}|\xi_1+\xi_2|^{-\frac{2}{n-2}}|\xi_1|^{\frac{2}{n-2}}\widehat{P_N j}(\xi_1) \widehat{P_M \tilde k}(\xi_2),
$$
is an order one symbol. The result cited above yields
$$
\||\nabla|^{-\frac{2}{n-2}}\pi_{h,l}\{(|\nabla|^{\frac{2}{n-2}} j)(|\nabla|^{-\frac{2}{n-2}}k)\}\|_{L_x^{\frac{2n(n-2)}{n^2+n-10}}}
     \lesssim \|j\|_{L_x^{\frac{2n(n-2)}{n^2-3n-2}}} \||\nabla|^{-\frac{2}{n-2}}k\|_{L_x^{\frac{n}{2}}}.
$$
Finally, Sobolev embedding dictates the estimate
$\||\nabla|^{-\frac{2}{n-2}}k\|_{L_x^{\frac{n}{2}}}\lesssim \|k\|_{ L_x^{\frac{n(n-2)}{2(n-1)}}}$.
\end{proof}

Thus, we are left with the task of estimating
$$
\||\nabla|^{\frac{2}{n-2}}P_{> 1/4}|u_{lo}|^{\frac{4}{n-2}}\|_{\infty,\frac{n(n-2)}{2(n-1)}}.
$$
Note that $P_{> 1/4}|u_{lo}|^{\frac{4}{n-2}}\in \dot{\Lambda}^{\frac{n(n-2)}{2(n-1)}}_{\frac{4}{n-2}}$, that is,
$P_{> 1/4}|u_{lo}|^{\frac{4}{n-2}}$ is homogeneous H\"older continuous of order $\frac{4}{n-2}$ in
$L_x^{\frac{n(n-2)}{2(n-1)}}$. Indeed, as $\nabla u_{lo}\in L_t^\infty L_x^{\frac{2n}{n-1}}$ (by Bernstein), we have
$$
\|u_{lo}^{(h)}(t)-u_{lo}(t)\|_{L_x^{\frac{2n}{n-1}}} \lesssim
|h|\|\nabla u_{lo}\|_{\infty,\frac{2n}{n-1}} \lesssim \eta_2 |h|,
$$
where $u^{(h)}$ denotes the translation $u^{(h)}(x):= u(x-h)$. As the map $z\mapsto |z|^{\frac{4}{n-2}}$ is
H\"older continuous of order $\frac{4}{n-2}$, we see that
$$
\bigl\|\bigl(|u_{lo}|^{\frac{4}{n-2}}\bigr)^{(h)}(t)-|u_{lo}|^{\frac{4}{n-2}}(t)\bigr\|_{L_x^{\frac{n(n-2)}{2(n-1)}}}
\lesssim \eta_2^{\frac{4}{n-2}}|h|^{\frac{4}{n-2}},
$$
which implies $P_{> 1/4}|u_{lo}|^{\frac{4}{n-2}}\in \dot{\Lambda}^{\frac{n(n-2)}{2(n-1)}}_{\frac{4}{n-2}}$.
Furthermore, as $P_{> 1/4}|u_{lo}|^{\frac{4}{n-2}}$ is restricted to high frequencies, the Besov characterization of the
homogeneous H\"older continuous functions (see Chapter~VI in \cite{stein:large}) yields
$$
|\nabla|^{\frac{2}{n-2}}P_{> 1/4}|u_{lo}|^{\frac{4}{n-2}}\in L_t^\infty L_x^{\frac{n(n-2)}{2(n-1)}}.
$$
Indeed, for
$F_0:=P_{> 1/4}|u_{lo}|^{\frac{4}{n-2}}$, we have $F_0\in \dot{\Lambda}^{\frac{n(n-2)}{2(n-1)}}_{\frac{4}{n-2}}$ iff
for all dyadic $N$'s we have
$$
N^{\frac{4}{n-2}}\|P_{N}F_0\|_{L_x^{\frac{n(n-2)}{2(n-1)}}}\lesssim \eta_2^{\frac{4}{n-2}}.
$$
Hence,
\begin{align*}
\||\nabla|^{\frac{2}{n-2}}F_0\|_{L_x^{\frac{n(n-2)}{2(n-1)}}}
&\lesssim \sum_{N>1/4}N^{\frac{2}{n-2}}\|P_N F_0\|_{L_x^{\frac{n(n-2)}{2(n-1)}}}
\lesssim \eta_2^{\frac{4}{n-2}}\sum_{N>1/4}N^{\frac{2}{n-2}}N^{-\frac{4}{n-2}}\\
&\lesssim \eta_2^{\frac{4}{n-2}}.
\end{align*}
Thus,
$$
\||\nabla|^{-\frac{2}{n-2}}(bP_{> 1/4}|u_{lo}|^{\frac{4}{n-2}})\|_{2,\frac{2n(n-2)}{n^2+n-10}}
\lesssim 2C_1\eta_1^{\frac{1}{4}}\eta_2^{\frac{4}{n-2}}\leq \frac{C_1}{100}\eta_1^{\frac{1}{4}}.
$$

Putting everything together, we find that
$$
\||\nabla|^{-\frac{2}{n-2}}b\|_{2,\frac{2n(n-2)}{n^2-3n-2}}
\leq C_1\eta_1^{\frac{1}{4}}.
$$

Therefore $T\in \Omega_1$. This concludes the proof of Proposition \ref{Scontrol}.
\end{proof}

\subsection{FLIM: the error terms} In this section we use the control on $u_{lo}$ and $u_{hi}$
that Proposition \ref{Scontrol} won us to bound the terms appearing on the right-hand side of Proposition \ref{etprop}.
For the rest of this section $n\geq6$ and all spacetime norms are taken on $I_0\times\R^n$.

Consider \eqref{et1}. Using \eqref{diff2} and H\"older, we estimate
\begin{align*}
\eqref{et1}
&\lesssim \eta_2\bigl\{\||u_{hi}|^2|u_{lo}|^{\frac{4}{n-2}}\chi_{\{|u_{hi}|\ll|u_{lo}|\}}\|_{L_{t,x}^1}
   +\||u_{hi}|^{\frac{n+2}{n-2}}u_{lo}\chi_{\{|u_{lo}|\ll|u_{hi}|\}}\|_{L_{t,x}^1}\bigr\} \\
&\lesssim \eta_2\||u_{hi}|^2|u_{lo}|^{\frac{4}{n-2}}\|_{L_{t,x}^1}
   \lesssim \eta_2\bigl\{\||g|^2|u_{lo}|^{\frac{4}{n-2}}\|_{L_{t,x}^1}+\||b|^2|u_{lo}|^{\frac{4}{n-2}}\|_{L_{t,x}^1}\bigr\} \\
&\lesssim \eta_2\bigl\{\|g\|^2_{2,\frac{2n}{n-2}}\|u_{lo}\|_{\infty,\frac{2n}{n-2}}^{\frac{4}{n-2}}
   +\|b\|^2_{\frac{2n}{n-2},\frac{2n^2}{(n+1)(n-2)}}\|u_{lo}\|_{\frac{2n}{n-2},\frac{4n^2}{(n+2)(n-2)}}^{\frac{4}{n-2}}\bigr\}.
\end{align*}
For $n\geq 6$, an application of Bernstein yields
$$
\|u_{lo}\|_{\frac{2n}{n-2},\frac{4n^2}{(n+2)(n-2)}}\lesssim \|u_{lo}\|_{\frac{2n}{n-2},\frac{2n^2}{(n-2)^2}}\lesssim \|u_{lo}\|_{\dot{S}^1}
$$
and hence, by Proposition \ref{Scontrol} and \eqref{lowfreqsmall},
\begin{align*}
\eqref{et1}
&\lesssim \eta_2\bigl\{\|g\|_{\dot{S}^0}^2\|u_{lo}\|_{\infty ,\frac{2n}{n-2}}^{\frac{4}{n-2}}
   +\|b\|^2_{\frac{2n}{n-2},\frac{2n^2}{(n+1)(n-2)}}\|u_{lo}\|_{\dot{S}^1}^{\frac{4}{n-2}}\bigr\}
\ll \eta_1.
\end{align*}

Consider next the error term \eqref{et2}. Replacing the projection $P_{lo}$ by the positive-kernel operator $P_{lo}'$ (see Subsection~1.2 for the
definition and properties of $P_{lo}'$) and splitting $u_{hi}=g+b$, we estimate
\begin{align*}
\eqref{et2}
&\lesssim \eta_2\bigl\{ \|g P_{lo}'(|g|^{\frac{n+2}{n-2}})\|_{L^1_{t,x}}+ \|g P_{lo}'(|b|^{\frac{n+2}{n-2}})\|_{L^1_{t,x}}
   +\|b P_{lo}'(|g|^{\frac{n+2}{n-2}})\|_{L^1_{t,x}}\\
&\quad + \|b P_{lo}'(|b|^{\frac{n+2}{n-2}})\|_{L^1_{t,x}}\bigr\}.
\end{align*}
By Proposition \ref{Scontrol} and Bernstein, we estimate
\begin{align*}
\|g P_{lo}'(|g|^{\frac{n+2}{n-2}})\|_{L^1_{t,x}}
  &\lesssim \|g\|^2_{2,\frac{2n}{n-2}}\|g\|_{\infty,\frac{2n}{n-2}}^{\frac{4}{n-2}}
  \lesssim C_1^{\frac{2n}{n-2}}\eta_2^{\frac{4}{n-2}},\\
\|g P_{lo}'(|b|^{\frac{n+2}{n-2}})\|_{L^1_{t,x}}
  &\lesssim \|g\|_{2,\frac{2n}{n-2}}\|P_{lo}'(|b|^{\frac{n+2}{n-2}})\|_{2,\frac{2n}{n+2}}
    \lesssim \|g\|_{\dot{S}^0}\|P_{lo}'(|b|^{\frac{n+2}{n-2}})\|_{2,\frac{2n}{n+3}}\\
  &\lesssim \|g\|_{\dot{S}^0}\|b\|_{\frac{2(n+2)}{n-2},\frac{2n(n+2)}{(n-2)(n+3)}}^{\frac{n+2}{n-2}}
    \lesssim C_1^{\frac{2n}{n-2}}\eta_2^{\frac{2}{n-2}}\eta_1^{\frac{1}{4}},\\
\|b P_{lo}'(|g|^{\frac{n+2}{n-2}})\|_{L^1_{t,x}}
  &\lesssim \|b\|_{\frac{2n}{n-2},\frac{2n^2}{(n+1)(n-2)}}\|P_{lo}'(|g|^{\frac{n+2}{n-2}})\|_{\frac{2n}{n+2},\frac{2n^2}{n^2+n+2}} \\
  &\lesssim \|b\|_{\frac{2n}{n-2},\frac{2n^2}{(n+1)(n-2)}}\|P_{lo}'(|g|^{\frac{n+2}{n-2}})\|_{\frac{2n}{n+2},\frac{2n^2(n-2)}{(n^2-2n+4)(n+2)}} \\
  &\lesssim \|b\|_{\frac{2n}{n-2},\frac{2n^2}{(n+1)(n-2)}}\|g\|_{\frac{2n}{n-2},\frac{2n^2}{n^2-2n+4}}^{\frac{n+2}{n-2}}\\
  &\lesssim \|b\|_{\frac{2n}{n-2},\frac{2n^2}{(n+1)(n-2)}}\|g\|_{\dot{S}^0}^{\frac{n+2}{n-2}}
    \lesssim C_1^{\frac{2n}{n-2}}\eta_1^{\frac{n-2}{4n}}\eta_2^{\frac{2(n+2)}{(n-2)^2}},\\
\|b P_{lo}'(|b|^{\frac{n+2}{n-2}})\|_{L^1_{t,x}}
  &\lesssim \|b\|_{\frac{2n}{n-2},\frac{2n^2}{(n+1)(n-2)}}\|P_{lo}'(|b|^{\frac{n+2}{n-2}})\|_{\frac{2n}{n+2},\frac{2n^2}{n^2+n+2}} \\
  &\lesssim \|b\|_{\frac{2n}{n-2},\frac{2n^2}{(n+1)(n-2)}}\|P_{lo}'(|b|^{\frac{n+2}{n-2}})\|_{\frac{2n}{n+2},\frac{2n^2}{(n+1)(n+2)}} \\
  &\lesssim \|b\|_{\frac{2n}{n-2},\frac{2n^2}{(n+1)(n-2)}}^{\frac{2n}{n-2}}
    \lesssim C_1^{\frac{2n}{n-2}}\eta_1^{\frac{1}{2}}.
\end{align*}
Hence,
$$
\eqref{et2}
   \lesssim \eta_2 C_1^{\frac{2n}{n-2}}\bigl\{\eta_2^{\frac{4}{n-2}}+\eta_2^{\frac{2}{n-2}}\eta_1^{\frac{1}{4}}
       +\eta_1^{\frac{n-2}{4n}}\eta_2^{\frac{2(n+2)}{(n-2)^2}}+\eta_1^{\frac{1}{2}}\bigr\}
   \ll\eta_1.
$$

We turn next to the error term \eqref{et3}. Decomposing again $u_{hi}=g+b$, we estimate
$$
\eqref{et3}
\lesssim \eta_2\bigl\{\|g P_{hi}\bigl(|u_{lo}|^{\frac{4}{n-2}}u_{lo}\bigr)\|_{L_{t,x}^1} +\|b P_{hi}\bigl(|u_{lo}|^{\frac{4}{n-2}}u_{lo}\bigr)\|_{L_{t,x}^1}\bigr\}.
$$
By H\"older, Bernstein, and Proposition \ref{Scontrol}, we estimate
\begin{align*}
\|g P_{hi}\bigl(|u_{lo}|^{\frac{4}{n-2}}u_{lo}\bigr)\|_{L_{t,x}^1}
&\lesssim \|g\|_{2,\frac{2n}{n-2}}\|\nabla u_{lo}\|_{2,\frac{2n}{n-2}}\|u_{lo}\|_{\infty,\frac{2n}{n-2}}^{\frac{4}{n-2}}\\
&\lesssim \|g\|_{\dot{S}^0} \| u_{lo}\|_{\dot{S}^1} \|u_{lo}\|_{\infty,\frac{2n}{n-2}}^{\frac{4}{n-2}}
  \lesssim C_1^2\eta_2^{\frac{2}{n-2}} \eta_2^{\frac{4}{(n-2)^2}} \eta_2^{\frac{4}{n-2}}, \\
\|b P_{hi}\bigl(|u_{lo}|^{\frac{4}{n-2}}u_{lo}\bigr)\|_{L_{t,x}^1}
&\lesssim \|b\|_{\frac{2n}{n-2},\frac{2n^2}{(n+1)(n-2)}}\|\nabla u_{lo}\|_{\frac{2n}{n-2},\frac{2n^2}{n^2-2n+4}}
    \|u_{lo}\|_{\frac{2n}{n-2},\frac{8n^2}{(3n-2)(n-2)}}^{\frac{4}{n-2}}.
\end{align*}
For $n\geq 6$, Bernstein dictates
\begin{align}\label{uloinS^1}
\|u_{lo}\|_{\frac{2n}{n-2},\frac{8n^2}{(3n-2)(n-2)}}
   \lesssim \|u_{lo}\|_{\frac{2n}{n-2},\frac{2n^2}{(n-2)^2}}
   \lesssim \|u_{lo}\|_{\dot{S}^1}
\end{align}
and hence
$$
\|b P_{hi}\bigl(|u_{lo}|^{\frac{4}{n-2}}u_{lo}\bigr)\|_{L_{t,x}^1}
  \lesssim \|b\|_{\frac{2n}{n-2},\frac{2n^2}{(n+1)(n-2)}}\|u_{lo}\|_{\dot{S}^1}^{\frac{n+2}{n-2}}
  \lesssim C_1^{\frac{2n}{n-2}}\eta_1^{\frac{n-2}{4n}}\eta_2^{\frac{4}{(n-2)^2}\cdot\frac{n+2}{n-2}}.
$$
Thus $\eqref{et3}\ll \eta_1$.

Consider now \eqref{et4}. Decomposing $u_{hi}=g+b$ and applying H\"older, we estimate
$$
\eqref{et4}
  \lesssim \eta_2^2\|\nabla u_{lo}\|_{2,\frac{2n}{n-2}}
    \bigl\{\|g |u_{lo}|^{\frac{4}{n-2}}\|_{2,\frac{2n}{n+2}}+\|b|u_{lo}|^{\frac{4}{n-2}}\|_{2,\frac{2n}{n+2}}\bigr\}.
$$
Using again Proposition \ref{Scontrol}, we estimate
\begin{align*}
\|g |u_{lo}|^{\frac{4}{n-2}}\|_{2,\frac{2n}{n+2}}
&\lesssim \|g\|_{2,\frac{2n}{n-2}}\|u_{lo}\|_{\infty,\frac{2n}{n-2}}^{\frac{4}{n-2}}
  \lesssim \|g\|_{\dot{S}^0}\|u_{lo}\|_{\infty,\frac{2n}{n-2}}^{\frac{4}{n-2}}
  \lesssim C_1\eta_2^{\frac{6}{n-2}},
\end{align*}
and
\begin{align}\label{bulo^4/n-2}
\|b|u_{lo}|^{\frac{4}{n-2}}\|_{2,\frac{2n}{n+2}}
&\lesssim \|b\|_{\frac{2n}{n-2},\frac{2n^2}{(n+1)(n-2)}} \|u_{lo}\|_{\frac{4n}{n-2},\frac{8n^2}{(3n+2)(n-2)}}^{\frac{4}{n-2}}\notag \\
&\lesssim \|b\|_{\frac{2n}{n-2},\frac{2n^2}{(n+1)(n-2)}} \|u_{lo}\|_{\dot{S}^1}^{\frac{4}{n-2}}
 \lesssim C_1^{\frac{n+2}{n-2}}\eta_1^{\frac{n-2}{4n}}\eta_2^{\frac{16}{(n-2)^3}},
\end{align}
where we also used the fact that for $n\geq 6$, an application of Bernstein yields
\begin{align}\label{ulo in 4n/n-2}
\|u_{lo}\|_{\frac{4n}{n-2},\frac{8n^2}{(3n+2)(n-2)}}
   \lesssim \|u_{lo}\|_{\frac{4n}{n-2},\frac{2n^2}{n^2-3n+2}}
   \lesssim \|u_{lo}\|_{\dot{S}^1},
\end{align}
Thus
$$
\eqref{et4}
  \lesssim \eta_2^2 \eta_2^{\frac{4}{(n-2)^2}}\bigl\{C_1^2\eta_2^{\frac{6}{n-2}} +C_1^{\frac{2n}{n-2}}\eta_1^{\frac{n-2}{4n}}\eta_2^{\frac{16}{(n-2)^3}}\bigr\}
  \ll \eta_1.
$$

We turn now towards the error term \eqref{et5}, which we estimate by
$$
\eqref{et5}
 \lesssim \eta_2^2 \bigl\{\|\nabla u_{lo} |g|^{\frac{n+2}{n-2}}\|_{L_{t,x}^1}+\|\nabla u_{lo} |b|^{\frac{n+2}{n-2}}\|_{L_{t,x}^1}\bigr\}.
$$
By H\"older and Proposition \ref{Scontrol},
\begin{align*}
\|\nabla u_{lo} |g|^{\frac{n+2}{n-2}}\|_{L_{t,x}^1}
&\lesssim \|\nabla u_{lo}\|_{2,\frac{2n}{n-2}}\|g\|_{2,\frac{2n}{n-2}}\|g\|_{\infty,\frac{2n}{n-2}}^{\frac{4}{n-2}}
  \lesssim \|u_{lo}\|_{\dot{S}^1}\|g\|_{\dot{S}^0}\|g\|_{\dot{S}^1}^{\frac{4}{n-2}}\\
&\lesssim C_1^{\frac{2n}{n-2}}\eta_2^{\frac{4}{(n-2)^2}}\eta_2^{\frac{2}{n-2}},\\
\|\nabla u_{lo} |b|^{\frac{n+2}{n-2}}\|_{L_{t,x}^1}
&\lesssim \|\nabla u_{lo}\|_{\frac{2n}{n-2},\frac{2n^2}{n^2-3n-2}}\|b\|^{\frac{n+2}{n-2}}_{\frac{2n}{n-2},\frac{2n^2}{(n+1)(n-2)}}\\
&\lesssim C_1^{\frac{2n}{n-2}}\eta_2^{\frac{4}{(n-2)^2}}\eta_1^{\frac{n-2}{4n}\cdot\frac{n+2}{n-2}},
\end{align*}
where we applied Bernstein to estimate
$$
\|\nabla u_{lo}\|_{\frac{2n}{n-2},\frac{2n^2}{n^2-3n-2}}
 \lesssim \|\nabla u_{lo}\|_{\frac{2n}{n-2},\frac{2n^2}{n^2-2n+4}}\lesssim \|u_{lo}\|_{\dot{S}^1}.
$$
Hence,
$$
\eqref{et5}
  \lesssim \eta_2^2C_1^{\frac{2n}{n-2}}\bigl\{\eta_2^{\frac{4}{(n-2)^2}}\eta_2^{\frac{2}{n-2}}+\eta_2^{\frac{4}{(n-2)^2}}\eta_1^{\frac{n+2}{4n}}\bigr\}
  \ll \eta_1.
$$

Consider now the error term \eqref{et6}. We estimate
\begin{align}\label{et6split}
\|\nabla P_{lo}\bigl(|u|^{\frac{4}{n-2}}u\bigr)u_{hi}\|_{1,1}
&\lesssim \|\nabla P_{lo}\bigl(|u_{hi}|^{\frac{4}{n-2}}u_{hi}\bigr)u_{hi}\|_{1,1}+\|\nabla P_{lo}\bigl(|u_{lo}|^{\frac{4}{n-2}}u_{lo}\bigr)u_{hi}\|_{1,1}\notag\\
&\quad + \|\nabla P_{lo}\bigl(|u|^{\frac{4}{n-2}}u-|u_{lo}|^{\frac{4}{n-2}}u_{lo}-|u_{hi}|^{\frac{4}{n-2}}u_{hi}\bigr)u_{hi}\|_{1,1}.
\end{align}
Using Bernstein to drop the derivative in front of the projection $P_{lo}$, we recognize in the first term on the
right-hand side of \eqref{et6split} the error term \eqref{et2}. Hence, by the previous computations, we have
\begin{align}
\eta_2^2\|\nabla P_{lo}\bigl(|u_{hi}|^{\frac{4}{n-2}}u_{hi}\bigr)u_{hi}\|_{L_{t,x}^1}\leq \eta_2 \eqref{et2} \ll\eta_1. \label{et6-1}
\end{align}
To estimate the second term on the right-hand side of \eqref{et6split}, we decompose $u_{hi}=g+b$ and use
\eqref{lowfreqsmall}, \eqref{uloinS^1}, and Proposition \ref{Scontrol} to control the two resulting terms as follows:
\begin{align*}
\|\nabla P_{lo}\bigl(|u_{lo}|^{\frac{4}{n-2}}u_{lo}\bigr)g\|_{L_{t,x}^1}
&\lesssim \|\nabla u_{lo}\|_{2,\frac{2n}{n-2}}\|u_{lo}\|_{\infty,\frac{2n}{n-2}}^{\frac{4}{n-2}}\|g\|_{2,\frac{2n}{n-2}}
  \lesssim \eta_2^{\frac{4}{n-2}}\|u_{lo}\|_{\dot{S}^1}\|g\|_{\dot{S}^0}\\
&\lesssim C_1^2\eta_2^{\frac{6}{n-2}}\eta_2^{\frac{4}{(n-2)^2}},\\
\|\nabla P_{lo}\bigl(|u_{lo}|^{\frac{4}{n-2}}u_{lo}\bigr)b\|_{L_{t,x}^1}
&\lesssim \|\nabla u_{lo}\|_{\frac{2n}{n-2},\frac{2n^2}{n^2-2n+4}}\| u_{lo}\|_{\frac{2n}{n-2},\frac{8n^2}{(3n-2)(n-2)}}^{\frac{4}{n-2}}\|b\|_{\frac{2n}{n-2},\frac{2n^2}{(n+1)(n-2)}}\\
&\lesssim \|u_{lo}\|_{\dot{S}^1}^{\frac{n+2}{n-2}}\|b\|_{\frac{2n}{n-2},\frac{2n^2}{(n+1)(n-2)}}
  \lesssim C_1^{\frac{2n}{n-2}}\eta_1^{\frac{n-2}{4n}}\eta_2^{\frac{4(n+2)}{(n-2)^3}}.
\end{align*}
Thus,
\begin{align}\label{et6-2}
\eta_2^2\|\nabla P_{lo}\bigl(|u_{lo}|^{\frac{4}{n-2}}u_{lo}\bigr)u_{hi}\|_{L_{t,x}^1}\ll \eta_1.
\end{align}
To estimate the third term on the right-hand side of \eqref{et6split}, we first use Bernstein to drop the derivative in
front of $P_{lo}$ and then replace the projection $P_{lo}$ by the positive-kernel operator $P_{lo}'$ sharing the same $L_x^p$-mapping
and Bernstein properties as $P_{lo}$; using \eqref{diff2}, we obtain the bound
\begin{align*}
\|P_{lo}'\bigl(|u_{lo}|^{\frac{4}{n-2}}|u_{hi}|&\chi_{\{|u_{hi}|\leq |u_{lo}|\}}\bigr)u_{hi}\|_{L_{t,x}^1}
+\|P_{lo}'\bigl(|u_{hi}|^{\frac{4}{n-2}}|u_{lo}|\chi_{\{|u_{lo}|< |u_{hi}|\}}\bigr)u_{hi}\|_{L_{t,x}^1}\\
&\lesssim \|P_{lo}'\bigl(|u_{lo}|^{\frac{4}{n-2}}|u_{hi}|\bigr)u_{hi}\|_{L_{t,x}^1}.
\end{align*}
Decomposing $u_{hi}=g+b$ and using \eqref{lowfreqsmall}, \eqref{bulo^4/n-2}, \eqref{ulo in 4n/n-2},
and Proposition~\ref{Scontrol}, we further estimate the third term on the right-hand side of \eqref{et6split} by
the sum of the following four terms:
\begin{align*}
\|P_{lo}'\bigl(|u_{lo}|^{\frac{4}{n-2}}|g|\bigr)g\|_{L_{t,x}^1}
&\lesssim \|g\|_{2,\frac{2n}{n-2}}^2\|u_{lo}\|_{\infty, \frac{2n}{n-2}}^{\frac{4}{n-2}}
\lesssim \eta_2^{\frac{4}{n-2}}\|g\|_{\dot{S}^0}^2
\lesssim C_1^2\eta_2^{\frac{8}{n-2}}\\
\|P_{lo}'\bigl(|u_{lo}|^{\frac{4}{n-2}}|b|\bigr)g\|_{L_{t,x}^1}
&\lesssim \|g\|_{2,\frac{2n}{n-2}}\|b |u_{lo}|^{\frac{4}{n-2}}\|_{2,\frac{2n}{n+2}}
\lesssim C_1^{\frac{2n}{n-2}} \eta_2^{\frac{2}{n-2}}\eta_1^{\frac{n-2}{4n}}\eta_2^{\frac{16}{(n-2)^3}}\\
\|P_{lo}'\bigl(|u_{lo}|^{\frac{4}{n-2}}|g|\bigr)b\|_{L_{t,x}^1}
&\lesssim \|b\|_{\frac{2n}{n-2},\frac{2n^2}{(n+1)(n-2)}}\|g\|_{2,\frac{2n}{n-2}}\|u_{lo}\|_{\frac{4n}{n-2},\frac{8n^2}{(3n+2)(n-2)}}^{\frac{4}{n-2}}\\
&\lesssim \|b\|_{\frac{2n}{n-2},\frac{2n^2}{(n+1)(n-2)}}\|g\|_{\dot{S}^0}\|u_{lo}\|_{\dot{S}^1}^{\frac{4}{n-2}}\\
&\lesssim C_1^{\frac{2n}{n-2}} \eta_1^{\frac{n-2}{4n}}\eta_2^{\frac{2}{n-2}}\eta_2^{\frac{16}{(n-2)^3}}\\
\|P_{lo}'\bigl(|u_{lo}|^{\frac{4}{n-2}}|b|\bigr)b\|_{L_{t,x}^1}
&\lesssim \|P_{lo}'\bigl(|u_{lo}|^{\frac{4}{n-2}}|b|\bigr)b\|_{1,\frac{n^2}{n^2+n-6}}\\
&\lesssim\|b\|_{\frac{2n}{n-2},\frac{2n^2}{(n+1)(n-2)}}^2\|u_{lo}\|_{\frac{2n}{n-2},\frac{2n^2}{(n-2)^2}}^{\frac{4}{n-2}}\\
&\lesssim C_1^{\frac{2n}{n-2}} \eta_1^{\frac{n-2}{2n}}\eta_2^{\frac{16}{(n-2)^3}}.
\end{align*}
Hence,
\begin{align}\label{et6-3}
\eta_2^2\|\nabla P_{lo}\bigl(|u|^{\frac{4}{n-2}}u-|u_{lo}|^{\frac{4}{n-2}}u_{lo}-|u_{hi}|^{\frac{4}{n-2}}u_{hi}\bigr)u_{hi}\|_{L_{t,x}^1}\ll\eta_1.
\end{align}
Collecting \eqref{et6split} through \eqref{et6-3}, we obtain
$$\eqref{et6}\ll \eta_1.$$

We turn now to the error terms \eqref{et7} through \eqref{et9}. We notice they are of the form
$\langle |u_{hi}|^2*\frac{1}{|x|},f\rangle$ where
$$
f=\left\{ \begin{aligned}
|u_{lo}|^{\frac{n+2}{n-2}}|u_{hi}| \ \ &{\rm in\ } \eqref{et7}, \\
|u_{lo}||u_{hi}|^{\frac{n+2}{n-2}} \ \ &{\rm in\ } \eqref{et8}, \\
|P_{lo}\bigl(|u_{hi}|^{\frac{4}{n-2}}u_{hi}\bigr)u_{hi}| \ \ &{\rm in\ } \eqref{et9}. \end{aligned}\right.
$$

Let us first note that $u_{hi}\in L_t^3L_x^{\frac{6n}{3n-4}}$. Indeed, by Strichartz,
\begin{align}\label{gL^3!!}
\|g\|_{3,\frac{6n}{3n-4}}\lesssim \|g\|_{\dot{S}^0} \lesssim C_1 \eta_2^{\frac{2}{n-2}},
\end{align}
while
\begin{align}\label{bL^3!!}
\|b\|_{3,\frac{6n}{3n-4}}
\lesssim\|b\|_{\frac{2n}{n-2},\frac{2n^2}{(n+1)(n-2)}}^{\frac{2n}{3(n-2)}}\|b\|_{\infty,\frac{2n}{n-2}}^{\frac{n-6}{3(n-2)}}
\lesssim C_1\eta_1^{\frac{1}{6}}.
\end{align}
Thus,
\begin{align}\label{uhiL^3!!}
\|u_{hi}\|_{3,\frac{6n}{3n-4}}\lesssim C_1 \eta_1^{\frac{1}{6}}.
\end{align}

As $u_{hi}\in L_t^3L_x^{\frac{6n}{3n-4}}$ and $u_{hi}\in L_t^{\infty}L_x^2$, we get
$|u_{hi}|^2\in L_t^3L_x^{\frac{3n}{3n-2}}$.
Therefore, by Hardy-Littlewood-Sobolev, $|u_{hi}|^2*\frac{1}{|x|}\in L_t^3L_x^{3n}$. Moreover,
\begin{align}\label{?!?!}
\bigl\langle |u_{hi}|^2*\frac{1}{|x|},f\bigr\rangle
&\lesssim \bigl\||u_{hi}|^2*|x|^{-1}\bigr\|_{3,3n}\|f\|_{\frac{3}{2},\frac{3n}{3n-1}}
 \lesssim \|u_{hi}\|_{3,\frac{6n}{3n-4}}\|u_{hi}\|_{\infty,2}\|f\|_{\frac{3}{2},\frac{3n}{3n-1}}\notag\\
&\lesssim C_1 \eta_1^{\frac{1}{6}}\eta_2\|f\|_{\frac{3}{2},\frac{3n}{3n-1}}.
\end{align}

Consider the case of \eqref{et7}, that is, $f=|u_{lo}|^{\frac{n+2}{n-2}}|u_{hi}|$. By H\"older, \eqref{slf},
and \eqref{uhiL^3!!}, we estimate
\begin{align*}
\||u_{lo}|^{\frac{n+2}{n-2}}u_{hi}\|_{\frac{3}{2},\frac{3n}{3n-1}}
&\lesssim \|u_{hi}\|_{3,\frac{6n}{3n-4}}\|u_{lo}\|_{\infty,\frac{2n}{n-2}}^{\frac{4}{n-2}}\|u_{lo}\|_{3,\frac{6n}{3n-10}}\\
&\lesssim C_1 \eta_1^{\frac{1}{6}}\eta_2^{\frac{4}{n-2}}\|u_{lo}\|_{\dot{S}^{1}}
 \lesssim C_1^2\eta_1^{\frac{1}{6}}\eta_2^{\frac{4}{n-2}}\eta_2^{\frac{4}{(n-2)^2}}.
\end{align*}
Thus, by \eqref{?!?!} and the above computation, $ \eqref{et7}\ll \eta_1$.

Consider next the error term \eqref{et8}, that is, $f=|u_{lo}||u_{hi}|^{\frac{n+2}{n-2}}$. By \eqref{slf},
\eqref{uhiL^3!!}, and the conservation of energy, we estimate
\begin{align*}
\|u_{lo}|u_{hi}|^{\frac{n+2}{n-2}}\|_{\frac{3}{2},\frac{3n}{3n-1}}
&\lesssim \|u_{hi}\|_{3,\frac{6n}{3n-4}}\|u_{hi}\|_{\infty,\frac{2n}{n-2}}^{\frac{4}{n-2}}\|u_{lo}\|_{3,\frac{6n}{3n-10}}
 \lesssim C_1 \eta_1^{\frac{1}{6}}\|u_{lo}\|_{\dot{S}^{1}}\\
&\lesssim C_1^2 \eta_1^{\frac{1}{6}}\eta_2^{\frac{4}{(n-2)^2}}.
\end{align*}
Hence, by \eqref{?!?!} and the above computation, $\eqref{et8}\ll \eta_1$.

The last error term left to consider is \eqref{et9}; in this case we use Bernstein, \eqref{uhiL^3!!}, and the
conservation of energy to estimate
\begin{align*}
\|P_{lo}\bigl(|u_{hi}|^{\frac{4}{n-2}}u_{hi}\bigr)u_{hi}\|_{\frac{3}{2},\frac{3n}{3n-1}}
&\lesssim \|u_{hi}\|_{3,\frac{6n}{3n-4}}\|P_{lo}\bigl(|u_{hi}|^{\frac{4}{n-2}}u_{hi}\bigr)\|_{3,\frac{6n}{3n+2}}\\
&\lesssim C_1 \eta_1^{\frac{1}{6}}\|P_{lo}\bigl(|u_{hi}|^{\frac{4}{n-2}}u_{hi}\bigr)\|_{3,\frac{6n}{3n+8}}\\
&\lesssim C_1 \eta_1^{\frac{1}{6}}\|u_{hi}\|_{3,\frac{6n}{3n-4}}\|u_{hi}\|_{\infty,\frac{2n}{n-2}}^{\frac{4}{n-2}}\\
&\lesssim C_1^2 \eta_1^{\frac{1}{3}}.
\end{align*}
Thus, considering \eqref{?!?!}, we obtain $\eqref{et9}\ll\eta_1$.

Hence all the error terms \eqref{et1} through \eqref{et9} are bounded by $\eta_1$. Upon rescaling,
this concludes the proof of Proposition \ref{flim}.

As a consequence of Proposition~\ref{Scontrol} and scaling, we obtain the following:
\begin{corollary}\label{cor uhi}
Let $n\geq 6$, $u$ a minimal energy blowup solution to \eqref{schrodinger equation}, and $N_*<c(\eta_2)N_{min}$.
Then, we can decompose $P_{\geq N_*}u=g+b$ such that
\begin{align}
&\|g\|_{\dot{S}^{0}(I_0\times \R^n)}\lesssim \eta_{2}^{\frac{2}{n-2}}N_*^{-1}, \label{gS0 scaled} \\
&\|g\|_{\dot{S}^{1}(I_0\times \R^n)}\lesssim 1, \label{gS1 scaled} \\
&\||\nabla|^{-\frac{2}{n-2}}b\|_{L_t^2L_x^{\frac{2n(n-2)}{n^2-3n-2}}(I_0\times \R^n)}\lesssim \eta_1^{\frac{1}{4}}N_*^{-\frac{3}{2}}. \label{bS scaled}
\end{align}
Moreover, under scaling, \eqref{b in 2n/n-2}, \eqref{b in 2(n+2)/(n-2)}, and \eqref{uhiL^3!!} become the following
\begin{gather}
\|b\|_{L_t^{\frac{2n}{n-2}}L_x^{\frac{2n^2}{(n+1)(n-2)}}(I_0\times \R^n)}\lesssim \bigl(\eta_1^{\frac{1}{4}}N_*^{-\frac{3}{2}}\bigr)^{\frac{n-2}{n}},\\
\|b\|_{L_t^{\frac{2(n+2)}{n-2}}L_x^{\frac{2n(n+2)}{(n-2)(n+3)}}(I_0\times \R^n)}\lesssim \bigl(\eta_1^{\frac{1}{4}}N_*^{-\frac{3}{2}}\bigr)^{\frac{n-2}{n+2}},\label{b in 2(n+2)/(n-2) scaled}\\
\|P_{\geq N_*}u\|_{L_t^3L_x^{\frac{6n}{3n-4}}(I_0\times\R^n)}\lesssim \eta_1^{\frac{1}{6}}N_*^{-1}.\label{uhiL^3!!!}
\end{gather}
\end{corollary}

\begin{remark} \label{rem n=5}
In dimension $n=5$, the analogue of \eqref{uhiL^3!!!} is
\begin{align}\label{n=5}
\|P_{\geq N_*}u\|_{L_t^3L_x^{\frac{30}{11}}(I_0\times\R^5)}\lesssim \eta_1^{\frac{1}{3}}N_*^{-1}.
\end{align}
For details see \cite{my thesis}.
\end{remark}

%
%
%
%

\section{Preventing Energy Evacuation}
The purpose of this section is to prove

\begin{proposition}[Energy cannot evacuate to high frequencies]\label{energyevac} Suppose that $u$ is a minimal energy
blowup solution to \eqref{schrodinger equation}. Then for all $t\in I_{0}$,
\begin{align}
N(t)\leq C(\eta_4)N_{min}. \label{freqbound}
\end{align}
\end{proposition}

\subsection{The setup}
We normalize so that $N_{min}=1$. As $N(t) \in 2^{\Z}$, there exists $t_{min} \in I_0$ such that $N(t_{min})=N_{min}=1$.

At time $t=t_{min}$, we have a considerable amount of mass at medium frequencies:
\begin{align}
\|P_{c(\eta_0)<\cdot<C(\eta_0)}u(t_{min})\|_{L^2_x}
   \gtrsim c(\eta_0)\|P_{c(\eta_0)<\cdot<C(\eta_0)}u(t_{min})\|_{\dot{H}^1_x}
   \sim c(\eta_0). \label{massmedfreq}
\end{align}
However, by Bernstein, there is not much mass at frequencies higher than $C(\eta_0)$:
$$
\|P_{>C(\eta_0)}u(t_{min})\|_{L^2_x} \lesssim c(\eta_0).
$$

Let's assume for a contradiction that there exists $t_{evac} \in I_0$ such that $N(t_{evac})\gg C(\eta_4)$.
By time reversal symmetry we may assume $t_{evac} > t_{min}$. As for every $\eta_4 \leq \eta \leq \eta_0$ and all
$t \in I_0$, $\|P_{<c(\eta)N(t)}u(t)\|_{\dot{H}^1_x} \leq \eta$, we see that by choosing $C(\eta_4)$ sufficiently large,
at time $t=t_{evac}$ there is very little energy at low and medium frequencies:
\begin{align}
\|P_{< \eta_4^{-1}}u(t_{evac})\|_{\dot{H}^1_x} \leq \eta_4. \label{goodcontrollowfreq}
\end{align}

We define $u_{lo}=P_{<\eta_{3}^{10n}}u$ and $u_{hi}=P_{\geq \eta_{3}^{10n}}u$. Then by \eqref{massmedfreq},
\begin{align}
\|u_{hi}(t_{min})\|_{L^2_x} \geq \eta_1. \label{highmass}
\end{align}
Suppose we could show that a big portion of the mass sticks around until time $t=t_{evac}$, i.e.,
\begin{align}
\|u_{hi}(t_{evac})\|_{L^2_x} \geq \frac{1}{2}\eta_1. \label{highmassasmp}
\end{align}
Then, since by Bernstein
$$
\|P_{>C(\eta_1)}u_{hi}(t_{evac})\|_{L^2_x}\leq c(\eta_1),
$$
the triangle inequality would imply
$$
\|P_{\leq C(\eta_1)}u_{hi}(t_{evac})\|_{L^2_x} \geq \frac{1}{4}\eta_1.
$$
Another application of Bernstein would give
$$
\|P_{\leq C(\eta_1)}u(t_{evac})\|_{\dot{H}^1_x} \gtrsim c(\eta_1, \eta_3),
$$
which would contradict \eqref{goodcontrollowfreq} if $\eta_4$ were chosen sufficiently small.

It therefore remains to show \eqref{highmassasmp}. In order to prove it we assume that there exists a time $t_{*}$
such that $t_{min} \leq t_{*} \leq t_{evac}$ and
\begin{align}
\inf_{t_{min} \leq t \leq t_{*}}\|u_{hi}(t)\|_{L^2_x} \geq \frac{1}{2} \eta_1. \label{contasmp}
\end{align}
We will show that this can be bootstrapped to
\begin{align}
\inf_{t_{min} \leq t \leq t_{*}}\|u_{hi}(t)\|_{L^2_x} \geq \frac{3}{4} \eta_1. \label{bootstrapmass}
\end{align}
Hence, $\{t_{*} \in [t_{min}, t_{evac}] : \eqref{contasmp} \text{ holds} \}$ is both open and closed in
$[t_{min}, t_{evac}]$ and \eqref{highmassasmp} holds.

In order to show that \eqref{contasmp} implies \eqref{bootstrapmass}, we treat the $L^2_x$-norm of $u_{hi}$
as an almost conserved quantity. Define
$$
L(t)=\int_{\R^n}|u_{hi}(t,x)|^2dx.
$$
By \eqref{highmass} we have $L(t_{min}) \geq \eta_1^2$. Hence,
by the Fundamental Theorem of Calculus it suffices to show that
$$
\int_{t_{min}}^{t_{*}}|\partial_{t}L(t)|dt \leq \frac{1}{100}\eta_1^2.
$$
As
\begin{align*}
\partial_{t}L(t)
  &=2 \int_{\R^n} \{P_{hi}(|u|^{\frac{4}{n-2}}u), u_{hi}\}_{m}dx \\
  &=2 \int_{\R^n} \{P_{hi}(|u|^{\frac{4}{n-2}}u)-|u_{hi}|^{\frac{4}{n-2}}u_{hi}, u_{hi}\}_{m}dx,
\end{align*}
we need to show
\begin{align}
\int_{t_{min}}^{t_{*}} \Bigl| \int_{\R^n} \{P_{hi}(|u|^{\frac{4}{n-2}}u)-|u_{hi}|^{\frac{4}{n-2}}u_{hi}, u_{hi}\}_{m}dx\Bigr| dt \leq \frac{1}{100}\eta_1^2. \label{last}
\end{align}

In order to prove \eqref{last} we need to control the various interactions between low, medium, and high frequencies.
In the next section we will develop the tools that will make this goal possible.

\subsection{Spacetime estimates on low, medium, and high frequencies}

Remember that the frequency-localized interaction Morawetz inequality implies that for $N<c(\eta_2)N_{min}$,
\begin{align}
\int_{t_{min}}^{t_{evac}} \int_{\R^n}\int_{\R^n} \frac{|P_{\geq N}u(t,y)|^2|P_{\geq N}u(t,x)|^2}{|x-y|^3} dxdydt \lesssim \eta_1 N^{-3}. \label{conseqflim}
\end{align}

This estimate is useful for medium and high frequencies; however it is extremely bad for low frequencies since
$N^{-3}$ gets increasingly larger as $N \rightarrow 0$. We therefore need to develop better estimates in this case.
As $u_{\leq \eta_3}$ has extremely small energy at $t=t_{evac}$ (see \eqref{goodcontrollowfreq}),
we expect it to have small energy at all times in $[t_{min}, t_{evac}]$. Of course, there is energy leaking from
the high frequencies to the low frequencies, but the interaction Morawetz estimate limits this leakage. Indeed, we have

\begin{proposition}\label{S1lf}
Under the assumptions above,
\begin{align}
\|P_{\leq N}u\|_{\dot{S}^1([t_{min}, t_{evac}]\times \R^n)}
    \lesssim \eta_4 + \max\{\eta_3^{-\frac{3}{2}}N^{\frac{3}{2}},\eta_3^{-\frac{n+2}{n-2}} N^{\frac{n+2}{n-2}}\}\label{S^1lf}
\end{align}
for all $N \leq \eta_3$.
\end{proposition}

\begin{remark}
One should think of the $\eta_4$ factor on the right-hand side of \eqref{S^1lf} as the energy coming from
the low modes of $u(t_{evac})$, and the
$\max\{\eta_3^{-\frac{3}{2}}N^{\frac{3}{2}},\eta_3^{-\frac{n+2}{n-2}} N^{\frac{n+2}{n-2}}\}$
term as the energy coming from the high frequencies of $u(t)$ for $t_{min} \leq t \leq t_{evac}$. The two possible bounds,
$\eta_3^{-\frac{3}{2}}N^{\frac{3}{2}}$ and $\eta_3^{-\frac{n+2}{n-2}} N^{\frac{n+2}{n-2}}$, are a consequence of the
scaling that $g$ and $b$ obey. Note that for $n\geq 10$, $\eta_3^{-\frac{n+2}{n-2}} N^{\frac{n+2}{n-2}}$ is the larger
term.
\end{remark}

\begin{proof}
Consider the set
\begin{align*}
\Omega=\{t \in [t_{min}, t_{evac}) : \|P_{\leq N}u\|_{\dot{S}^1([t, t_{evac}]\times \R^n)}
 \leq C_0 \eta_4 + \eta_0 \max\{\eta_3^{-\frac{3}{2}}N^{\frac{3}{2}},\ &\eta_3^{-\frac{n+2}{n-2}} N^{\frac{n+2}{n-2}}\}, \\
  &\forall N \leq \eta_3\},
\end{align*}
where $C_0$ is a large constant to be chosen later and not depending on any of the $\eta$'s.

Our goal is to show that $t_{min} \in \Omega$. First, we will show that $t \in \Omega$ for $t$ close to $t_{evac}$.
Indeed, from Strichartz we get
\begin{align*}
\|P_{\leq N}u\|_{\dot{S}^1([t, t_{evac}]\times \R^n)}
    \lesssim & \|\nabla P_{\leq N}u\|_{L_t^{\infty}L_x^2([t, t_{evac}]\times \R^n)}
               + \|\nabla P_{\leq N}u\|_{L^2_tL^{\frac{2n}{n-2}}_x([t, t_{evac}]\times \R^n)} \\
    \lesssim & \|\nabla P_{\leq N}u(t_{evac})\|_{L_x^2} + C|t_{evac}-t| \|\nabla \partial_{t}P_{\leq N}u\|_{L_t^{\infty}L_x^2(I_0\times \R^n)} \\
             & + |t_{evac}-t|^{\frac{1}{2}} \|\nabla P_{\leq N}u\|_{L^{\infty}_tL^{\frac{2n}{n-2}}_x(I_0\times \R^n)}.
\end{align*}
The last two norms are finite and proportional to $N$ (as can easily be seen from Bernstein and the conservation
of energy), so \eqref{goodcontrollowfreq} implies
$$
\|P_{\leq N}u\|_{\dot{S}^1([t, t_{evac}]\times \R^n)}
   \lesssim \eta_4 + N^2|t_{evac}-t| + N|t_{evac}-t|^{\frac{1}{2}}.
$$
Thus $t \in \Omega$ provided $|t_{evac}-t|$ is sufficiently small and $C_0$ is chosen sufficiently large.

Now fix $t \in \Omega$; then for all $N \leq \eta_3$ we have
\begin{align}\label{ipoteza}
\|P_{\leq N}u\|_{\dot{S}^1([t, t_{evac}]\times \R^n)}
       \leq  C_0 \eta_4 + \eta_0 \max\{\eta_3^{-\frac{3}{2}}N^{\frac{3}{2}},\eta_3^{-\frac{n+2}{n-2}} N^{\frac{n+2}{n-2}}\}.
\end{align}
We will show that in fact,
\begin{align}
\|P_{\leq N}u\|_{\dot{S}^1([t, t_{evac}]\times \R^n)}
       \leq \frac{1}{2} C_0 \eta_4 + \frac{1}{2} \eta_0 \max\{\eta_3^{-\frac{3}{2}}N^{\frac{3}{2}},\eta_3^{-\frac{n+2}{n-2}} N^{\frac{n+2}{n-2}}\} \label{bootstraplast}
\end{align}
holds for all $ N \leq \eta_3$. Thus, $\Omega$ is both open and closed in $[t_{min},t_{evac}]$ and we get
$t_{min} \in \Omega$, as desired.

Throughout the rest of the proof all spacetime norms will be on the slab $[t, t_{evac}]\times \R^n$. Fix $N\leq \eta_3$;
by Strichartz,
\begin{align*}
\|P_{\leq N}u\|_{\dot{S}^1}
\lesssim \|P_{\leq N}u(t_{evac})\|_{\dot{H}^1_x}+\|\nabla P_{\leq N}\bigl(|u|^{\frac{4}{n-2}}u\bigr)\|_{2,\frac{2n}{n+2}}.
\end{align*}

By \eqref{goodcontrollowfreq}, we have
$$
\|P_{\leq N}u(t_{evac})\|_{\dot{H}^1_x}\lesssim \eta_4,
$$
which is acceptable if $C_0$ is chosen sufficiently large.

To handle the nonlinearity, we use the triangle inequality to estimate
\begin{align}\label{dec00}
\|\nabla P_{\leq N}\bigl(|u|^{\frac{4}{n-2}}u\bigr)\|_{2,\frac{2n}{n+2}}
&\leq \|\nabla P_{\leq N}\bigl(|u_{<\eta_4}|^{\frac{4}{n-2}}u_{<\eta_4}\bigr)\|_{2,\frac{2n}{n+2}} \notag\\
& \quad + \|\nabla P_{\leq N}\bigl(|u|^{\frac{4}{n-2}}u-|u_{<\eta_4}|^{\frac{4}{n-2}}u_{<\eta_4}\bigr)\|_{2,\frac{2n}{n+2}}.
\end{align}

Using our bootstrap hypothesis, i.e., \eqref{ipoteza}, and dropping the projection $P_{\leq N}$, we estimate the first
term on the right-hand side of \eqref{dec0} as follows:
\begin{align}\label{1!!}
\|\nabla P_{\leq N}\bigl(|u_{< \eta_4}|^{\frac{4}{n-2}}u_{<\eta_4} \bigr)\|_{2,\frac{2n}{n+2}}
&\lesssim \|\nabla u_{<\eta_4}\|_{2,\frac{2n}{n-2}}\|u_{<\eta_4}\|_{\infty,\frac{2n}{n-2}}^{\frac{4}{n-2}}\notag\\
&\lesssim \|u_{<\eta_4}\|_{\dot{S}^1}^{\frac{n+2}{n-2}}\notag\\
&\lesssim \bigl(C_0\eta_4+\eta_0\max\{\eta_3^{-\frac{3}{2}}\eta_4^{\frac{3}{2}},\eta_3^{-\frac{n+2}{n-2}} \eta_4^{\frac{n+2}{n-2}}\}\bigr)^{\frac{n+2}{n-2}}\notag\\
&\lesssim \eta_4,
\end{align}
which again is acceptable provided $C_0$ is sufficiently large.

We turn now to the second term on the right-hand side of \eqref{dec00}. By Bernstein,
\begin{align*}
\|\nabla P_{\leq N}\bigl(|u|^{\frac{4}{n-2}}u-|u_{<\eta_4}|^{\frac{4}{n-2}}&u_{<\eta_4}\bigr)\|_{2,\frac{2n}{n+2}}\\
&\lesssim N \|P_{\leq N}\bigl(|u|^{\frac{4}{n-2}}u-|u_{<\eta_4}|^{\frac{4}{n-2}}u_{<\eta_4}\bigr)\|_{2,\frac{2n}{n+2}}.
\end{align*}
Replacing the projection $P_{\leq N}$ by the positive-kernel operator $P_{\leq N}'$ (for the definition and properties of $P_{\leq N}'$
see Subsection~1.2) and using \eqref{diff1}, we further estimate
\begin{align}\label{dec0}
\|\nabla P_{\leq N}\bigl(|u|^{\frac{4}{n-2}}u-|u_{<\eta_4}|^{\frac{4}{n-2}}u_{<\eta_4}\bigr)\|_{2,\frac{2n}{n+2}}
&\lesssim N\| P_{\leq N}'\bigl(|u_{<\eta_4}|^{\frac{4}{n-2}}|u_{\geq \eta_4}|\bigr)\|_{2,\frac{2n}{n+2}}\notag\\
&\quad + N\| P_{\leq N}'\bigl(|u_{\geq \eta_4}|^{\frac{n+2}{n-2}}\bigr)\|_{2,\frac{2n}{n+2}}.
\end{align}
Decomposing $u_{\geq \eta_4} = u_{\eta_4\leq\cdot\leq \eta_3} + u_{>\eta_3}$, we estimate
\begin{align}\label{dec1}
N\| P_{\leq N}'\bigl(|u_{\geq \eta_4}|^{\frac{n+2}{n-2}}\bigr)\|_{2,\frac{2n}{n+2}}
&\lesssim N\| P_{\leq N}'\bigl(|u_{\eta_4\leq\cdot\leq\eta_3}|^{\frac{n+2}{n-2}}\bigr)\|_{2,\frac{2n}{n+2}}\notag\\
&\quad+N\| P_{\leq N}'\bigl(|u_{> \eta_3}|^{\frac{n+2}{n-2}}\bigr)\|_{2,\frac{2n}{n+2}}
\end{align}
and
\begin{align}\label{dec2}
N\| P_{\leq N}'\bigl(|u_{<\eta_4}|^{\frac{4}{n-2}}|u_{\geq \eta_4}|\bigr)\|_{2,\frac{2n}{n+2}}
&\lesssim N\| P_{\leq N}'\bigl(|u_{<\eta_4}|^{\frac{4}{n-2}}|u_{\eta_4\leq\cdot\leq\eta_3}|\bigr)\|_{2,\frac{2n}{n+2}}\notag\\
&\quad+ N\| P_{\leq N}'\bigl(|u_{<\eta_4}|^{\frac{4}{n-2}}|u_{> \eta_3}|\bigr)\|_{2,\frac{2n}{n+2}}.
\end{align}
Consider the first term on the right-hand side of \eqref{dec1}. By Bernstein,
\begin{align*}
N\|P_{\leq N}'\bigl(|u_{\eta_4\leq \cdot\leq\eta_3}|^{\frac{n+2}{n-2}}\bigr)\|_{2,\frac{2n}{n+2}}
&\lesssim N^{\frac{n+2}{n-2}}\| P_{\leq N}'\bigl(|u_{\eta_4\leq \cdot\leq\eta_3}|^{\frac{n+2}{n-2}}\bigr)\|_{2,\frac{2n(n-2)}{n^2+4}}\\
&\lesssim N^{\frac{n+2}{n-2}}\|u_{\eta_4\leq \cdot\leq\eta_3}\|_{\frac{2(n+2)}{n-2},\frac{2n(n+2)}{n^2+4}}^{\frac{n+2}{n-2}}.
\end{align*}
Using \eqref{ipoteza}, we get
\begin{align}
\|u_{\eta_4\leq \cdot\leq\eta_3}\|_{\frac{2(n+2)}{n-2},\frac{2n(n+2)}{n^2+4}}
&\lesssim \sum_{\eta_4\leq M\leq\eta_3}\|P_M u\|_{\frac{2(n+2)}{n-2},\frac{2n(n+2)}{n^2+4}}\notag\\
&\lesssim \sum_{\eta_4\leq M\leq\eta_3}M^{-1} \|\nabla P_M u\|_{\frac{2(n+2)}{n-2},\frac{2n(n+2)}{n^2+4}}\notag\\
&\lesssim \sum_{\eta_4\leq M\leq\eta_3}M^{-1} \| P_M u\|_{\dot{S}^1}\notag\\
&\lesssim \sum_{\eta_4\leq M\leq\eta_3}M^{-1}(C_0\eta_4+\eta_0\max\{\eta_3^{-\frac{3}{2}}M^{\frac{3}{2}},\eta_3^{-\frac{n+2}{n-2}} M^{\frac{n+2}{n-2}}\}\bigr)\notag\\
&\lesssim \eta_0\eta_3^{-1} \label{medfreq!}
\end{align}
and hence,
\begin{align}\label{dec1-1}
N\| P_{\leq N}'\bigl(|u_{\eta_4\leq \cdot\leq\eta_3}|^{\frac{n+2}{n-2}}\bigr)\|_{2,\frac{2n}{n+2}}
&\lesssim \eta_0^{\frac{n+2}{n-2}}\eta_3^{-\frac{n+2}{n-2}}N^{\frac{n+2}{n-2}}\notag\\
&\lesssim \eta_0^{\frac{n+2}{n-2}}\max\{\eta_3^{-\frac{3}{2}}N^{\frac{3}{2}},\eta_3^{-\frac{n+2}{n-2}}N^{\frac{n+2}{n-2}}\}.
\end{align}
To estimate the second term on the right-hand side of \eqref{dec1}, we further decompose\footnote{Of course, the
decomposition holds in dimensions $n\geq 6$. To cover the case $n=5$, we make use of Remark~\ref{rem n=5} and treat
$u_{>\eta_3}$ in the same manner as $b$ is treated below.}
$u_{>\eta_3}=g+b$ according to Corollary~\ref{cor uhi}. Using again the positivity of the operator $P_{\leq N}'$,
we estimate
\begin{align*}
N\|P_{\leq N}'\bigl(|u_{>\eta_3}|^{\frac{n+2}{n-2}}\bigr)\|_{2,\frac{2n}{n+2}}
&\lesssim N\|P_{\leq N}'\bigl(|g|^{\frac{n+2}{n-2}}\bigr)\|_{2,\frac{2n}{n+2}}
+ N\|P_{\leq N}'\bigl(|b|^{\frac{n+2}{n-2}}\bigr)\|_{2,\frac{2n}{n+2}}.
\end{align*}
By Bernstein and Corollary \ref{cor uhi}, we get
\begin{align*}
N\|P_{\leq N}'\bigl(|g|^{\frac{n+2}{n-2}}\bigr)\|_{2,\frac{2n}{n+2}}
&\lesssim N^{\frac{n+2}{n-2}}\| P_{\leq N}'\bigl(|g|^{\frac{n+2}{n-2}}\bigr)\|_{2,\frac{2n(n-2)}{n^2+4}}
\lesssim N^{\frac{n+2}{n-2}}\|g\|_{\frac{2(n+2)}{n-2},\frac{2n(n+2)}{n^2+4}}^{\frac{n+2}{n-2}}\\
&\lesssim N^{\frac{n+2}{n-2}} \bigl(\eta_{2}^{\frac{2}{n-2}}\eta_3^{-1}\bigr)^{\frac{n+2}{n-2}}
\end{align*}
and
\begin{align*}
N\|P_{\leq N}'\bigl(|b|^{\frac{n+2}{n-2}}\bigr)\|_{2,\frac{2n}{n+2}}
&\lesssim N^{\frac{3}{2}}\| P_{\leq N}'\bigl(|b|^{\frac{n+2}{n-2}}\bigr)\|_{2,\frac{2n}{n+3}}
\lesssim N^{\frac{3}{2}} \|b\|_{\frac{2(n+2)}{n-2},\frac{2n(n+2)}{(n-2)(n+3)}}^{\frac{n+2}{n-2}}\\
&\lesssim N^{\frac{3}{2}}\eta_1^{\frac{1}{4}}\eta_3^{-\frac{3}{2}}.
\end{align*}
Thus,
\begin{align}\label{dec1-2}
N\|P_{\leq N}'\bigl(|u_{>\eta_3}|^{\frac{n+2}{n-2}}\bigr)\|_{2,\frac{2n}{n+2}}
&\lesssim \eta_1^{\frac{1}{4}}\max\{\eta_3^{-\frac{3}{2}}N^{\frac{3}{2}},\eta_3^{-\frac{n+2}{n-2}}N^{\frac{n+2}{n-2}}\}\notag\\
&\lesssim \eta_0^{\frac{n+2}{n-2}}\max\{\eta_3^{-\frac{3}{2}}N^{\frac{3}{2}},\eta_3^{-\frac{n+2}{n-2}}N^{\frac{n+2}{n-2}}\}.
\end{align}
By \eqref{dec1-1} and \eqref{dec1-2}, we get
\begin{align}\label{dec0-1}
N\| P_{\leq N}'\bigl(|u_{\geq \eta_4}|^{\frac{n+2}{n-2}}\bigr)\|_{2,\frac{2n}{n+2}}
\lesssim \eta_0^{\frac{n+2}{n-2}}\max\{\eta_3^{-\frac{3}{2}}N^{\frac{3}{2}},\eta_3^{-\frac{n+2}{n-2}}N^{\frac{n+2}{n-2}}\}.
\end{align}

We turn now to the first term on the right-hand side of \eqref{dec2}. By \eqref{ipoteza} and \eqref{medfreq!},
\begin{align*}
N\| P_{\leq N}'\bigl(|u_{<\eta_4}|^{\frac{4}{n-2}}&|u_{\eta_4\leq\cdot\leq\eta_3}|\bigr)\|_{2,\frac{2n}{n+2}}\\
&\lesssim N\|u_{\eta_4\leq\cdot\leq\eta_3}\|_{\frac{2(n+2)}{n-2},\frac{2n(n+2)}{n^2+4}}\|u_{<\eta_4}\|_{\frac{2(n+2)}{n-2},\frac{2(n+2)}{n-2}}^{\frac{4}{n-2}}\\
&\lesssim \eta_0 \eta_3^{-1} N  \|u_{<\eta_4}\|_{\dot{S}^1}^{\frac{4}{n-2}}\\
&\lesssim \eta_0 \eta_3^{-1} N \bigl(C_0\eta_4+\eta_0\max\{\eta_3^{-\frac{3}{2}}\eta_4^{\frac{3}{2}},\eta_3^{-\frac{n+2}{n-2}} \eta_4^{\frac{n+2}{n-2}}\}\bigr)^{\frac{4}{n-2}}\\
&\lesssim \eta_0 \eta_4^{\frac{4}{n-2}} \eta_3^{-1} N.
\end{align*}
Treating the cases $N<\eta_4$ and $\eta_4\leq N\leq \eta_3$ separately, one easily sees that
\begin{align}\label{ineq}
\eta_0 \eta_4^{\frac{4}{n-2}} \eta_3^{-1} N
\lesssim \eta_0 \bigl(\eta_4+\eta_0\max\{\eta_3^{-\frac{3}{2}}N^{\frac{3}{2}},\eta_3^{-\frac{n+2}{n-2}} N^{\frac{n+2}{n-2}}\}\bigr).
\end{align}
Hence,
\begin{align}\label{dec2-1}
N\| P_{\leq N}'\bigl(|u_{<\eta_4}|^{\frac{4}{n-2}}&|u_{\eta_4\leq\cdot\leq\eta_3}|\bigr)\|_{2,\frac{2n}{n+2}} \notag\\
&\lesssim \eta_0 \bigl(\eta_4+\eta_0\max\{\eta_3^{-\frac{3}{2}}N^{\frac{3}{2}},\eta_3^{-\frac{n+2}{n-2}} N^{\frac{n+2}{n-2}}\}\bigr).
\end{align}
To estimate the second term on the right-hand side of \eqref{dec2}, in dimensions $n\geq 6$ we decompose
$u_{>\eta_3}=g+b$ according to Corollary~\ref{cor uhi} and use the triangle inequality and the positivity of
$P_{\geq N}'$ to bound
\begin{align*}
N\| P_{\leq N}'\bigl(|u_{<\eta_4}|^{\frac{4}{n-2}}|u_{> \eta_3}|\bigr)\|_{2,\frac{2n}{n+2}}
&\lesssim N\| P_{\leq N}'\bigl(|u_{<\eta_4}|^{\frac{4}{n-2}}|g|\bigr)\|_{2,\frac{2n}{n+2}}\\
&\quad +N\| P_{\leq N}'\bigl(|u_{<\eta_4}|^{\frac{4}{n-2}}|b|\bigr)\|_{2,\frac{2n}{n+2}}.
\end{align*}
By Bernstein, Corollary \ref{cor uhi}, and \eqref{ipoteza}, we get
\begin{align*}
N\| P_{\leq N}'\bigl(|u_{<\eta_4}|^{\frac{4}{n-2}}|g|\bigr)\|_{2,\frac{2n}{n+2}}
&\lesssim N\|g\|_{\frac{2(n+2)}{n-2},\frac{2n(n+2)}{n^2+4}}\|u_{<\eta_4}\|_{\frac{2(n+2)}{n-2},\frac{2(n+2)}{n-2}}^{\frac{4}{n-2}}\\
&\lesssim N\|g\|_{\dot{S}^0}\|u_{<\eta_4}\|_{\dot{S}^1}^{\frac{4}{n-2}}
\lesssim N \eta_2^{\frac{2}{n-2}}\eta_3^{-1}\eta_4^{\frac{4}{n-2}}
\end{align*}
and
\begin{align*}
N\| P_{\leq N}'\bigl(|u_{<\eta_4}|^{\frac{4}{n-2}}|b|\bigr)\|_{2,\frac{2n}{n+2}}
&\lesssim N N^{\frac{n-6}{2(n+2)}}\| P_{\leq N}'\bigl(|u_{<\eta_4}|^{\frac{4}{n-2}}|b|\bigr)\|_{2,\frac{2n(n+2)}{n^2+5n-2}}\\
&\lesssim N^{\frac{3n-2}{2(n+2)}}\|u_{<\eta_4}\|_{\frac{2(n+2)}{n-2},\frac{2(n+2)}{n-2}}^{\frac{4}{n-2}}\|b\|_{\frac{2(n+2)}{n-2},\frac{2n}{n-1}}\\
&\lesssim \eta_4^{\frac{4}{n-2}}N^{\frac{3n-2}{2(n+2)}}\|b\|_{\frac{2n}{n-2},\frac{2n^2}{(n+1)(n-2)}}^{\frac{n}{n+2}}\|b\|_{\infty,2}^{\frac{2}{n+2}}\\
&\lesssim \eta_4^{\frac{4}{n-2}}N^{\frac{3n-2}{2(n+2)}}(\eta_1^{\frac{1}{4}}\eta_3^{-\frac{3}{2}})^{\frac{n-2}{n+2}}(\eta_2^{\frac{2}{n-2}}\eta_3^{-1})^{\frac{2}{n+2}}\\
&\lesssim \eta_1 \eta_4^{\frac{4}{n-2}}N^{\frac{3n-2}{2(n+2)}}\eta_3^{-\frac{3n-2}{2(n+2)}}.
\end{align*}
Thus, by treating the cases $N<\eta_4$ and $\eta_4\leq N\leq \eta_3$ separately, one sees that in dimensions $n\geq 6$,
\begin{align*}
N\| P_{\leq N}'\bigl(|u_{<\eta_4}|^{\frac{4}{n-2}}|u_{> \eta_3}|\bigr)\|_{2,\frac{2n}{n+2}}
\lesssim \eta_1 \bigl(\eta_4+\eta_0\max\{\eta_3^{-\frac{3}{2}}\eta_4^{\frac{3}{2}},\eta_3^{-\frac{n+2}{n-2}} \eta_4^{\frac{n+2}{n-2}}\}\bigr).
\end{align*}
In dimension $n=5$, by Remark~\ref{rem n=5} and \eqref{ipoteza} we have
\begin{align*}
N\| P_{\leq N}'\bigl(|u_{<\eta_4}|^{\frac{4}{n-2}}|u_{> \eta_3}|\bigr)\|_{2,\frac{2n}{n+2}}
&= N\| P_{\leq N}'\bigl(|u_{<\eta_4}|^{\frac{4}{3}}|u_{> \eta_3}|\bigr)\|_{2,\frac{10}{7}}\\
&\lesssim N \|u_{>\eta_3}\|_{3,\frac{30}{11}}\|u_{<\eta_4}\|^{\frac{4}{3}}_{8,4}\\
&\lesssim N \eta_1^{\frac{1}{3}}\eta_3^{-1} \|u_{<\eta_4}\|_{\dot{S}^1}^{\frac{4}{3}}\\
&\lesssim \eta_1^{\frac{1}{3}} \eta_4^{\frac{4}{3}} N\eta_3^{-1}\lesssim\eta_1\eta_4.
\end{align*}
Hence, for all $n\geq 5$,
\begin{align}\label{dec2-2}
N\| P_{\leq N}'\bigl(|u_{<\eta_4}|^{\frac{4}{n-2}}|u_{> \eta_3}|\bigr)\|_{2,\frac{2n}{n+2}}
\lesssim \eta_1 \bigl(\eta_4+\eta_0\max\{\eta_3^{-\frac{3}{2}}\eta_4^{\frac{3}{2}},\eta_3^{-\frac{n+2}{n-2}} \eta_4^{\frac{n+2}{n-2}}\}\bigr).
\end{align}
By \eqref{dec2-1} and \eqref{dec2-2}, we obtain
\begin{align}\label{dec0-2}
N\| P_{\leq N}'\bigl(|u_{<\eta_4}|^{\frac{4}{n-2}}|u_{\geq \eta_4}|\bigr)\|_{2,\frac{2n}{n+2}}
\lesssim \eta_0 \bigl(\eta_4+\eta_0\max\{\eta_3^{-\frac{3}{2}}\eta_4^{\frac{3}{2}},\eta_3^{-\frac{n+2}{n-2}} \eta_4^{\frac{n+2}{n-2}}\}\bigr).
\end{align}
By \eqref{dec0}, \eqref{dec0-1}, and \eqref{dec0-2}, we get
\begin{align}\label{2!!}
\|\nabla P_{\leq N}\bigl(|u|^{\frac{4}{n-2}}u-|u_{<\eta_4}|^{\frac{4}{n-2}}u_{<\eta_4}\bigr)\|_{2,\frac{2n}{n+2}}
\lesssim \eta_4 + \eta_0^{1+}\max\{\eta_3^{-\frac{3}{2}}N^{\frac{3}{2}},\eta_3^{-\frac{n+2}{n-2}}N^{\frac{n+2}{n-2}}\}.
\end{align}
By \eqref{dec00}, \eqref{1!!}, and \eqref{2!!},
\begin{align*}
\|\nabla P_{\leq N}\bigl(|u|^{\frac{4}{n-2}}u\bigr)\|_{2,\frac{2n}{n+2}}
\lesssim \eta_4 + \eta_0^{1+}\max\{\eta_3^{-\frac{3}{2}}N^{\frac{3}{2}},\eta_3^{-\frac{n+2}{n-2}}N^{\frac{n+2}{n-2}}\}.
\end{align*}
Hence, \eqref{bootstraplast} holds for $C_1$ sufficiently large and the proof of Lemma \ref{S1lf} is complete.

\end{proof}

\subsection{Controlling the localized mass increment.} We now have good enough control over low, medium, and high
frequencies to prove \eqref{last}.  Writing
\begin{align*}
P_{hi}(|u|^{\frac{4}{n-2}}u)-|u_{hi}|^{\frac{4}{n-2}}u_{hi}
  &=P_{hi}(|u|^{\frac{4}{n-2}}u-|u_{hi}|^{\frac{4}{n-2}}u_{hi}-|u_{lo}|^{\frac{4}{n-2}}u_{lo})\\
  &\quad - P_{lo}(|u_{hi}|^{\frac{4}{n-2}}u_{hi})+ P_{hi}(|u_{lo}|^{\frac{4}{n-2}}u_{lo}),
\end{align*}
we see that we only have to consider the following terms
\begin{align}
\int_{t_{min}}^{t_{*}}\Bigl|\int_{\R^n} \overline{u_{hi}}P_{hi}(|u|^{\frac{4}{n-2}}u-|u_{hi}|^{\frac{4}{n-2}}u_{hi}-|u_{lo}|^{\frac{4}{n-2}}u_{lo})dx\Bigr|dt \label{1}, \\
\int_{t_{min}}^{t_{*}}\Bigl|\int_{\R^n} \overline{u_{hi}}P_{lo}(|u_{hi}|^{\frac{4}{n-2}}u_{hi})dx\Bigr|dt \label{2}, \\
\int_{t_{min}}^{t_{*}}\Bigl|\int_{\R^n} \overline{u_{hi}}P_{hi}(|u_{lo}|^{\frac{4}{n-2}}u_{lo})dx\Bigr|dt \label{3}.
\end{align}

For the remaining of this section all spacetime norms will be on the slab $[t_{min}, t_*]\times\R^n$. Consider \eqref{1}.
By \eqref{diff2}, we estimate
\begin{align*}
\eqref{1}
&\lesssim \||u_{hi}|^2|u_{lo}|^{\frac{4}{n-2}}\chi_{\{|u_{hi}|\ll|u_{lo}|\}}\|_{1,1}+\||u_{hi}|^{\frac{n+2}{n-2}}u_{lo}\chi_{\{|u_{lo}|\ll|u_{hi}|\}}\|_{1,1}\\
&\lesssim \||u_{hi}|^2|u_{lo}|^{\frac{4}{n-2}}\|_{1,1}
 \lesssim \||u_{>\eta_3}|^2|u_{lo}|^{\frac{4}{n-2}}\|_{1,1}+\||u_{\eta_3^{10n}\leq\cdot\leq \eta_3}|^2|u_{lo}|^{\frac{4}{n-2}}\|_{1,1}.
\end{align*}
Taking $N_*=\eta_3$ in Corollary~\ref{cor uhi}, we decompose $u_{>\eta_3}=g+b$ and estimate
$$
\||u_{>\eta_3}|^2|u_{lo}|^{\frac{4}{n-2}}\|_{1,1}
\lesssim \||g|^2|u_{lo}|^{\frac{4}{n-2}}\|_{1,1}+\||b|^2|u_{lo}|^{\frac{4}{n-2}}\|_{1,1}.
$$
Using H\"older, Bernstein, Corollary~\ref{cor uhi}, and Proposition \ref{S1lf}, we get
\begin{align*}
\||g|^2|u_{lo}|^{\frac{4}{n-2}}\|_{1,1}
&\lesssim \|g\|_{2,\frac{2n}{n-2}}^2 \|u_{lo}\|_{\infty,\frac{2n}{n-2}}^{\frac{4}{n-2}}\\
&\lesssim \bigl(\eta_2^{\frac{2}{n-2}}\eta_3^{-1}\bigr)^2\bigl(\eta_4+\max\{\eta_3^{-\frac{3}{2}}\eta_3^{10n\frac{3}{2}},\eta_3^{-\frac{n+2}{n-2}}\eta_3^{10n\frac{n+2}{n-2}}\}\bigr)^{\frac{4}{n-2}}\\
&\ll\eta_1^2,\\
\||b|^2|u_{lo}|^{\frac{4}{n-2}}\|_{1,1}
&\lesssim \|b\|_{\frac{2n}{n-2},\frac{2n^2}{(n+1)(n-2)}}^2\|u_{lo}\|_{\frac{2n}{n-2},\frac{4n^2}{(n+2)(n-2)}}^{\frac{4}{n-2}}\\
&\lesssim \bigl(\eta_1^{\frac{1}{4}}\eta_3^{-\frac{3}{2}}\bigr)^{\frac{2(n-2)}{n}}\eta_3^{10(n-6)}\|u_{lo}\|_{\dot{S}^1}^{\frac{4}{n-2}}\\
&\lesssim \bigl(\eta_1^{\frac{1}{2}}\eta_3^{-3}\bigr)^{\frac{n-2}{n}}\eta_3^{10(n-6)}\bigl(\eta_4+\max\{\eta_3^{-\frac{3}{2}}\eta_3^{10n\frac{3}{2}},\eta_3^{-\frac{n+2}{n-2}}\eta_3^{10n\frac{n+2}{n-2}}\}\bigr)^{\frac{4}{n-2}}\\
&\ll \eta_1^2,
\end{align*}
where in the last sequence of inequalities we used the fact that for $n\geq 6$, Bernstein dictates
$$
\|u_{lo}\|_{\frac{2n}{n-2},\frac{4n^2}{(n+2)(n-2)}}^{\frac{4}{n-2}}
\lesssim \eta_3^{10(n-6)}\|u_{lo}\|_{\frac{2n}{n-2},\frac{2n^2}{(n-2)^2}}^{\frac{4}{n-2}}
\lesssim \eta_3^{10(n-6)}\|u_{lo}\|_{\dot{S}^1}.
$$
To cover the case $n=5$, we estimate
$$
\||u_{>\eta_3}|^2 |u_{lo}|^{\frac{4}{n-2}}\|_{1,1}=\||u_{>\eta_3}|^2 |u_{lo}|^{\frac{4}{3}}\|_{1,1}
\lesssim\|u_{>\eta_3}\|_{3,\frac{30}{11}}^2\|u_{lo}\|_{4,5}^{\frac{4}{3}}.
$$
In dimension $n=5$, one easily checks that the pair $(4,\tfrac{5}{2})$
is Schr\"odinger admissible and hence, by Sobolev embedding,
\begin{align}\label{ulo12/n-2}
\|u_{lo}\|_{4,5}
\lesssim \|\nabla u_{lo}\|_{4,\frac{5}{2}}
\lesssim \|u_{lo}\|_{\dot{S}^1}.
\end{align}
Thus, by \eqref{ulo12/n-2}, Remark \ref{rem n=5} (with $N_*=\eta_3$), and Proposition~\ref{S1lf}
(with $n=5$), we get
\begin{align*}
\||u_{>\eta_3}|^2 |u_{lo}|^{\frac{4}{3}}\|_{1,1}
\lesssim \eta_1^{\frac{2}{3}}\eta_3^{-2}\bigl(\eta_4+\eta_3^{-\frac{3}{2}}\eta_3^{75}\bigr)^{\frac{4}{3}}
\ll \eta_1^2.
\end{align*}
Hence, in all dimensions $n\geq 5$ we have
$$
\||u_{>\eta_3}|^2|u_{lo}|^{\frac{4}{n-2}}\|_{1,1}\ll \eta_1^2.
$$
Next, by H\"older,
$$
\||u_{\eta_3^{10n}\leq\cdot\leq \eta_3}|^2|u_{lo}|^{\frac{4}{n-2}}\|_{1,1}
\lesssim \|u_{\eta_3^{10n}\leq\cdot\leq \eta_3}\|_{2,\frac{2n}{n-2}}^2\|u_{lo}\|_{\infty,\frac{2n}{n-2}}^{\frac{4}{n-2}}.
$$
Using Proposition \ref{S1lf}, we estimate
\begin{align}
\|u_{\eta_3^{10n}\leq\cdot\leq \eta_3}\|_{2,\frac{2n}{n-2}}
&\lesssim \sum_{\eta_3^{10n}\leq N\leq
\eta_3}\|u_N\|_{2,\frac{2n}{n-2}}
 \lesssim \sum_{\eta_3^{10n}\leq N\leq \eta_3}N^{-1}\|u_N\|_{\dot{S}^1}\notag\\
&\lesssim \sum_{\eta_3^{10n}\leq N\leq \eta_3}N^{-1}\bigl(\eta_4+\max\{\eta_3^{-\frac{3}{2}}N^{\frac{3}{2}},\eta_3^{-\frac{n+2}{n-2}}N^{\frac{n+2}{n-2}}\}\bigr)
 \lesssim \eta_3^{-1}\label{umed!}
\end{align}
and hence,
$$
\||u_{\eta_3^{10n}\leq\cdot\leq \eta_3}|^2|u_{lo}|^{\frac{4}{n-2}}\|_{1,1}
\lesssim \eta_3^{-2}\bigl(\eta_4+\max\{\eta_3^{-\frac{3}{2}}\eta_3^{10n\frac{3}{2}},\eta_3^{-\frac{n+2}{n-2}}\eta_3^{10n\frac{n+2}{n-2}}\}\bigr)^{\frac{4}{n-2}}
\ll \eta_1^2.
$$
Therefore,
$$
\eqref{1}\ll \eta_1^2.
$$

To estimate \eqref{2}, we write
$$
\eqref{2}=\int_{t_{min}}^{t_*}\Bigl|\int \overline{P_{lo}(u_{hi})}|u_{hi}|^{\frac{4}{n-2}}u_{hi}dx\Bigr|dt.
$$
As $P_{lo}(u_{hi})=P_{hi}(u_{lo})$ satisfies all the estimates that $u_{lo}$ and $u_{hi}$ satisfy, we see by the previous
analysis that
$$
\eqref{2}\lesssim \||u_{hi}|^2|u_{lo}|^{\frac{4}{n-2}}\|_{1,1}\ll\eta_1^2.
$$

We consider next \eqref{3}. We estimate
\begin{align*}
\eqref{3}
&\lesssim \|u_{hi}P_{hi}\bigl(|u_{lo}|^{\frac{4}{n-2}}u_{lo}\bigr)\|_{1,1}\\
&\lesssim \|u_{>\eta_3}P_{hi}\bigl(|u_{lo}|^{\frac{4}{n-2}}u_{lo}\bigr)\|_{1,1}
  +\|u_{\eta_3^{10n}\leq\cdot\leq \eta_3}P_{hi}\bigl(|u_{lo}|^{\frac{4}{n-2}}u_{lo}\bigr)\|_{1,1}.
\end{align*}
Using Bernstein to place a derivative in front of $P_{hi}$, we get
$$
\|u_{\eta_3^{10n}\leq\cdot\leq \eta_3}P_{hi}\bigl(|u_{lo}|^{\frac{4}{n-2}}u_{lo}\bigr)\|_{1,1}
\lesssim \eta_3^{-10n}\|u_{\eta_3^{10n}\leq\cdot\leq \eta_3}\|_{2,\frac{2n}{n-2}}\|\nabla u_{lo}\|_{2,\frac{2n}{n-2}}\|u_{lo}\|_{\infty,\frac{2n}{n-2}}^{\frac{4}{n-2}}.
$$
By \eqref{umed!} and Proposition \ref{S1lf}, we obtain
\begin{align*}
\|u_{\eta_3^{10n}\leq\cdot\leq \eta_3}P_{hi}\bigl(&|u_{lo}|^{\frac{4}{n-2}}u_{lo}\bigr)\|_{1,1}\\
&\lesssim \eta_3^{-10n}\eta_3^{-1}\bigl(\eta_4+\max\{\eta_3^{-\frac{3}{2}}\eta_3^{10n\frac{3}{2}},\eta_3^{-\frac{n+2}{n-2}}\eta_3^{10n\frac{n+2}{n-2}}\}\bigr)^{\frac{n+2}{n-2}}\\
&\ll \eta_1^2.
\end{align*}
On the other hand, decomposing $u_{>\eta_3}=g+b$ in dimensions $n\geq 6$ according to Corollary~\ref{cor uhi}, we estimate
$$
\|u_{>\eta_3}P_{hi}\bigl(|u_{lo}|^{\frac{4}{n-2}}u_{lo}\bigr)\|_{1,1}
\lesssim \|gP_{hi}\bigl(|u_{lo}|^{\frac{4}{n-2}}u_{lo}\bigr)\|_{1,1}
  +\|bP_{hi}\bigl(|u_{lo}|^{\frac{4}{n-2}}u_{lo}\bigr)\|_{1,1}.
$$
By Bernstein, \eqref{uloinS^1}, Corollary \ref{cor uhi}, and Proposition \ref{S1lf}, we have
\begin{align*}
\|gP_{hi}\bigl(|u_{lo}|^{\frac{4}{n-2}}u_{lo}\bigr)\|_{1,1}
&\lesssim \eta_3^{-10n}\|g\|_{2,\frac{2n}{n-2}}\|\nabla u_{lo}\|_{2,\frac{2n}{n-2}}\|u_{lo}\|_{\infty,\frac{2n}{n-2}}^{\frac{4}{n-2}}\\
&\lesssim \eta_3^{-10n}\|g\|_{\dot{S}^0}\| u_{lo}\|_{\dot{S}^1}^{\frac{n+2}{n-2}}\\
&\lesssim \eta_3^{-10n}\eta_2^{\frac{2}{n-2}}\eta_3^{-1}\bigl(\eta_4+\max\{\eta_3^{-\frac{3}{2}}\eta_3^{10n\frac{3}{2}},\eta_3^{-\frac{n+2}{n-2}}\eta_3^{10n\frac{n+2}{n-2}}\}\bigr)^{\frac{n+2}{n-2}}\\
&\ll \eta_1^2
\end{align*}
and
\begin{align*}
\|bP_{hi}\bigl(|u_{lo}&|^{\frac{4}{n-2}}u_{lo}\bigr)\|_{1,1}\\
&\lesssim \eta_3^{-10n}\|b\|_{\frac{2n}{n-2},\frac{2n^2}{(n+1)(n-2)}}\|\nabla u_{lo}\|_{\frac{2n}{n-2},\frac{2n^2}{n^2-2n+4}}\|u_{lo}\|_{\frac{2n}{n-2},\frac{8n^2}{(3n-2)(n-2)}}^{\frac{4}{n-2}}\\
&\lesssim \eta_3^{-5(n+6)}\bigl(\eta_1^{\frac{1}{4}}\eta_3^{-\frac{3}{2}}\bigr)^{\frac{n-2}{n}}\|u_{lo}\|_{\dot{S}^1}^{\frac{n+2}{n-2}}\\
&\lesssim \eta_3^{-5(n+6)}\bigl(\eta_1^{\frac{1}{4}}\eta_3^{-\frac{3}{2}}\bigr)^{\frac{n-2}{n}}\bigl(\eta_4+\max\{\eta_3^{-\frac{3}{2}}\eta_3^{10n\frac{3}{2}},\eta_3^{-\frac{n+2}{n-2}}\eta_3^{10n\frac{n+2}{n-2}}\}\bigr)^{\frac{n+2}{n-2}}\\
&\ll \eta_1^2.
\end{align*}
To cover the case $n=5$, we use Bernstein (to add a derivative in front of $P_{hi}$), Remark \ref{rem n=5},
\eqref{ulo12/n-2}, and Proposition~\ref{S1lf} (with $n=5$) to estimate instead
\begin{align*}
\|u_{>\eta_3} P_{hi}\bigl(|u_{lo}|^{\frac{4}{3}}u_{lo}\bigr)\|_{1,1}
&\lesssim \|u_{>\eta_3}\|_{3,\frac{15}{11}}\|\nabla u_{lo}\|_{3,\frac{15}{11}}\|u_{lo}\|_{4,5}^{\frac{4}{3}}\\
&\lesssim \eta_1^{\frac{1}{3}}\eta_3^{-1} \|u_{lo}\|_{\dot{S}^1}^{\frac{7}{3}}\\
&\lesssim \eta_1^{\frac{1}{3}}\eta_3^{-1}\bigl(\eta_4+\eta_3^{-\frac{3}{2}}\eta_3^{75}\bigr)^{\frac{7}{3}}
\ll\eta_1^2.
\end{align*}
Thus,
$$
\|u_{>\eta_3}P_{hi}\bigl(|u_{lo}|^{\frac{4}{n-2}}u_{lo}\bigr)\|_{1,1} \ll \eta_1^2
$$
for all $n\geq 5$ and hence
$$
\eqref{3}\ll \eta_1^2.
$$

Therefore \eqref{last} holds and this concludes the proof of Proposition \ref{energyevac}.

%
%
%
%

\section{The contradiction argument}

We now have all the information we need about a minimal energy blowup solution to conclude the contradiction argument.
Corollary \ref{lemma freq loc} shows that it is localized in frequency and Proposition \ref{lemma physical concentration}
that it concentrates in space. The interaction Morawetz inequality provides good control over the high-frequency
part of $u$ in $L^3_tL^{\frac{6n}{3n-4}}_x$ (see Corollary~\ref{cor uhi} and Remark~\ref{rem n=5}).
By the arguments in the previous section we have also excluded the last enemy by showing
that the solution can't shift its energy from low modes to high modes causing the $L^{\frac{2(n+2)}{n-2}}_{t,x}$-norm
to blow up while the $L^3_tL^{\frac{6n}{3n-4}}_x$-norm remains bounded. Hence, $N(t)$ must remain within a bounded
set $[N_{min},N_{max}]$, where $N_{max}\leq C(\eta_4)N_{min}$ and $N_{min}>0$.
Combining all these (and relying again on the interaction Morawetz inequality), we will derive the desired contradiction.
We begin with

\begin{lemma}\label{lemma integral bound on N}
For any minimal energy blowup solution to \eqref{schrodinger equation}, we have
\begin{equation}\label{integral bound on N}
\int_{I_0} N(t)^{-1} dt \lesssim C(\eta_1,\eta_2) N_{min}^{-3}.
\end{equation}
In particular, as $N(t) \leq C(\eta_4) N_{min}$ for all $t\in I_0$, we have
\begin{equation}\label{length I bound}
|I_0| \lesssim C(\eta_1, \eta_2, \eta_4) N_{min}^{-2}.
\end{equation}
\end{lemma}

\begin{proof}
By \eqref{uhiL^3!!!} and \eqref{n=5}, in all dimensions $n\geq 5$ we have
\begin{equation*}
\int_{I_0} \Bigl(\int_{\R^n} |P_{\geq N_*}u|^{\frac{6n}{3n-4}} dx \Bigr)^{\frac{3n-4}{2n}} dt \lesssim \eta_1^{\frac{1}{2}} N_*^{-3}
\end{equation*}
for all $N_* < c(\eta_2)N_{min}$.  Let $N_* = c(\eta_2) N_{min}$ and rewrite the above estimate as
\begin{equation}\label{int bound low}
\int_{I_0} \Bigl(\int_{\R^n} |P_{\geq N_*}u|^{\frac{6n}{3n-4}} dx \Bigr)^{\frac{3n-4}{2n}} dt  \lesssim C(\eta_1,\eta_2)N_{min}^{-3}.
\end{equation}
On the other hand, by Bernstein and the conservation of energy,
\begin{equation}\label{int bound high}
\begin{split}
\int_{|x-x(t)| \leq C(\eta_1)/N(t)}|P_{<N_*}u(t)&|^{\frac{6n}{3n-4}} dx\\
& \lesssim C(\eta_1)N(t)^{-n} \|P_{<N_*} u(t)\|_{L^\infty_x}^{\frac{6n}{3n-4}} \\
& \lesssim C(\eta_1)N(t)^{-n}N(t)^{\frac{3n(n-2)}{3n-4}} c(\eta_2)\|P_{<N_*}u(t)\|_{L^{\frac{2n}{n-2}}_x}^{\frac{6n}{3n-4}} \\
& \lesssim c(\eta_2) N(t)^{-\frac{2n}{3n-4}}.
\end{split}
\end{equation}
By \eqref{physical conc lp}, we also have
$$
\int_{|x-x(t)| \leq C(\eta_1)/N(t)}|u(t)|^{\frac{6n}{3n-4}} dx \gtrsim c(\eta_1) N(t)^{-\frac{2n}{3n-4}}.
$$
Combining this estimate with \eqref{int bound high} and using the triangle inequality, we find
\begin{equation*}
c(\eta_1) N(t)^{-\frac{2n}{3n-4}} \lesssim \int_{|x-x(t)| \leq C(\eta_1)/N(t)}|P_{\geq N_*}u(t,x)|^{\frac{6n}{3n-4}} dx.
\end{equation*}
Integrating over $I_0$ and comparing with \eqref{int bound low} proves \eqref{integral bound on N}.
\end{proof}

We can now (finally!) conclude the contradiction argument. It remains to prove that
$\|u\|_{L^{\frac{2(n+2)}{n-2}}_{t,x}(I_0 \times \R^n)} \lesssim C(\eta_0,\eta_1,\eta_2,\eta_3,\eta_4)$, which
contradicts \eqref{HUGE} for $\eta_5$ sufficiently small and which we expect since
the bound \eqref{length I bound} shows that the interval $I_0$ is not long enough to allow the
$L^{\frac{2(n+2)}{n-2}}_{t,x}$-norm of $u$ to grow too large. Indeed, we have

\begin{proposition}\label{lemma noncon imp spacetime bdd}
$$
\|u\|_{L^{\frac{2(n+2)}{n-2}}_{t,x}(I_0 \times \R^n)} \lesssim C(\eta_0,\eta_1,\eta_2,\eta_4).
$$
\end{proposition}

\begin{proof}
We normalize $N_{min} = 1$.  Let $\delta = \delta(\eta_0,N_{max}) > 0$ be a small number to be chosen momentarily.
Partition $I_0$ into $O(|I_0|/\delta)$ subintervals $I_1,\dots,I_J$ with $|I_j| \leq \delta$.  Let $t_j \in I_j$.
As $N(t_j) \leq N_{max}$, Corollary \ref{lemma freq loc} yields
\begin{align*}
\|P_{\geq C(\eta_0)N_{max}}u(t_j)\|_{\ho_x} \leq \eta_0.
\end{align*}
Let $\tilde u(t) = e^{i(t-t_j)\Delta}P_{<C(\eta_0)N_{max}}u(t_j)$ be the free evolution of the low and medium frequencies
of $u(t_j)$. The above bound becomes
\begin{align}\label{close 1}
\|u(t_j) - \tilde u(t_j)\|_{\ho_x} \leq \eta_0.
\end{align}
Moreover, by Remark~\ref{redundant}, \eqref{close 1} implies
\begin{align}\label{close 2}
\Bigl(\sum_{N}\| e^{i(t-t_j)\Delta}P_N\nabla (u(t_j) - \tilde u(t_j))\|_{L_t^{\frac{2(n+2)}{n-2}}L_x^{\frac{2n(n+2)}{n^2+4}}(I_j\times\R^n)}^2\Bigr)^{\frac{1}{2}}
\lesssim \eta_0.
\end{align}
By Bernstein, Sobolev embedding, and conservation of energy, we get
\begin{align*}
\|\tilde u(t) \|_{L^{\frac{2(n+2)}{n-2}}_{x}}
&\lesssim C(\eta_0,N_{max}) \|\tilde u(t_j)\|_{L^{\frac{2n}{n-2}}_{x}}
\lesssim C(\eta_0, N_{max})\|\tilde u(t_j)\|_{\ho_x}\\
&\lesssim C(\eta_0,N_{max})
\end{align*}
for all $t \in I_j$, so
\begin{align}\label{tilde u small}
\|\tilde u \|_{L^{\frac{2(n+2)}{n-2}}_{t,x}(I_j \times\R^n)} \lesssim C(\eta_0, N_{max}) \delta^{\frac{n-2}{2(n+2)}}.
\end{align}
Similarly, we have
\begin{align*}
\|\nabla (|\tilde u(t)|^{\frac{4}{n-2}} \tilde u(t))\|_{L^{\frac{2n}{n+2}}_x}
& \lesssim \|\nabla \tilde u(t) \|_{L^{\frac{2n}{n-2}}_x} \|\tilde u(t)\|_{L^{\frac{2n}{n-2}}_x}^{\frac{4}{n-2}}\\
&\lesssim C(\eta_0, N_{max}) \|\nabla \tilde u(t) \|_{L^2_x} \|\tilde u(t) \|_{\ho_x}^{\frac{4}{n-2}}\\
&\lesssim C(\eta_0, N_{max}) \|\tilde u(t) \|_{\ho_x}^{\frac{n+2}{n-2}}\\
&\lesssim C(\eta_0, N_{max}),
\end{align*}
which shows that
\begin{align}\label{e small}
\|\nabla(|\tilde u|^{\frac{4}{n-2}} \tilde u)\|_{L^2_t L^{\frac{2n}{n+2}}_x (I_j \times \R^n)} \lesssim C(\eta_0, N_{max}) \delta ^{1/2}.
\end{align}
By \eqref{close 1} through \eqref{e small}, conservation of energy, and
Lemma \ref{lemma long time} with $e=-|\tilde u|^{\frac{4}{n-2}}\tilde u$, we see that
$$
\|u\|_{L^{\frac{2(n+2)}{n-2}}_{t,x}(I_j \times \R^n)} \lesssim 1,
$$
provided $\delta$ and $\eta_0$ are chosen small enough.  Summing these bounds in $j$ and using \eqref{length I bound}, we get
$$
\|u\|_{L^{\frac{2(n+2)}{n-2}}_{t,x}(I_0 \times \R^n)} \lesssim C(\eta_0, N_{max}) |I_0|
 \lesssim C(\eta_0,\eta_1,\eta_2, \eta_4).
$$
\end{proof}

%
%
%
%

\appendix
\section{Fractional derivatives of fractional powers}

In this appendix, we show how a characterization of Sobolev spaces due to Strichartz, \cite{Strich}, can be used
to prove results of chain-rule type.  This extends \cite{ChW:fractional chain rule} from $C^1$ functions to merely H\"older continuous functions.

The results in this section were worked out in collaboration with Rowan Killip.

Strichartz proved that for all Schwartz functions $f$, $1<p<\infty$, and $0<s<1$,
\begin{equation}\label{E:Strich}
\bigl\||\nabla|^s f \bigr\|_{L^p_x} \approx \bigl\| \mathcal{D}_s(f) \bigr\|_{L^p_x},
\end{equation}
where
\begin{equation}\label{E:Ddef}
\mathcal{D}_s(f)(x) = \biggl( \int_0^\infty \biggl|\int_{|y|<1} \bigl|f(x+ry)-f(x)\bigr| \,dy\biggr|^2 \frac{dr}{r^{1+2s}}\biggr)^{1/2}.
\end{equation}
This extended earlier work of Stein; see the discussion in \cite[Ch. V, \S6.13]{stein:small}.

\begin{proposition}\label{fdfp}
Let $F$ be a H\"older continuous function of order $0<\alpha<1$.  Then, for every $0<\sigma<\alpha$, $1<p<\infty$,
and $\tfrac{\sigma}{\alpha}<s<1$ we have
\begin{align}\label{fdfp2}
\bigl\| |\nabla|^\sigma F(u)\bigr\|_p
\lesssim \bigl\||u|^{\alpha-\frac{\sigma}{s}}\bigr\|_{p_1} \bigl\||\nabla|^s u\bigr\|^{\frac{\sigma}{s}}_{\frac{\sigma}{s}p_2},
\end{align}
provided $\tfrac{1}{p}=\tfrac{1}{p_1} +\tfrac{1}{p_2}$ and $(1-\frac\sigma{\alpha s})p_1>1$.
\end{proposition}

\begin{proof}
The result will follow from the pointwise inequality
\begin{equation}\label{E:pointwise}
\mathcal{D}_\sigma(F(u))(x) \lesssim \bigl[ M(|u|^\alpha)(x) \bigr]^{1-\frac\sigma{\alpha s}} \bigl[ \mathcal{D}_s(u)(x) \bigr]^{\frac\sigma s},
\end{equation}
where $M$ denotes the Hardy--Littlewood maximal function.

As $F$ is $\alpha$-H\"older continuous,
$$
|F(u(x+ry))-F(u(x))| \lesssim  |u(x+ry)-u(x)|^\alpha \lesssim  |u(x+ry)|^\alpha + |u(x)|^\alpha.
$$
We use both estimates; the first one for small values of $r$, the second for large values of $r$.  The meaning
of `small' and `large' will be $x$-dependent.

For $r$ small, we apply H\"older's inequality:
\begin{align*}
\int_0^{A(x)} \biggl|\int_{|y|<1} &\bigl|F(u(x+ry))-F(u(x))\bigr| \,dy\biggr|^2 \frac{dr}{r^{1+2\sigma}} \\
&\lesssim \int_0^{A(x)} \biggl|\int_{|y|<1} \bigl|u(x+ry)-u(x)\bigr|^\alpha \,dy\biggr|^2 \frac{dr}{r^{1+2\sigma}} \\
&\lesssim \int_0^{A(x)} \biggl|\int_{|y|<1} \bigl|u(x+ry)-u(x)\bigr| \,dy\biggr|^{2\alpha} \frac{dr}{r^{1+2\sigma}} \\
&\lesssim \bigl[A(x)\bigr]^{2(s\alpha-\sigma)}\left( \int_0^{A(x)} \biggl|\int_{|y|<1} \bigl|u(x+ry)-u(x)\bigr| \,dy\biggr|^2 \frac{dr}{r^{1+2s}} \right)^\alpha \\
&\lesssim \bigl[A(x)\bigr]^{2(s\alpha-\sigma)} \bigl[ \mathcal{D}_s(u)(x) \bigr]^{2\alpha}.
\end{align*}
Note that the penultimate step requires $s\alpha-\sigma>0$.

For large $r$, we first note that
\begin{align*}
\int_{|y|<1} \bigl|u(x+ry)\bigr|^\alpha \,dy \lesssim M\bigl(|u|^\alpha \bigr)(x)
\end{align*}
because the left-hand side is essentially the average of $|u|^\alpha$ over the $r$-ball centered at $x$.
Consequently,
\begin{align*}
\int_{A(x)}^\infty \biggl|\int_{|y|<1} &\bigl|F(u(x+ry))-F(u(x))\bigr| \,dy\biggr|^2 \frac{dr}{r^{1+2\sigma}} \\
&\lesssim \int_{A(x)}^\infty \biggl|\int_{|y|<1} \bigl|u(x+ry)|^\alpha+|u(x)\bigr|^\alpha \,dy\biggr|^2 \frac{dr}{r^{1+2\sigma}} \\
&\lesssim \int_{A(x)}^\infty \frac{dr}{r^{1+2\sigma}} \bigl[ M\bigl(|u|^\alpha \bigr)(x)\bigr]^2  \\
&\lesssim \bigl[A(x)\bigr]^{-2\sigma} \bigl[ M\bigl(|u|^\alpha \bigr)(x)\bigr]^2.
\end{align*}

Choosing $A(x)=[M(|u|^\alpha)(x)]^{\frac1{s\alpha}} [\mathcal{D}_s(u)(x)]^{-\frac{1}{s}}$ leads immediately to \eqref{E:pointwise}.
The proposition follows from H\"older's inequality and boundedness of the maximal operator; the latter requires $(1-\frac\sigma{\alpha s})p_1>1$.
\end{proof}

%
%
%
%

\section{}

In this appendix we prove the existence and uniqueness of local $\dot{S}^0\bigcap\dot{S}^1$ solutions to the initial value problem
\begin{align}\label{eq for g}
\begin{cases}
(i\partial_t+\Delta)g=G + P_{hi}F(u_{lo}) + P_{hi}\bigl(F(u_{lo}+g)-F(g)-F(u_{lo})\bigr)\\
g(t_0)=u_{hi}(t_0),
\end{cases}
\end{align}
where the function $F$ represents the energy-critical nonlinearity, $u_{hi}=P_{>1}u$ and $u_{lo}=P_{\leq 1}u$
are as in subsection 5.2, and $G=G_{med}+G_{vhi}$ with
\begin{align*}
G_{med}&=\sum_{1<N<\eta_2^{-100}}\tilde{P_N} \bigl(\chi_{\{|P_NF(g)|\leq 1/N\}} P_NF(g)\bigr),\\
G_{vhi}&=|\nabla|^{-1}P_{\geq\eta_2^{-100}}\bigl(\chi_{\{||\nabla| F(g)|\leq 1\}}|\nabla| F(g) \bigr),
\end{align*}
and, in each case, $\chi$ represents a smooth cutoff to the set indicated.

In other words, we need to prove that the integral equation
\begin{align*}
g(t)&={}e^{i(t-t_0)\Delta}u_{hi}(t_0) + \int_{t_0}^t e^{i(t-s)\Delta} G(s)ds\\
   &\quad +\int_{t_0}^t e^{i(t-s)\Delta}P_{hi}\bigl(F(u_{lo}+g)-F(g)\bigr)(s)ds,
\end{align*}
admits an $\dot{S}^0\bigcap\dot{S}^1$ solution on a small interval $I:=[t_0, T]\subset I_0$, task which we accomplish
by proving convergence (in the appropriate spaces) of the iterates
\begin{align*}
g^{(1)}(t):=e^{i(t-t_0)\Delta}u_{hi}(t_0)
\end{align*}
and, for $m\geq 1$,
\begin{align}\label{recurrence}
g^{(m+1)}(t)&:={}e^{i(t-t_0)\Delta}u_{hi}(t_0) + \int_{t_0}^t e^{i(t-s)\Delta} G^{(m)}(s)ds \notag\\
    &\quad +\int_{t_0}^t e^{i(t-s)\Delta}P_{hi}\bigl(F(u_{lo}+g^{(m)})-F(g^{(m)})\bigr)(s)ds,
\end{align}
where
\begin{align*}
G^{(m)}&=G_{med}^{(m)}+G_{vhi}^{(m)}\\
&=\sum_{1<N<\eta_2^{-100}}\tilde{P_N} \bigl(\chi_{\{|P_NF(g^{(m)})|\leq 1/N\}} P_NF(g^{(m)})\bigr)\\
&\quad +|\nabla|^{-1}P_{\geq\eta_2^{-100}}\bigl(\chi_{\{||\nabla| F(g^{(m)})|\leq 1\}}|\nabla| F(g^{(m)}) \bigr).
\end{align*}

By Lemma~\ref{GB} (specifically \eqref{F(g)}, \eqref{nabla F(g)}, \eqref{Gvhi N0}, \eqref{Gvhi N1}, \eqref{Gmed N0},
and \eqref{Gmed N1}), we have
\begin{align}
\|G_{med}^{(m)}\|_{L_{t,x}^{\frac{2(n+2)}{n+4}}(\ir)}
&\lesssim \|g^{(m)}\|_{\dot{S}^0(\ir)}^{\frac{n+4}{n}}\label{A Gmed m}\\
\|\nabla G_{med}^{(m)}\|_{L_{t,x}^{\frac{2(n+2)}{n+4}}(\ir)}
&\lesssim \log(\tfrac{1}{\eta_2}) \bigl(\|g^{(m)}\|_{\dot{S}^1(\ir)}\|g^{(m)}\|_{\dot{S}^0(\ir)}^{\frac{4}{n-2}}\bigr)^{\frac{(n-2)(n+4)}{n(n+2)}}\label{A deriv Gmed m}\\
\|G_{vhi}^{(m)}\|_{L_{t,x}^{\frac{2(n+2)}{n+4}}(\ir)}
&\lesssim \eta_2^{100} \bigl(\|g^{(m)}\|_{\dot{S}^1(\ir)}\|g^{(m)}\|_{\dot{S}^0(\ir)}^{\frac{4}{n-2}}\bigr)^{\frac{(n-2)(n+4)}{n(n+2)}}\label{A Gvhi m}\\
\|\nabla G_{vhi}^{(m)}\|_{L_{t,x}^{\frac{2(n+2)}{n+4}}(\ir)}
&\lesssim  \bigl(\|g^{(m)}\|_{\dot{S}^1(\ir)}\|g^{(m)}\|_{\dot{S}^0(\ir)}^{\frac{4}{n-2}}\bigr)^{\frac{(n-2)(n+4)}{n(n+2)}}.\label{A deriv Gvhi m}
\end{align}

To simplify notation, we introduce the norm $\|\cdot\|_W$ defined on the slab $I\times\R^n$ as
$$
\|f\|_W=\|f\|_{W(I\times\R^n)}:=\|\nabla f\|_{L_t^{\frac{2(n+2)}{n-2}}L_x^{\frac{2n(n+2)}{n^2+4}}(I\times\R^n)}.
$$
Note that
$$
\|f\|_W\lesssim \|f\|_{\dot{S}^1(I\times\R^n)}.
$$
The first step is to take $I$ sufficiently small such that
\begin{align}\label{uhi W}
\|u_{hi}\|_{W} \lesssim \eta_2^{\eps}
\end{align}
for some $0<\eps\ll 1$; the dominated convergence theorem shows that this is indeed possible.

Throughout the rest of the proof all spacetime norms will be on $\ir$.
Using Bernstein, \eqref{diff2}, \eqref{diff3}, and our hypotheses on $u_{lo}$ (more precisely, \eqref{lowfreqsmall}
and \eqref{uloass}), we estimate
\begin{align}\label{A F(ulo)}
\|P_{hi}F(u_{lo})\|_{2,\frac{2n}{n+2}}
&\lesssim \|\nabla P_{hi}F(u_{lo})\|_{2,\frac{2n}{n+2}}
\lesssim \|\nabla u_{lo}\|_{2,\frac{2n}{n-2}}\|u_{lo}\|_{\infty,\frac{2n}{n-2}}^{\frac{4}{n-2}}\notag\\
&\lesssim \eta_2^{\frac{4(n-1)}{(n-2)^2}}
\end{align}
and
\begin{align}\label{A diff ulo and g}
\|P_{hi}\bigl(F(u_{lo}&+g^{(m)})-F(u_{lo})-F(g^{(m)})\bigr)\|_{2,\frac{2n}{n+2}}\notag\\
&\lesssim \|u_{lo}|g^{(m)}|^{\frac{4}{n-2}}\chi_{\{|u_{lo}|\leq |g^{(m)}|\}}\|_{2,\frac{2n}{n+2}}
   +\|g^{(m)}|u_{lo}|^{\frac{4}{n-2}}\chi_{\{|g^{(m)}|< |u_{lo}|\}}\|_{2,\frac{2n}{n+2}}\notag\\
&\lesssim \|g^{(m)}|u_{lo}|^{\frac{4}{n-2}}\|_{2,\frac{2n}{n+2}}
 \lesssim \|g^{(m)}\|_{2,\frac{2n}{n-2}}\|u_{lo}\|_{\infty,\frac{2n}{n-2}}^{\frac{4}{n-2}}\notag\\
&\lesssim \eta_2^{\frac{4}{n-2}} \|g^{(m)}\|_{\dot{S}^0},
\end{align}
\begin{align}\label{A deriv diff ulo and g}
\|\nabla P_{hi}\bigl(F(u_{lo}&+g^{(m)})-F(u_{lo})-F(g^{(m)})\bigr)\|_{2,\frac{2n}{n+2}}\notag\\
&\lesssim \|\nabla u_{lo}|g^{(m)}|^{\frac{4}{n-2}}\|_{2,\frac{2n}{n+2}}
   +\|\nabla g^{(m)}|u_{lo}|^{\frac{4}{n-2}}\|_{2,\frac{2n}{n+2}}\notag\\
&\lesssim \|\nabla u_{lo}\|_{2,\frac{2n}{n-2}}\|g^{(m)}\|_{\infty,\frac{2n}{n-2}}^{\frac{4}{n-2}}
   +\|\nabla g^{(m)}\|_{2,\frac{2n}{n-2}}\|u_{lo}\|_{\infty,\frac{2n}{n-2}}^{\frac{4}{n-2}}\notag\\
&\lesssim \eta_2^{\frac{4}{(n-2)^2}}\|g^{(m)}\|_{\dot{S}^1}^{\frac{4}{n-2}}
   +\eta_2^{\frac{4}{n-2}}\|g^{(m)}\|_{\dot{S}^1}.
\end{align}

Using the recurrence relation \eqref{recurrence}, Strichartz, and \eqref{A Gmed m} through
\eqref{A deriv diff ulo and g}, we estimate
\begin{align*}
\|g^{(m+1)}\|_{\dot{S}^0}
&\lesssim \|u_{hi}(t_0)\|_2+\|G^{(m)}\|_{\frac{2(n+2)}{n+4},\frac{2(n+2)}{n+4}} + \|P_{hi}F(u_{lo})\|_{2,\frac{2n}{n+2}}\\
&\quad +\|P_{hi}\bigl(F(u_{lo}+g^{(m)})-F(u_{lo})-F(g^{(m)})\bigr)\|_{2,\frac{2n}{n+2}}\\
&\lesssim \eta_2 + \|g^{(m)}\|_{\dot{S}^0(\ir)}^{\frac{n+4}{n}} +\eta_2^{100} \bigl(\|g^{(m)}\|_{\dot{S}^1(\ir)}\|g^{(m)}\|_{\dot{S}^0(\ir)}^{\frac{4}{n-2}}\bigr)^{\frac{(n-2)(n+4)}{n(n+2)}}\\
&\quad +\eta_2^{\frac{4(n-1)}{(n-2)^2}} +\eta_2^{\frac{4}{n-2}} \|g^{(m)}\|_{\dot{S}^0},
\end{align*}
\begin{align*}
\|g^{(m+1)}\|_{\dot{S}^1}
&\lesssim \|\nabla u_{hi}(t_0)\|_2+\|\nabla G^{(m)}\|_{\frac{2(n+2)}{n+4},\frac{2(n+2)}{n+4}} + \|\nabla P_{hi}F(u_{lo})\|_{2,\frac{2n}{n+2}}\\
&\quad +\|\nabla P_{hi}\bigl(F(u_{lo}+g^{(m)})-F(u_{lo})-F(g^{(m)})\bigr)\|_{2,\frac{2n}{n+2}}\\
&\lesssim 1 + \log(\tfrac{1}{\eta_2}) \bigl(\|g^{(m)}\|_{\dot{S}^1(\ir)}\|g^{(m)}\|_{\dot{S}^0(\ir)}^{\frac{4}{n-2}}\bigr)^{\frac{(n-2)(n+4)}{n(n+2)}}\\
&\quad +\bigl(\|g^{(m)}\|_{\dot{S}^1(\ir)}\|g^{(m)}\|_{\dot{S}^0(\ir)}^{\frac{4}{n-2}}\bigr)^{\frac{(n-2)(n+4)}{n(n+2)}} +\eta_2^{\frac{4(n-1)}{(n-2)^2}}\\
&\quad + \eta_2^{\frac{4}{(n-2)^2}}\|g^{(m)}\|_{\dot{S}^1}^{\frac{4}{n-2}} +\eta_2^{\frac{4}{n-2}}\|g^{(m)}\|_{\dot{S}^1},
\end{align*}
and
\begin{align*}
\|g^{(m+1)}\|_W
&\lesssim \|u_{hi}\|_W+\|\nabla G^{(m)}\|_{\frac{2(n+2)}{n+4},\frac{2(n+2)}{n+4}} + \|\nabla P_{hi}F(u_{lo})\|_{2,\frac{2n}{n+2}}\\
&\quad +\|\nabla P_{hi}\bigl(F(u_{lo}+g^{(m)})-F(u_{lo})-F(g^{(m)})\bigr)\|_{2,\frac{2n}{n+2}}\\
&\lesssim \eta_2^\eps + \log(\tfrac{1}{\eta_2}) \bigl(\|g^{(m)}\|_{\dot{S}^1(\ir)}\|g^{(m)}\|_{\dot{S}^0(\ir)}^{\frac{4}{n-2}}\bigr)^{\frac{(n-2)(n+4)}{n(n+2)}}\\
&\quad +\bigl(\|g^{(m)}\|_{\dot{S}^1(\ir)}\|g^{(m)}\|_{\dot{S}^0(\ir)}^{\frac{4}{n-2}}\bigr)^{\frac{(n-2)(n+4)}{n(n+2)}} +\eta_2^{\frac{4(n-1)}{(n-2)^2}}\\
&\quad + \eta_2^{\frac{4}{(n-2)^2}}\|g^{(m)}\|_{\dot{S}^1}^{\frac{4}{n-2}} +\eta_2^{\frac{4}{n-2}}\|g^{(m)}\|_{\dot{S}^1}.
\end{align*}
A simple inductive argument yields
\begin{align}
\|g^{(m)}\|_{\dot{S}^0}&\lesssim \eta_2^{\frac{2}{n-2}}, \label{A gm S0}\\
\|g^{(m)}\|_{\dot{S}^1}&\lesssim 1, \label{A gm S1}\\
\|g^{(m)}\|_W &\lesssim \eta_2^\eps, \label{A gm W}
\end{align}
for all $m\geq 1$ and provided $\eps$ is sufficiently small.

To prove that the initial value problem \eqref{eq for g} admits a local solution in $\dot{S}^0 \bigcap \dot{S}^1$,
it suffices to prove that the sequence $\{g^{(m)}\}_m$ converges strongly in $\dot{S}^0$ to some function $g$.
This function is guaranteed to lie in $\dot{S}^1$ for the following reasons: As the sequence $\{g^{(m)}\}_m$ stays
bounded in $\dot{S}^1$ (see \eqref{A gm S1}), it follows that $g^{(m)}$ converges weakly to $g$ in $\dot{S}^1$.
As weak limits are unique, we conclude that $g$ lies in both $\dot{S}^0$ and $\dot{S}^1$. Moreover, by Fatou,
\eqref{A gm S0}, and \eqref{A gm S1}, $g$ obeys the bounds
$$
\|g\|_{\dot{S}^0}\lesssim \eta_2^{\frac{2}{n-2}} \quad \text{and} \quad \|g\|_{\dot{S}^1}\lesssim 1.
$$

In what follows, we prove that the sequence $\{g^{(m)}\}_m$ is Cauchy in $\dot{S}^0$, which completes the proof of local
existence for the reasons just given. We start by considering differences of the form
$g^{(m+1)}-g^{(m)}$. Using the recurrence relation \eqref{recurrence} and Strichartz, for $m>1$ we bound
\begin{align}\label{ultima diff}
\|g^{(m+1)}-g^{(m)}\|_{\dot{S}^0}
&\lesssim \|G^{(m)}-G^{(m-1)}\|_{L_{t,x}^\frac{2(n+2)}{n+4}} + \|F(g^{(m)})-F(g^{(m-1)})\|_{L_{t,x}^\frac{2(n+2)}{n+4}}\notag\\
&\quad +\|F(u_{lo}+g^{(m)})-F(u_{lo}+g^{(m-1)})\|_{L_{t,x}^\frac{2(n+2)}{n+4}}.
\end{align}
By \eqref{diff1}, Sobolev embedding, \eqref{uloass}, and \eqref{A gm W}, we get
\begin{align}\label{A F(gm) OK}
\|F(g^{(m)})-F(&g^{(m-1)})\|_{L_{t,x}^\frac{2(n+2)}{n+4}}\notag\\
&\lesssim \|g^{(m)}-g^{(m-1)}\|_{L_{t,x}^\frac{2(n+2)}{n}}
    \Bigl(\|g^{(m)}\|_{L_{t,x}^\frac{2(n+2)}{n-2}}^{\frac{4}{n-2}}+\|g^{(m-1)}\|_{L_{t,x}^\frac{2(n+2)}{n-2}}^{\frac{4}{n-2}}\Bigr)\notag\\
&\lesssim \|g^{(m)}-g^{(m-1)}\|_{\dot{S}^0}\bigl(\|g^{(m)}\|_W^{\frac{4}{n-2}}+\|g^{(m-1)}\|_W^{\frac{4}{n-2}}\bigr)\notag\\
&\lesssim \eta_2^{\frac{4\eps}{n-2}}\|g^{(m)}-g^{(m-1)}\|_{\dot{S}^0}
\end{align}
and
\begin{align}\label{A F(ulo+gm) OK}
\|F(u_{lo}+g^{(m)})&-F(u_{lo}+g^{(m-1)})\|_{L_{t,x}^\frac{2(n+2)}{n+4}}\notag\\
&\lesssim \Bigl(\|g^{(m)}\|_{L_{t,x}^\frac{2(n+2)}{n-2}}+\|g^{(m-1)}\|_{L_{t,x}^\frac{2(n+2)}{n-2}}^{\frac{4}{n-2}}
   +\|u_{lo}\|_{L_{t,x}^\frac{2(n+2)}{n-2}}^{\frac{4}{n-2}}\Bigr)\notag\\
&\qquad \qquad \times \|g^{(m)}-g^{(m-1)}\|_{\frac{2(n+2)}{n},\frac{2(n+2)}{n}}\notag\\
&\lesssim \|g^{(m)}-g^{(m-1)}\|_{\dot{S}^0}\Bigl(\|g^{(m)}\|_W^{\frac{4}{n-2}}+\|g^{(m-1)}\|_W^{\frac{4}{n-2}}+\|u_{lo}\|_{\dot{S}^1}^{\frac{4}{n-2}}\Bigr)\notag\\
&\lesssim \eta_2^{\frac{4\eps}{n-2}}\|g^{(m)}-g^{(m-1)}\|_{\dot{S}^0},
\end{align}
again, assuming $\eps$ is sufficiently small, which amounts to taking $I$ sufficiently small.

We are left with the task of estimating
\begin{align*}
\|G^{(m)}-G^{(m-1)}\|_{L_{t,x}^\frac{2(n+2)}{n+4}}
\leq \|G_{med}^{(m)}-G_{med}^{(m-1)}\|_{L_{t,x}^\frac{2(n+2)}{n+4}} + \|G_{vhi}^{(m)}-G_{vhi}^{(m-1)}\|_{L_{t,x}^\frac{2(n+2)}{n+4}}.
\end{align*}

We consider first differences coming from medium frequencies. For $1<N<\eta_2^{-100}$, denote
$$
f^N:=\chi_{\{|P_NF|\leq 1/N\}} P_NF.
$$
Then,
\begin{align*}
\|G_{med}^{(m)}-G_{med}^{(m-1)}\|_{L_{t,x}^\frac{2(n+2)}{n+4}}
&\lesssim \sum_{1<N<\eta_2^{-100}} \|\tilde{P_N}\bigl(f^N(g^{(m)})-f^N(g^{(m-1)})\bigr)\|_{L_{t,x}^\frac{2(n+2)}{n+4}}\\
&\lesssim \sum_{1<N<\eta_2^{-100}} \|f^N(g^{(m)})-f^N(g^{(m-1)})\|_{L_{t,x}^\frac{2(n+2)}{n+4}}.
\end{align*}
By the Fundamental Theorem of Calculus, we write
\begin{align*}
f^N(g^{(m)})-f^N(g^{(m-1)})= (g^{(m)}-g^{(m-1)})\int_0^1 (f^N_z+f^N_{\bar z})(g^{(m)}_\theta)d \theta,
\end{align*}
where for $\theta \in [0,1]$ we define $g^{(m)}_\theta := g^{(m)} + \theta( g^{(m)}-g^{(m-1)})$.
As
$$
f^N_z=\chi_{\{|P_NF|\leq 1/N\}} P_NF_z + \tilde{\chi}_{\{|P_NF|\sim 1/N\}} N P_NF
$$
and
$$
f^N_{\bar z}=\chi_{\{|P_NF|\leq 1/N\}} P_NF_{\bar z} + \tilde{\chi}_{\{|P_NF|\sim 1/N\}} N P_NF
$$
for some smooth cutoff $\tilde{\chi}$ to the set indicated, Minkowski's inequality and Sobolev embedding yield
\begin{align*}
\|f^N(g^{(m)})&-f^N(g^{(m-1)})\|_{L_{t,x}^\frac{2(n+2)}{n+4}}\\
&\lesssim \sup_{\theta \in [0,1]} \bigl(\|g^{(m)}_\theta \|_{L_{t,x}^{\frac{2(n+2)}{n-2}}}^{\frac{4}{n-2}}+\|\tilde{\chi}_{\{|P_NF(g^{(m)}_\theta)|\sim 1/N\}} N P_NF( g^{(m)}_\theta)\|_{L_{t,x}^{\frac{n+2}{2}}}\bigr)\\
&\qquad \qquad \qquad \times \|g^{(m)}-g^{(m-1)}\|_{L_{t,x}^{\frac{2(n+2)}{n}}}\\
&\lesssim \bigl(\|g^{(m)}\|_W^{\frac{4}{n-2}}+\|g^{(m-1)}\|_W^{\frac{4}{n-2}}\bigr)\|g^{(m)}-g^{(m-1)}\|_{\dot{S}^0}\\
&\quad +\sup_{\theta \in [0,1]}\|\tilde{\chi}_{\{|P_NF(g^{(m)}_\theta)|\sim 1/N\}} N P_NF( g^{(m)}_\theta)\|_{L_{t,x}^{\frac{n+2}{2}}}\|g^{(m)}-g^{(m-1)}\|_{\dot{S}^0}\\
&\lesssim \eta_2^{\frac{4\eps}{n-2}}\|g^{(m)}-g^{(m-1)}\|_{\dot{S}^0}\\
&\quad +\sup_{\theta \in [0,1]}\|\tilde{\chi}_{\{|P_NF(g^{(m)}_\theta)|\sim 1/N\}} N P_NF( g^{(m)}_\theta)\|_{L_{t,x}^{\frac{n+2}{2}}}\|g^{(m)}-g^{(m-1)}\|_{\dot{S}^0}.
\end{align*}
For $\theta \in [0,1]$, by Bernstein and \eqref{nabla F(g)} we estimate
\begin{align*}
\|\tilde{\chi}_{\{|P_NF(g^{(m)}_\theta)|\sim 1/N\}} &N P_NF( g^{(m)}_\theta)\|_{L_{t,x}^{\frac{n+2}{2}}}\\
&\sim \bigl\|\tilde{\chi}_{\{|P_NF(g^{(m)}_\theta)|\sim 1/N\}} \bigl[N P_NF( g^{(m)}_\theta)\bigr]^{\frac{4(n-2)}{n(n+2)}}\bigr\|_{L_{t,x}^{\frac{n+2}{2}}}\\
&\lesssim \|N P_NF( g^{(m)}_\theta)\|_{L_{t,x}^{\frac{2(n-2)}{n}}}^{\frac{4(n-2)}{n(n+2)}}\\
&\lesssim \|\nabla F( g^{(m)}_\theta)\|_{L_{t,x}^{\frac{2(n-2)}{n}}}^{\frac{4(n-2)}{n(n+2)}}\\
&\lesssim \bigl(\|g^{(m)}_\theta\|_{\dot{S}^1} \|g^{(m)}_\theta\|_{\dot{S}^0}^{\frac{4}{n-2}}\bigr)^{\frac{4(n-2)}{n(n+2)}}.
\end{align*}
Thus, by \eqref{A gm S0}, \eqref{A gm S1}, and taking $\eps$ sufficiently small, we get
\begin{align*}
\|f^N(g^{(m)})&-f^N(g^{(m-1)})\|_{L_{t,x}^\frac{2(n+2)}{n+4}}
\lesssim \eta_2^{\frac{4\eps}{n-2}}\|g^{(m)}-g^{(m-1)}\|_{\dot{S}^0},
\end{align*}
which yields
\begin{align}\label{A med diff OK}
\|G_{med}^{(m)}-G_{med}^{(m-1)}\|_{L_{t,x}^\frac{2(n+2)}{n+4}}
&\lesssim \sum_{1<N<\eta_2^{-100}} \eta_2^{\frac{4\eps}{n-2}}\|g^{(m)}-g^{(m-1)}\|_{\dot{S}^0} \notag \\
&\lesssim \log(\tfrac{1}{\eta_2})\eta_2^{\frac{4\eps}{n-2}}\|g^{(m)}-g^{(m-1)}\|_{\dot{S}^0} \notag \\
&\lesssim \eta_2^{\frac{3\eps}{n-2}}\|g^{(m)}-g^{(m-1)}\|_{\dot{S}^0}.
\end{align}

We consider last the difference coming from very high frequencies. More precisely, denoting
$$
f^{vhi}:=|\nabla|^{-1}P_{\geq \eta_2^{-100}}\bigl(\chi_{\{||\nabla|F|\leq 1\}}|\nabla|F\bigr),
$$
we have to estimate
\begin{align*}
\|G_{vhi}^{(m)}-G_{vhi}^{(m-1)}\|_{L_{t,x}^\frac{2(n+2)}{n+4}}
= \|f^{vhi}(g^{(m)})-f^{vhi}(g^{(m-1)})\|_{L_{t,x}^\frac{2(n+2)}{n+4}}.
\end{align*}
By the Fundamental Theorem of Calculus, we write
$$
f^{vhi}(g^{(m)})-f^{vhi}(g^{(m-1)})=(g^{(m)}-g^{(m-1)})\int_0^1 \nabla f^{vhi} (g^{(m)}_\theta)d \theta,
$$
with the convention that for $\theta \in [0,1]$, $g^{(m)}_\theta := g^{(m)} + \theta( g^{(m)}-g^{(m-1)})$.
By Minkowski's inequality, the boundedness of the Riesz potentials on $L_x^p$ for $1<p<\infty$, and \eqref{nabla F(g)},
we estimate
\begin{align*}
\|f^{vhi}(g^{(m)})&-f^{vhi}(g^{(m-1)})\|_{L_{t,x}^\frac{2(n+2)}{n+4}}\\
&\lesssim \|g^{(m)}-g^{(m-1)}\|_{L_{t,x}^{\frac{2(n+2)}{n}}} \sup_{\theta\in [0,1]}\|\chi_{\{||\nabla|F(g^{(m)}_\theta)|\leq 1\}}|\nabla|F(g^{(m)}_\theta)\|_{L_{t,x}^{\frac{n+2}{2}}}\\
&\lesssim \|g^{(m)}-g^{(m-1)}\|_{\dot{S}^0}\bigl\|\chi_{\{||\nabla|F(g^{(m)}_\theta)|\leq 1\}}\bigr[|\nabla|F(g^{(m)}_\theta)\bigr]^{\frac{4(n-2)}{n(n+2)}}\bigr\|_{L_{t,x}^{\frac{n+2}{2}}}\\
&\lesssim \|g^{(m)}-g^{(m-1)}\|_{\dot{S}^0} \||\nabla|F(g^{(m)}_\theta)\|_{L_{t,x}^{\frac{2(n-2)}{n}}}^{\frac{4(n-2)}{n(n+2)}}\\
&\lesssim \|g^{(m)}-g^{(m-1)}\|_{\dot{S}^0}\bigl(\|g^{(m)}_\theta\|_{\dot{S}^1} \|g^{(m)}_\theta\|_{\dot{S}^0}^{\frac{4}{n-2}}\bigr)^{\frac{4(n-2)}{n(n+2)}}.
\end{align*}
Thus, by \eqref{A gm S0}, \eqref{A gm S1}, and taking $\eps$ sufficiently small, we get
\begin{align}\label{A vhi diff OK}
\|G_{vhi}^{(m)}-G_{vhi}^{(m-1)}\|_{L_{t,x}^\frac{2(n+2)}{n+4}}
\lesssim \eta_2^{\frac{4\eps}{n-2}} \|g^{(m)}-g^{(m-1)}\|_{\dot{S}^0}.
\end{align}

Collecting \eqref{ultima diff}, \eqref{A F(gm) OK}, \eqref{A F(ulo+gm) OK}, \eqref{A med diff OK},
and \eqref{A vhi diff OK}, we obtain
$$
\|g^{(m+1)}-g^{(m)}\|_{\dot{S}^0}\lesssim \eta_2^{\frac{3\eps}{n-2}} \|g^{(m)}-g^{(m-1)}\|_{\dot{S}^0}
$$
for all $m>1$. Also, as by \eqref{recurrence}, \eqref{A Gmed m}, \eqref{A F(ulo)}, and \eqref{A diff ulo and g}, we have
\begin{align*}
\|g^{(2)}-g^{(1)}\|_{\dot{S}^0}
&\lesssim \|G^{(1)}\|_{\frac{2(n+2)}{n+4},\frac{2(n+2)}{n+4}} + \|P_{hi}F(u_{lo})\|_{2,\frac{2n}{n+2}}\\
&\quad +\|P_{hi}\bigl(F(u_{lo}+g^{(1)})-F(u_{lo})-F(g^{(1)})\bigr)\|_{2,\frac{2n}{n+2}}\\
&\lesssim  \|g^{(1)}\|_{\dot{S}^0(\ir)}^{\frac{n+4}{n}} +\eta_2^{100} \bigl(\|g^{(1)}\|_{\dot{S}^1(\ir)}\|g^{(1)}\|_{\dot{S}^0(\ir)}^{\frac{4}{n-2}}\bigr)^{\frac{(n-2)(n+4)}{n(n+2)}}\\
&\quad +\eta_2^{\frac{4(n-1)}{(n-2)^2}} +\eta_2^{\frac{4}{n-2}} \|g^{(1)}\|_{\dot{S}^0}\\
&\lesssim \eta_2^{\frac{n+4}{n}}+\eta_2^{100}\eta_2^{\frac{4(n+4)}{n(n+2)}}+\eta_2^{\frac{4(n-1)}{(n-2)^2}} +\eta_2^{\frac{n+2}{n-2}},
\end{align*}
we immediately obtain that the sequence $\{g^{(m)}\}_{m \in \N}$ is Cauchy in $\dot{S}^0$ and thus convergent to some
function $g\in \dot{S}^0$.

The uniqueness of $\dot{S}^0\cap\dot{S}^1$ solutions to \eqref{eq for g} is standard and based on the estimates above.
We skip the details.

%
%
%
%

\end{document}